# EXPANSION IN FOURIER PERIODIC SERIES, METHOD FOR SOLVING DIFFERENTIAL EQUATIONS

Arpad Török[1], Stoian Petrescu[2], Michel Feidt[3]

[1] PhD student, The Polytechnic University of Bucharest, Department of Engineering Thermodynamics, 313, Splaiul Independentei, 060042 Bucharest, Romania, e-mail: arpi_torok@yahoo.com
[2] Prof. Dr. Eng., Polytechnic University of Bucharest, Department of Engineering Thermodynamics, România
[3] Prof. Dr. Eng., L.E.M.T.A., U.R.A. C.N.R.S. 7563, Université de Loraine Nancy 12, avenue de la Foret de Haye, 54516 Vandœuvre-lès-Nancy, France

**Summary:** The expansion of functions of real variables in Taylor and Frobenius series (power series which are nonorthogonal, nonperiodic bases), in sinusoidal Fourier series (bases of orthogonal, periodic functions), in series of special functions (bases of orthogonal, nonperiodic functions) etc, is a commonly used method for solving a wide range of ordinary differential equations (ODEs) and partial differential equations (PDEs).

In this paper, based on an in-depth analysis of the properties of periodic sinusoidal Fourier series (SFS), we will be able to apply this procedure to a much broader category of ODEs (all linear, homogeneous and non-homogeneous equations with constant coefficients, a large category of linear and non-linear equations with variable coefficients, systems of ODEs, integro-differential equations, etc.). We will also extend this procedure and we use it to solve certain ODEs, on non-orthogonal periodic bases, represented by non sinusoidal periodic Fourier series (SFN).

## 1. Introduction

The ODE resolution methods proposed here use recently obtained results [1] in the field of functional analysis, concerning the expansion of real variable functions, defined over an interval [–L, L], in infinite series of periodic functions over the same interval, forming orthogonal, but also non-orthogonal bases. According to harmonic analysis (Fourier), any function $f(x)$, periodic over the interval [−L, L], which satisfies the Dirichlet conditions, can be developed into an infinite sum, known in the literature under the name trigonometric series (for which, for reasons highlighted in article [1], we used the name sinusoidal series). This series is formed by the components of a complete biortagonal base, composed by the unit function $1$, the fundamental harmonics $sin(\pi x/L)$-even and $cos(\pi x/L)$-odd, with the period $2L$ and the secondary harmonics $sin(n\pi x/L)$ and $cos(n\pi x/L)$, with the period $2L/n$, for $n \in \mathbf{N}^+$. The coefficients of this expansion (Fourier coefficients) can be calculated using defined integrals (Euler formulas). This paper generalizes this statement by showing that the function $f(x)$ can also be developed in non-sinusoidal periodic series (SFN), consisting of the infinite sum of the weighted components of a *complete*, non-orthogonal base: the unit function $1$, the fundamental quasi-harmonics $g(x)$-even and $h(x)$-odd, periodic with the period $2L$, with mean value zero over the definition interval and the secondary quasi-harmonics, defined on [−L, L], $g_n(x)=g(nx)$ and $h_n(x)=h(nx)$, with the period $2L/n$, for $n \in \mathbf{N}^+$. The fundamental quasi-harmonics $g(x)$ and $h(x)$ can be any real function, of real variable, which admits on the interval [−L, L] expansions in sinusoidal series. We obtain the coefficients $A_n$ and $B_n$ of the expansion in SFN of the function $f(x)$, using the algebraic relationships between the Fourier coefficients of the expansions in SFS of the functions $f(x)$, $g(x)$ and $h(x)$.Thus, any function $f(x):[x_1, x_2] \rightarrow \mathbf{R}$, $T$-periodic ($T=x_2−x_1$), of $L^2$-space (that is, square-integrable), can be represented by the sum:

$\widehat{f}(x) = f_0 + \sum_{n=1}^{\infty} A_n [\overline{g}_n(x) - g_0] + \sum_{n=1}^{\infty} B_n \overline{h}_n(x)$ , where :

$f_0 = \dfrac{1}{T} \int_{x_1}^{x_2} f(x) dx$ and $g_0 = \dfrac{1}{T} \int_{x_1}^{x_2} g(x) dx$ are the mean values over the interval $[x_1, x_2]$,

$\overline{g}_n(x)$ et $\overline{h}_n(x)$ ($n=1, 2, 3, ..., \infty$) are **Fourier-functions** (continuous functions, precisely the expansions in SFS of the functions $g(nx)$ and $h(nx)$, defined on $[-T/n, T/n]$.

From this general result, in the paper [1] are also described some particular situations:

- $\widehat{f}(x) = f_0 + \sum_{n=-\infty}^{\infty} C_n [\overline{g}_n(x) - g_0]$ , where $g(x)$ has both components: even and odd

- $\widehat{f}(x) = f_0 + \sum_{n=1}^{\infty} A_n \mathbf{S}[g(x)]_L + \sum_{n=1}^{\infty} B_n \mathbf{C}[g(x+L/2)]_L$ where the function $g(x)$ is a function of $L^2$ defined on the interval $[0, L/2]$, $\mathbf{S}[g(x)]_L$ and $\mathbf{C}[g(x+L/2)]_L$ are the functions called quasi-sinusoids, derived from it

- $\widetilde{f}(x) = f_0 + \sum_{n=1}^{\infty} [A_n^0 \Phi_n(x) + B_n^0 \Psi_n(x)]$ where $\Phi_n(x)$ and $\Psi_n(x)$ are orthogonal functions, generated by the Fourier-functions $\overline{g}(x) - g_0$ and $\overline{h}(x)$ by an orthogonalisation process.

## 1.1. Properties of expansions in periodic Fourier series (the general case)

The expansion in SFN of a function $f(x)$ is obtained from its expansion in SFS, by a redistribution and a grouping of the terms of expansion, so as to obtain the expansions in SFS of the periodic functions $g(x)$, $h(x)$ and corresponding quasi-harmonics $g_n(x)$, $h_n(x)$. Consequently, this expansion also benefits from the properties of convergence, differentiability and integrability, similar to those of SFS [1-9]. So:

- let $f(x)$ be a *2L*-periodic function, continuous in the interval $[-L, L]$. Its Fourier expansion $f(x) = f_0 + \sum_{n=1}^{\infty} A_n g_n(x) + \sum_{n=1}^{\infty} B_n h_n(x)$, where $g_0=0$, sinusoidal or not, convergent or not, can be integrated term by term, between all integration limits:

$$F(x) = \int_0^x f(x) dx = C + f_0 x + \sum_{n=1}^{\infty} A_n G_n(x) dx + \sum_{n=1}^{\infty} B_n H_n(x) dx, \qquad (1.1)$$

where $G(x)$ and $H(x)$ are the primitives of $g(x)$, respectively $h(x)$, and $C$ is an integration constant which depends on the expansion coefficients. After the replacements $x = \sum_{n=1}^{\infty} C_n G_n(x)$

(the functions $x$ and $G(x)$ are odd) and $C = \sum_{n=1}^{\infty} B_n H_n'(0) + F(0)$, results an expansion in SFN of the primitive $F(x)$, in a base, most often different from that of the function $f(x)$.

- let $f(x)$ be a *2L*-periodic function, continuous in the interval $[-L, L]$, with $f(-L)=f(L)$ and with the derivative $f'(x)$ smooth by pieces in this interval. The Fourier expansion, sinusoidal or not, of the function $f'(x)$, can be obtained by deriving term by term the Fourier expansion of the function $f(x)$. The series obtained converges punctually towards $f'(x)$ in all the points of continuity and towards $[f'(x)+f'(-x)]/2$ in those of discontinuity.

If $f(x) = f_0 + \sum_{n=1}^{\infty} A_n g_n(x) + \sum_{n=1}^{\infty} B_n h_n(x)$ , so: $\widehat{f}'(x) = \sum_{n=1}^{\infty} A_n g_n'(x) + \sum_{n=1}^{\infty} B_n h_n'(x)$ \qquad (1.2)

In this case too, the basis of the non-sinusoidal expansion of the derivative differs from that of the function $f(x)$.

The condition $f(-L)=f(L)$ means that the number of problems in which formula (1.2) can be useful is quite small, but it can be avoided if the jump from point $x=L$ (as well as any other jump of the odd component) is compensated by a jump in the opposite direction (by subtracting another odd function which makes an identical jump in the same point). In the following example, by this process, the odd component $f_o$ of the function $f(x)$ is decomposed into a sum of the differentiable function $f_{os}$ (for which $f(-L)=f(L)$) and the odd ramp function $f_r=x\cdot f_o(L)/L \rightarrow f_{os}=f_o-f_r$. Therefore:

$$\frac{d}{dx}f_o(x) = \frac{d}{dx}\left[f_{os}(x) + \frac{f_o(L)}{L}x\right] = \frac{d}{dx}f_{os}(x) + \frac{f_o(L)}{L},$$

relation which allows us to find an expression for the expansion of the derivative $f'(x)$ of the function $f(x)$, for all the categories of the functions which satisfy the other conditions. The expansion in SFN of the function $f_{os}=f_o-f_r$ allows us to calculate the coefficients of the expansion in series of the derivative.

The general case, that of the expansion of functions in non-sinusoidal periodic bases, highlights the fact that the element $I=1$ of the base has a particular character. It is part of all periodic bases, it does not change while the other components of the base (with zero mean value over the definition interval) change after integration, or after derivation. Its coefficient is calculated by a definite integral and not by algebraic relations. $I=1$ is an even function, but for $f_0=0$, it simultaneously receives an odd character too. The derivatives and primitives of all the even functions (including the function $I\cdot f_0=f_0$) are odd functions and, conversely, those of the odd functions (including the function $I\cdot 0=0$) are even functions. By deriving any even function $f_0$ we obtain the odd function $0$ and vice versa, by integrating the odd function $0$, we obtain any even function $f_0$ (the integration constant).

### 1.2. Properties of expansions in sinusoidal Fourier series

As in the general case, in the case of expansions in SFS, over the interval $[-L, L]$, for $\omega_n=n\omega_0=n\pi/L$, because $\bar{x}=2\sum_{n=1}^{\infty}\frac{(-1)^{n+1}}{\omega_n}\sin(\omega_n x)$, we can write [9]:

For $\bar{f}(x) = f_0 + \bar{f}_e + \bar{f}_o = f_0 + \sum_{n=1}^{\infty}a_n\cos(\omega_n x) + \sum_{n=1}^{\infty}b_n\sin(\omega_n x),$ \hfill (1.3)

we have: $f_o(L)=-f_o(-L)$, $f_e(L)=f_e(-L)$, $f_e(L)=\dfrac{f(L)+f(-L)}{2}$, $f_o(L)=\dfrac{f(L)-f(-L)}{2}$,

$f(0) = f_0 + \sum_{n=1}^{\infty}a_n$,

$\bar{f}'(x)=\overline{\Phi}(x)=\Phi_0+\overline{\Phi}_e+\overline{\Phi}_o=\Phi_0\cdot 1+\sum_{n=1}^{\infty}\alpha_n\cdot\cos(\omega_n x)+\sum_{n=1}^{\infty}\beta_n\cdot\sin(\omega_n x),$ in which: \hfill (1.4)

$\Phi_0=\dfrac{f_o(L)}{L}=\lim_{n\to\infty}2\left[(-1)^{n+1}\omega_n b_n\right]$, $\alpha_n=\omega_n\left(b_n+2(-1)^n\dfrac{f_o(L)}{\omega_n L}\right)$, $\beta_n=-a_n\omega_n$, $\Phi(0)=\Phi_0+\sum_{n=1}^{\infty}\alpha_n$

$\left(\int_0^x f(x)dx\right)^{SFS}=\overline{F}(x)-F(0)=F_{00}+\overline{F}_e+\overline{F}_o=F_{00}\cdot 1+\sum_{n=1}^{\infty}A_n\cdot\cos(\omega_n x)+\sum_{n=1}^{\infty}B_n\cdot\sin(\omega_n x),$ in which:

$F_{00}=\sum_{n=1}^{\infty}\dfrac{b_n}{\omega_n}$, $A_n=-\dfrac{b_n}{\omega_n}$, $B_n=\dfrac{a_n}{\omega_n}+\dfrac{2(-1)^{n+1}}{\omega_n}\dfrac{f_0(L)}{L}$, $F(0)=F_0+\sum_{n=1}^{\infty}A_n=F_0-F_{00}$ \hfill (1.5)

Also valid are the relationships: $\bar{f}(x)=\int_0^x\overline{\Phi}(x)+\dfrac{f_o(L)}{L}+\sum_{n=1}^{\infty}a_n$ , $f_0=\dfrac{F_o(L)}{L}$ and

$$\int_a^b \bar{f}(x)dx = f_0(b-a) + \sum_{n=1}^{\infty} \frac{a_n(\sin \omega_n b - \sin \omega_n a) - b_n(\cos \omega_n b - \cos \omega_n a)}{n\omega_0} \qquad (1.6)$$

Thus, we obtained a series of relations to calculate the mean values and the expansion coefficients of the derivative and primitive functions of the first rank, from the values of the expansion coefficients in SFS of the function $f(x)$ and the values of the function in the limit points $f(-L)$ and $f(L)$. These relationships also allow us to compute, step by step, the expressions of derivatives and higher-order primitives, after having calculated the values of these functions at the limit of the interval, the mean values and the values of the coefficients of their expansion. The new relations can be used to solve differential and integro-differential equations of higher rank, to calculate definite or indefinite integrals, etc. As an example, here is the derivation and integration of the expansion of the function $f(x)=e^x$ for the interval $[-\pi, \pi]$ [8, 9]:

$$\bar{f}(x) = e^x = \frac{\sinh \pi}{\pi} + \sum_{n=1}^{\infty}\left[\frac{2\sinh \pi}{\pi}\cdot\frac{(-1)^n}{1+n^2}\cos(nx) - \frac{2n\cdot\sinh \pi}{\pi}\cdot\frac{(-1)^n}{1+n^2}\sin(nx)\right].$$

$$\Phi_0 = \frac{\sinh \pi}{\pi}, \ \alpha_n = n\left(-\frac{\sinh \pi}{\pi}\cdot\frac{2n(-1)^n}{1+n^2} + 2(-1)^n\frac{\sinh \pi}{n\cdot\pi}\right) = \frac{2\sinh \pi}{\pi}\cdot\frac{(-1)^n}{1+n^2}, \ \beta_n = -n\frac{2\sinh \pi}{\pi}\cdot\frac{(-1)^n}{1+n^2},$$

$$F_{00} = \sum_{n=1}^{\infty}\frac{b_n}{\omega_n} = \frac{2n\cdot\sinh \pi}{n\cdot\pi}\sum_{n=1}^{\infty}\frac{(-1)^{n+1}}{1+n^2} = \frac{2\sinh \pi}{\pi}\cdot\frac{1}{2}\left(1-\frac{\pi}{\sinh \pi}\right) = \frac{\sinh \pi}{\pi}-1, \ F_0 = \frac{\sinh \pi}{\pi} \qquad (1.7)$$

$$A_n = \frac{-1}{n}\left(-\frac{2n\cdot\sinh \pi}{\pi}\cdot\frac{(-1)^n}{1+n^2}\right), \ B_n = \frac{1}{n}\left[\frac{2\sinh \pi}{\pi}\cdot\frac{(-1)^n}{1+n^2} - 2(-1)^n\frac{\sinh \pi}{\pi}\right] = -\frac{2n\cdot\sinh \pi}{\pi}\cdot\frac{(-1)^n}{1+n^2}$$

So: $\overline{\Phi}(x) = \left[\left(e^x\right)'\right]^{SFS} = e^x$ and $\overline{F}(x) = \left(\int_0^x e^x dx\right)^{SFS} + F(0) = e^x$

To get (1.7), we start from $\bar{f}(0) = 1 = \frac{\sinh \pi}{\pi}\left[1 + 2\sum_{n=1}^{\infty}\frac{(-1)^n}{1+n^2}\right] \rightarrow \sum_{n=1}^{\infty}\frac{(-1)^n}{1+n^2} = \frac{1}{2}\left(\frac{\pi}{\sinh \pi}-1\right)$

We can notice that for $f_0=0$ and $f_o(L)=0$ (without discontinuities of the odd component), both the integration and the derivation are made term by term: $\alpha_n = b_n\omega_n$, $\beta_n = -a_n\omega_n$, $A_n = b_n/\omega_n$, $B_n = -a_n/\omega_n$. If $f_o(x)$ has discontinuities at the ends of the interval $[-L, L]$ (or inside), they will cause, during the derivation, the appearance of an average value $\Phi_0\neq 0$ and a modification as a consequence of the component $\Phi_e$. During integration, the effect of discontinuities is transmitted to the odd component $F_o$. The derivatives of all the functions $f(x)+C$ have the same expression, and the return to the initial function, by integration, is ensured by the relation $F(0) = F_0 + \sum_{n=1}^{\infty} A_n$, which for $C=0$ leads to $F(0)=0$.

On the real axis, the odd periodic function $f(x)$, discontinuous at the ends of the interval $[-L, L]$ (Fig.1A), is the sum between the continuous function $f_C(x)$ (Fig.1D) and the "scale" function $f_H(x)$ (a succession of negative Heaviside "steps") (Fig.1E). The derivatives of these functions are the periodic functions $f'(x) = \Xi_{2\pi}(x+L) \stackrel{def}{=} \sum_{n=-\infty}^{\infty}\delta_{2n\pi}(x+L)$ (the Dirac comb, Fig.1B), whose expansion in SFS is $\overline{\Xi} = \frac{1}{\pi}\sum_{n=1}^{\infty}(-1)^n\cos(n\omega_0 x)$ [6], multiplied by the coefficient $f_o(L)/L$ (Fig.1C). The expansions in SFS of these two derivatives are therefore divergent, but their sum is convergent.

If the function $f(x)$ is even, the derivation of jumps within the interval generates two Dirac combs of opposite sign, which cancel each other, and if $f(x)$ has finite discontinuities in a finite number of points in the interval [- *L*, *L*], this is reflected in the position and amplitude of the corresponding Dirac pulses.

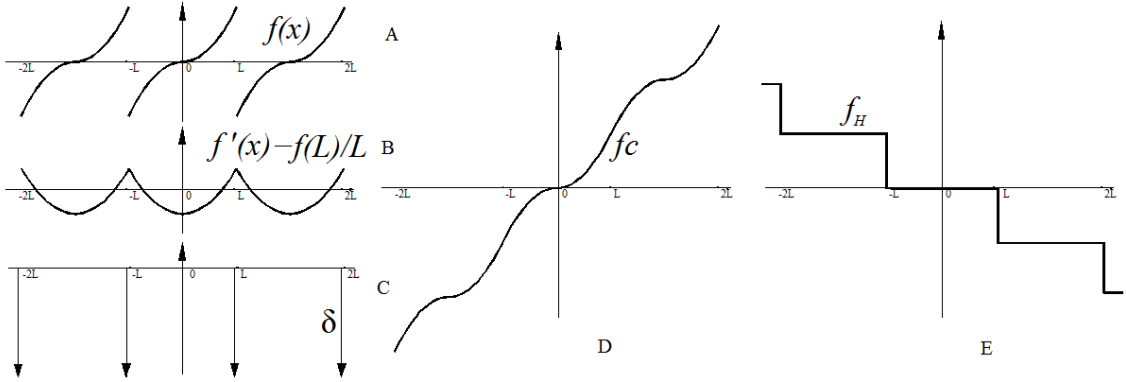

*Fig.1. The derivation of the expansion in SFS of a periodic function with discontinuities*

### 1.3. The Fourier series expansion of the product of any two periodic functions

Are known the relationships to determine the coefficients of the expansion in SFS of the product *p(x)q(x)* of any two periodic square-integrable functions, (which makes that their product is also a function of this type), defined on the interval [−*L*, *L*] [9]. The calculation of these coefficients is possible if the Fourier coefficients of the expansion in series of each of the two functions are known: one must calculate the defined integrals of the type $\int_{-L}^{L} p_e \cos\frac{n\pi x}{L}\,dx$ and $\int_{-L}^{L} p_o \sin\frac{n\pi x}{L}\,dx$. If a certain approximation is allowed, the calculation can be done by numerical methods.

In the case of ODEs, one or both of the terms of the product *p(x)q(x),* can be even one of the functions *y(x), y'(x), y''(x), ∫y(x)dx,* etc., with unknown Fourier coefficients. The replacement of these products, when they appear in a differential equation, with their expansion in Fourier series, or with the Fourier sums $S_N$ who approximates them, leads to a system of *2N+1 (N→∞)* of algebraic equations. Solving them, we can find the coefficients $f_0$, $a_1, ..., a_N, b_1, ..., b_N$ of the Fourier series of the expansion of the function *y(x)*. For

$$y(x) = y_0 + y_e + y_o = y_0 + \sum_{n=1}^{\infty} A_n \cos\frac{n\pi x}{L} + \sum_{n=1}^{\infty} B_n \sin\frac{n\pi x}{L} \text{ and}$$

$$p(x) = p_0 + p_e + p_o = p_0 + \sum_{n=1}^{\infty} C_n \cos\frac{n\pi x}{L} + \sum_{n=1}^{\infty} D_n \sin\frac{n\pi x}{L}$$

(where *$y_{e0}=p_{e0}=0$*), the coefficients of Fourier expansion of the product *yp* are:

$$(yp)^{SFS} = y_0 p_0 + (y_0 \overline{p}_e + p_0 \overline{y}_e + \overline{y}_e \overline{p}_e + \overline{y}_o \overline{p}_o) + (y_0 \overline{p}_o + p_0 \overline{y}_o + \overline{y}_e \overline{p}_o + \overline{y}_o \overline{p}_e) =$$

$$= \overline{P}(x) = P_0 + \overline{P}_e + \overline{Q}_o = P_0 + \sum_{n=1}^{\infty} P_n \cos\frac{n\pi x}{L} + \sum_{n=1}^{\infty} Q_n \sin\frac{n\pi x}{L} \qquad \text{where:} \qquad (1.8)$$

$$P_0 = \frac{1}{2L}\int_{-L}^{L}(yp)^{SFS}\,dx = y_0 p_0 + \frac{1}{2L}\int_{-L}^{L}\overline{y}_e\overline{p}_e\,dx + \frac{1}{2L}\int_{-L}^{L}\overline{y}_o\overline{p}_o\,dx = y_0 p_0 + \frac{1}{2}\left(\sum_{i=1}^{\infty}A_i C_i + \sum_{l=1}^{\infty}B_l D_l\right) \qquad (1.9)$$

$$P_n = \frac{1}{L}\int_{-L}^{L}(y_0\overline{p}_e + p_0\overline{y}_e + \overline{y}_e\overline{p}_e + \overline{y}_o\overline{p}_o)\cos\frac{n\pi x}{L}\,dx = \frac{y_0}{L}\int_{-L}^{L}p_e\cos\frac{n\pi x}{L}\,dx + \frac{p_0}{L}\int_{-L}^{L}y_e\cos\frac{n\pi x}{L}\,dx +$$

$$+ \frac{1}{L}\int_{-L}^{L}\left(\sum_{i=1}^{\infty}A_i\cos\frac{l\pi x}{L}\right)\left(\sum_{m=1}^{\infty}C_m\cos\frac{m\pi x}{L}\right)\cos\frac{n\pi x}{L}\,dx + \frac{1}{L}\int_{-L}^{L}\left(\sum_{i=1}^{\infty}B_i\sin\frac{l\pi x}{L}\right)\left(\sum_{m=1}^{\infty}D_m\sin\frac{m\pi x}{L}\right)\cos\frac{n\pi x}{L}\,dx$$

$$\rightarrow P_n = y_0 C_n + p_0 A_n + \frac{1}{2}\sum_{l=1}^{\infty} A_l\big(C_{n+l} + C_{|n-l|}\big) + \frac{1}{2}\sum_{l=1}^{\infty} B_l\big[D_{n+l} + \text{sgn}(l-n)D_{|n-l|}\big]$$
(1.10a)

$$Q_n = \frac{1}{L}\int_{-L}^{L}(y_0\overline{p}_o + p_0\overline{y}_o + \overline{y}_e\overline{p}_o + \overline{y}_o\overline{p}_e)\sin\frac{n\pi x}{L}\,dx = \frac{y_0}{L}\int_{-L}^{L} p_o\sin\frac{n\pi x}{L}\,dx + \frac{p_0}{L}\int_{-L}^{L} y_o\sin\frac{n\pi x}{L}\,dx +$$

$$+\frac{1}{L}\int_{-L}^{L}\left(\sum_{l=1}^{\infty} A_l\cos\frac{l\pi x}{L}\right)\left(\sum_{m=1}^{\infty} D_m\sin\frac{m\pi x}{L}\right)\sin\frac{n\pi x}{L}\,dx + \frac{1}{L}\int_{-L}^{L}\left(\sum_{l=1}^{\infty} B_l\sin\frac{l\pi x}{L}\right)\left(\sum_{m=1}^{\infty} C_m\cos\frac{m\pi x}{L}\right)\sin\frac{n\pi x}{L}\,dx$$

$$\rightarrow Q_n = y_0 D_n + p_0 B_n + \frac{1}{2}\sum_{l=1}^{\infty} A_l\big[D_{n+l} - \text{sgn}(l-n)D_{|n-l|}\big] - \frac{1}{2}\sum_{l=1}^{\infty} B_l\big(C_{n+l} - C_{|n-l|}\big)$$
(1.10b)

In these relationships, we have considered: $A_0 = B_0 = C_0 = D_0 = 0$ (fo $l=n$). The relations are the same as in [9], in which $A_0 = 2y_0$ and $C_0 = 2p_0$, but here we have chosen another method to achieve the result, using intermediate calculations such as:

$$\int_{-L}^{L}\sum_{l=1}^{\infty}\left(A_l\cos\frac{l\pi x}{L}\right)\sum_{m=1}^{\infty}\left(C_m\cos\frac{m\pi x}{L}\right)\cos\frac{n\pi x}{L}\,dx = \int_{-L}^{L}\sum_{l=1}^{\infty}\left(A_l\cos\frac{l\pi x}{L}\right)\sum_{m=1}^{\infty}\frac{C_m}{2}\left[\cos\frac{(m-n)\pi x}{L} + \cos\frac{(m+n)\pi x}{L}\right]dx =$$

$$= \int_{-L}^{L}\sum_{l=1}^{\infty}\left(A_l\cos\frac{l\pi x}{L}\right)\left[\sum_{m=n+1}^{\infty}\frac{C_m}{2}\cos\frac{(m-n)\pi x}{L} + \frac{C_m}{2} + \sum_{m=1}^{\infty}\frac{C_m}{2}\cos\frac{(m-n)\pi x}{L} + \sum_{m=1}^{\infty}\frac{C_m}{2}\cos\frac{(m+n)\pi x}{L}\right]dx =$$

$$= \frac{L}{2}\left(\sum_{l=1}^{n-1} A_l C_{n-l} + \sum_{l=1}^{\infty} A_l C_{l+n} + \sum_{l=n+1}^{\infty} A_l C_{l-n}\right) = \frac{L}{2}\sum_{l=1}^{\infty} A_l\big(C_{l+n} + C_{|l-n|}\big)$$

Obviously, the following formulas are also valid:

$$P_n = p_0 A_n + y_0 C_n + \frac{1}{2}\sum_{l=1}^{\infty} C_l\big(A_{n+l} + A_{|n-l|}\big) + \frac{1}{2}\sum_{l=1}^{\infty} D_l\big(B_{n+l} + \text{sgn}(l-n)B_{|n-l|}\big)$$
(1.11)

$$Q_n = p_0 B_n + y_0 D_n + \frac{1}{2}\sum_{l=1}^{\infty} C_l\big(B_{n+l} - \text{sgn}(l-n)B_{|n-l|}\big) - \frac{1}{2}\sum_{l=1}^{\infty} D_l\big(A_{n+l} - A_{|n-l|}\big)$$

### 1.4. The calculation of certain infinite numerical series

Each of the coefficients $P_n$ and $Q_n$ of relations 1.9 and 1.10 contain the expressions of infinite numerical series, the sum of which can be determined (if the coefficients $A_n$, $B_n$, $C_n$, and $D_n$ are known) by a volume, sometimes large, of calculations. Between each of these numerical series and the corresponding coefficient of the expansion in SFS of the product $yp$, there is a direct relationship, so that the sum of series can be determined by calculating the respective coefficient. Consequently, if the calculations to determine the Fourier coefficients of the product $yp$ are less voluminous, this calculation method becomes a preferable solution for the calculation of this quantity. To illustrate this, we will compare these calculation volumes for the products between the function:

$$\overline{y}(x) = y_0 + \sum_{n=1}^{\infty} A_n\cos\frac{n\pi x}{L} + \sum_{n=1}^{\infty} B_n\sin\frac{n\pi x}{L} = \overline{x} = \sum_{n=1}^{\infty} 2L\frac{(-1)^{n+1}}{n\pi}\sin\frac{n\pi x}{L} \quad \text{and the functions:}$$

$$\overline{p}_1(x) = p_{1.0} + \sum_{n=1}^{\infty} C_{1.n}\cos\frac{n\pi x}{L} + \sum_{n=1}^{\infty} D_{1.n}\sin\frac{n\pi x}{L} = \overline{x} = \sum_{n=1}^{\infty} 2L\frac{(-1)^{n+1}}{n\pi}\sin\frac{n\pi x}{L} \quad \text{and}$$

$$\overline{p}_2(x) = p_{2.0} + \sum_{n=1}^{\infty} C_{2.n}\cos\frac{n\pi x}{L} + \sum_{n=1}^{\infty} D_{2.n}\sin\frac{n\pi x}{L} = \big(x^2\big)^{SFS} = \frac{L^2}{3} + \sum_{n=1}^{\infty}\frac{4L^2}{\pi^2}\frac{(-1)^n}{n^2}\cos\frac{n\pi x}{L}$$

with the coefficients: $y_0 = 0$, $A_n = 0$, $B_n = 2L\dfrac{(-1)^{n+1}}{n\pi}$; $p_{1.0} = 0$, $C_{1.n} = 0$, $D_{1.n} = 2\dfrac{L}{\pi}\dfrac{(-1)^{n+1}}{n}$; $p_{2.0} = \dfrac{L^2}{3}$,

$C_{2.n} = \dfrac{4L^2(-1)^n}{n^2\pi^2}$, $D_{2.n} = 0$. Result the coefficients of the products:

$$\overline{P}_1 = P_{1.0} + \sum_{n=1}^{\infty} P_{1.n}\cos\frac{n\pi x}{L} + \sum_{n=1}^{\infty} Q_{1.n}\sin\frac{n\pi x}{L} = \big(yp_1\big)^{SFS} = \big(x^2\big)^{SFS} = \frac{L^2}{3} + \frac{4L^2}{\pi^2}\frac{(-1)^n}{n^2}\cos\frac{n\pi x}{L},$$

$$P_2 = P_{2,0} + \sum_{n=1}^{\infty} P_{2,n} \cos\frac{n\pi x}{L} + \sum_{n=1}^{\infty} Q_{2,n} \sin\frac{n\pi x}{L} = (yp_2)^{SFS} = (x^3)^{SFS} = (-1)^{n+1}\left[\frac{2L^3}{n\pi} - \frac{12L^3}{\pi^3 n^3}\right]\sin\frac{n\pi x}{L}$$

We will keep the meanings previously assigned: $A_0 = B_0 = C_0 = D_0 = 0$. Conform to 1.9−1.10:

$$P_{1,0} = \frac{1}{2}\sum_{n=1}^{\infty} B_n \frac{2L(-1)^{n+1}}{n\pi} \quad \text{For } y(x) = x: \to P_{1,0} = 2\sum_{n=1}^{\infty}\frac{L^2}{\pi^2}\frac{(-1)^{n+1}}{n}\frac{(-1)^{n+1}}{n} = \frac{2L^2}{\pi^2}\sum_{n=1}^{\infty}\frac{1}{n^2} = p_{2,0} \tag{1.12}$$

$$P_{1,n} = \frac{B_n}{2}\frac{2L}{\pi}\frac{(-1)^{2n+1}}{2n} + \frac{2L}{2\pi}\sum_{\substack{l=1\\l\neq n}}^{\infty}B_l\left[\frac{(-1)^{n+l+1}}{n+l} + \text{sgn}(l-n)\frac{(-1)^{n-l+1}}{|n-l|}\right] = \frac{L}{\pi}\left[-\frac{B_n}{2n} + \sum_{\substack{l=1\\l\neq n}}^{\infty}(-1)^{n+l}\frac{2lB_l}{n^2-l^2}\right] \tag{1.13}$$

$$\to P_{1,n} = \frac{L}{\pi}\left[\frac{L(-1)^n}{n^2\pi} + \frac{4L}{\pi}\sum_{\substack{l=1\\l\neq n}}^{\infty}\frac{(-1)^n}{l^2-n^2}\right] = (-1)^n\frac{L^2}{\pi^2}\left(\frac{1}{n^2} + 4\frac{3}{4n^2}\right) = 4\frac{(-1)^n L^2}{n^2\pi^2} = C_{2,n}$$

$$Q_{1,n} = \frac{L}{\pi}\left\{y_0\frac{2(-1)^{n+1}}{n} - \frac{A_n}{2n} + \sum_{\substack{l=1\\l\neq n}}^{\infty}A_l\left[\frac{(-1)^{n+l+1}}{n+l} - \text{sgn}(l-n)\frac{(-1)^{n-l+1}}{|n-l|}\right]\right\} \tag{1.14}$$

$$\to Q_{1,n} = \frac{L}{\pi}\left\{y_0\frac{2(-1)^{n+1}}{n} - \frac{A_n}{2n} + 2n\sum_{\substack{l=1\\l\neq n}}^{\infty}(-1)^{n+l+1}\frac{A_l}{n^2-l^2}\right\} = 0$$

$$P_{2,0} = y_0\frac{L^2}{3} + \frac{2L^2}{\pi^2}\sum_{n=1}^{\infty}\frac{(-1)^n A_n}{n^2} \quad \text{For } y(x) = x: \to P_{2,0} = 0 \tag{1.15}$$

$$P_{2,n} = y_0\frac{4L^2}{\pi^2}\frac{(-1)^n}{n^2} + \frac{A_n L^2}{3} + \frac{2L^2}{\pi^2}\left\{\frac{A_n}{4n^2} + \sum_{\substack{l=1\\l\neq n}}^{\infty}A_l\left[\frac{(-1)^{n+l}}{(n+l)^2} + \frac{(-1)^{n-l}}{(n-l)^2}\right]\right\} \to P_{n,2} = 0 \tag{1.16}$$

$$Q_{2,n} = \frac{L^2 B_n}{3} - \frac{2L^2}{\pi^2}\left\{\frac{(-1)^{2n}}{4n^2}B_n + \sum_{\substack{l=1\\l\neq n}}^{\infty}B_l\left[\frac{(-1)^{n+l}}{(n+l)^2} - \frac{(-1)^{n-l}}{(n-l)^2}\right]\right\} \tag{1.17}$$

$$\to Q_{2,n} = \frac{L^2}{3}\frac{2L(-1)^{n+1}}{n\pi} - \frac{2L^2}{\pi^2}\left[\frac{(-1)^{2n}}{4n^2}\frac{2L(-1)^{n+1}}{n\pi} + \sum_{\substack{l=1\\l\neq n}}^{\infty}(-1)^{2l+n+1}\frac{2L}{l\pi}\frac{4nl}{(n^2-l^2)^2}\right] =$$

$$= \frac{2(-1)^{n+1}L^3}{3n\pi} - (-1)^{n+1}\frac{4L^3}{\pi^3}\left[\frac{1}{4n^3} + 4n\left(\frac{11}{16n^4} - \frac{\pi^2}{12n^2}\right)\right] = (-1)^{n+1}\left(\frac{2L^3}{n\pi} - \frac{12L^3}{n^3\pi^3}\right)$$

To obtain these results, we collected data from the specialized literature [10]:

$$\sum_{n=1}^{\infty}\frac{1}{n^2} = \frac{\pi^2}{6}, \quad \sum_{\substack{l=1\\l\neq n}}^{\infty}\frac{1}{l^2-n^2} = \frac{3}{4n^2}, \quad \text{and for } \sum_{\substack{l=1\\l\neq n}}^{\infty}\frac{1}{(l^2-n^2)^2} = \frac{1}{4n^2}\left(\frac{11}{4n^2} - \frac{\pi^2}{3}\right) \text{ we started from:}$$

$$\sum_{l=0}^{\infty}\frac{1}{(l^2-a^2)^2} = \frac{1}{2a^4} + \frac{\pi}{4a^3}ctg(a\pi) + \frac{\pi^2}{4a^2}\cos ec(a\pi), \quad \text{where } a\neq 0, \; 1, \; 2,..., \; \infty. \text{ For } \forall n\in \mathbf{N}, \text{ if}$$

$a = n+\varepsilon \to n$, namely $\varepsilon\in \mathbf{R}\to 0$ and $\sin(a\pi) = \sin(n\pi+\varepsilon\pi) = \sin(\varepsilon\pi)\to 0$. So,

$$\sum_{l=0}^{\infty}\frac{1}{(l^2-a^2)^2} = \sum_{\substack{l=1\\l\neq n}}^{\infty}\frac{1}{(l^2-a^2)^2} + \frac{1}{a^4} + \frac{1}{4an}\left[\frac{1}{(n-a)^2} - \frac{1}{(n+a)^2}\right]$$

$$\lim_{\varepsilon\to 0}\frac{1}{4n(n\pm\varepsilon)}\frac{1}{\varepsilon^2} + \frac{\pi}{4(n\pm\varepsilon)^3}\cot(\pm\varepsilon\pi) + \frac{\pi^2}{4(n\pm\varepsilon)^2}\csc(\pm\varepsilon\pi) = \frac{\pi^2}{12n^2} - \frac{1}{4n^4}$$

## 2. The solution of an ODE, determined by calculating the coefficients of its expansion in sinusoidal Fourier series

To solve the linear differential equations with variable coefficients $a_i(x)$ when these coefficients are analytical functions at point $x=x_0$, we can use the method of expansion in Taylor series: $y(x) = \sum_{n=1}^{\infty} a_n (x - x_0)^n$, if $x_0$ is an ordinary point, or the Frobenius method:

$y(x) = \sum_{n=1}^{\infty} a_n (x - x_0)^{n+r}$, if $x_0$ is a regular singular point [4, 5]. The method, which makes it possible to find a valid solution on its convergence disk of radius $R$, is based on the term by term derivability property of these series and allows to find recurrence relations for the coefficients of the series. This is possible because the base on which the general solution is developed is the same as the base on which its derivatives are also developed. To find the particular solutions it is necessary to know the values of the function and its derivatives at a point of the interval $(-R, R)$. By this method, solving differential equations is transformed into an algebra problem, which most often involves finding the recurrence formulas.

Since SFS and SFN also enjoy the property of term by term differentiation, a similar method to that previously described, can be applied using these new types of series expansions, in equations involving periodic square-integrable functions, defined on an interval $[-L, L]$, if we know the corresponding boundary conditions.

The method of solving differential equations by determining the coefficients of the expansion in sinusoidal or non-sinusoidal series of the unknown functions is a particularly solid method, applicable to all types of differential and integrro-differential equations, linear and nonlinear, of partial differential equations, systems of such equations, whatever their order and whatever the complexity of the coefficients. The conditions required for the application of the method are not very restrictive and easy to fulfill, in particular for the situations encountered in physics and engineering. The method can be applied to large classes of such equations, and even more, the same equation can be solved using several types of expansion. The method can also be easily extended to functions of the complex domain C.

As for the power series expansion method, the derivation and integration operations performed on these series expansions are transformed into algebraic operations performed on the expansion coefficients. Therefore, the resolution of ODE is transformed into resolution of algebraic equations. Of all the bases used as support of expansion, the most advantageous are the sinusoids ($1, cos(nx), sin(nx)$), because of the ease with which the expansion coefficients can be calculated (Euler's formulas), due to the high convergence speed (which makes them ideal for numerical approximation methods) and the relative ease with which the resulting algebraic equations can be solved (the degree of difficulty is dictated by the order of equation and by the degree of non-linearity of the differential equation). Unlike the bases of positive integer powers which, through repeated derivations, "skate" along the elements of the base (leading to the generation of recurrent chains), the sinusoidal bases "oscillate" between the same elements, even and odd, of the base, leading to simpler algebraic equations.

### 2.1. Linear ODEs with *constant* coefficients

As we found in the previous section, on the interval $[-L, L]$, for

$$\bar{y}(x) = y_0 + \bar{y}_e + \bar{y}_o = y_0 + \sum_{n=1}^{\infty} A_n \cos \omega_n x + \sum_{n=1}^{\infty} B_n \sin \omega_n x \qquad (2.1)$$

$$\rightarrow \bar{y}'(x) = \sum_{n=1}^{\infty} \left[ B_n \omega_n + 2(-1)^n \frac{y(L) - y(-L)}{2L} \right] \cos \omega_n x - \sum_{n=1}^{\infty} (A_n \omega_n \sin \omega_n x) + \frac{y(L) - y(-L)}{2L} \qquad (2.2)$$

$$\rightarrow \overline{y}'' = \sum_{n=1}^{\infty}\left[\left(-A_n\omega_n^2 + 2(-1)^n\frac{y'_e(L)}{L}\right)\cos\omega_n x - \left(B_n\omega_n^2 + 2(-1)^n\omega_n\frac{y_o(L)}{L}\right)\sin\omega_n x\right] + \frac{y'_e(L)}{L} \quad (2.3)$$

As well, $y_o(L) = L\lim_{n\to\infty}2\left[(-1)^{n+1}\omega_n b_n\right] = 1/2\cdot\left[y(L) - y(-L)\right]$ \hfill (2.1a)

Further on, by successive derivations (at each derivation, the expansion coefficients are corrected taking into account the existence of the discontinuity points) we find the expressions of series expansions for higher order derivatives, when are known their values at the points at the limit of the definition interval. Therefore, the particular solution y(x) of any ordinary linear differential equation with constant coefficients, homogeneous or inhomogeneous, whatever the order of the equation, if it is defined over any interval [−L, L], or equivalent, can be determined if y(x) is square-integrable, by calculating the coefficients of its expansion in sinusoidal series, if the values of the function and its derivatives are known at the points of the ends of the interval (boundary conditions). The solutions of the equations can be found for other types of boundary conditions too.

To illustrate this, we will apply the method in the case of homogeneous linear equations with constant coefficients, over the interval [−π, π], for some simple equations:

**Example 2.1: *y' = a*,** \hfill (2.4)

with the boundary condition *y(0)=C*. For the chosen interval, we can write:

$$\sum_{n=1}^{\infty}\left[B_n n + 2(-1)^n\frac{y(\pi)-y(-\pi)}{2\pi}\right]\cos nx + \sum_{n=1}^{\infty}\left(-A_n n\sin nx\right) + \frac{y(\pi)-y(-\pi)}{2\pi} = a \quad (2.5)$$

$$\rightarrow A_n=0, \quad B_n = \frac{y(\pi)-y(-\pi)}{2\pi}\frac{2(-1)^{n+1}}{n}, \quad \frac{y(\pi)-y(-\pi)}{2\pi} = \frac{y_o(\pi)}{\pi} = a$$

For the average value $y_0$, there is no conditioning, so it can take any value K. So:

$$y(x) = a\frac{2(-1)^{n+1}}{n}\sin nx + K = ax + K \text{ . For } y=0 \rightarrow y(0)=K=C.$$

For $a=0 \rightarrow y(x)=K$

**Example 2.2: *y' = ay*,** \hfill (2.6)

with the boundary condition *y(π)=C*. From the relations (2.2) and (2.6), we obtain:

$$ay_0 = \frac{y_o(\pi)}{\pi}, \quad aB_n = -nA_n \rightarrow B_n = -nA_n/a \text{ , for } n=1, 2, 3, ...\infty \text{ and}$$

$$aA_n = -\frac{n^2}{a}A_n + 2(-1)^n\frac{y_o(\pi)}{\pi} \rightarrow A_n = 2\frac{(-1)^n a}{n^2+a^2}\frac{y_o(\pi)}{\pi} \rightarrow B_n = -2n\frac{(-1)^n}{n^2+a^2}\frac{y_o(\pi)}{\pi} \text{ and so:}$$

$$\overline{y}(x) = 2\frac{y_o(\pi)}{\pi}\left[\frac{1}{2a} + \sum_{n=1}^{\infty}(-1)^n\frac{a\cos nx - n\sin nx}{n^2+a^2}\right] \quad (2.7)$$

We recognize here, for *y_o(π)=sinhaπ*, the Fourier series expansion of the function *y(x)=e^{ax}*, for which *y(π)=e^{aπ}*, and for *y_o(π)=Ksinhaπ*, the Fourier series expansion of the function *y(x)=K·e^{ax}* (the general solution), for which *y(π)=Ke^{aπ}*, et *y(−π)=Ke^{−aπ}*.

In conclusion, for *y(π)=C → K=Ce^{−aπ}*, with the particular solution *y(x)=Ce^{a(x−π)}*.

During the resolution of the equation, the constraint $aB_n = -nA_n$ appeared. This means that $A_n=0$ imposes $B_n=0$ , so the solution cannot have a single component (it cannot be only even or odd). This fact is imposed even by the relation of equality (2.4). So $f_o(\pi)/\pi\neq0$, therefore, always $y_0\neq0$.

**Example 2.3: *y'' = a²y*,**

with *y(π)=K_1* and *y'(π)=K_2*. Because:

$$\overline{y}'(x) = \overline{y}'_e + \overline{y}'_o = \sum_{n=1}^{\infty}\left(-A_n n\sin nx\right) + \sum_{n=1}^{\infty}\left[B_n n + 2(-1)^n\frac{y_o(\pi)}{\pi}\right]\cos nx + \frac{y_o(\pi)}{\pi} \quad (2.8)$$

$$\overline{y}'' = \sum_{n=1}^{\infty}\left[\left(-A_n n^2 + 2(-1)^n \frac{y'_e(\pi)}{\pi}\right)\cos(nx) - \left(B_n n^2 + 2(-1)^n n \frac{y_o(\pi)}{\pi}\right)\sin(nx)\right] + \frac{y'_e(\pi)}{\pi} \qquad (2.9)$$

and by equality with $a^2 y$ we discover that:

$$y_0 = \frac{y'_e(\pi)}{a^2 \pi}, \ A_n = 2\frac{(-1)^n}{n^2 + a^2}\frac{y'_e(\pi)}{\pi} \text{ et } B_n = -2n\frac{(-1)^n}{n^2 + a^2}\frac{y_o(\pi)}{\pi}, \text{ so}$$

$$\overline{y}(x) = \frac{y'_e(\pi)}{a^2 \pi} + 2\sum_{n=1}^{\infty}\frac{(-1)^n}{n^2 + a^2}\left[\frac{y'_e(\pi)}{\pi}\cos nx - n\frac{y_o(\pi)}{\pi}\sin nx\right] \qquad (2.10)$$

For $y_o(\pi)=sinh(a\pi)$ and $y'_e(\pi)=a \cdot sinh(a\pi)$, we recognize the series expansion of the function: $\overline{y}(x) = \cosh(ax) + \sinh(ax) = e^{ax}$. But, because there are no constraints between the coefficients of the two components (even and odd), each of them can be a solution: $y_1=cosh(ax)$ and $y_2=sinh(ax)$ (for $f_o(\pi)=0$, respectively $f'_e(\pi)=0$). Therefore, the general solution of the equation is $y(x)=C_1cosh(ax)+C_2sinh(ax)$. The particular solution results from: $y(\pi)=C_1cosh(a\pi)+C_2sinh(a\pi)=K_1$ and $y'(\pi)=aC_1sinh(a\pi)+aC_2cosh(a\pi)=K_2$

**Example 2.4: $y'' = -a^2 y$**

Of the relationship (2.10), replacing $a^2$ with $-a^2$, we have:

$$\overline{y}(x) = -\frac{y'_e(\pi)}{a^2 \pi} + 2\frac{(-1)^n}{n^2 - a^2}\left[\frac{y'_e(\pi)}{\pi}\cos nx - n\frac{y_o(\pi)}{\pi}\sin nx\right],$$

in which, for $y_o(\pi)=sin(a\pi)$ and $y'_e(\pi)=-a \cdot sin(a\pi)$ we recognize, for $a \notin Z$, the series expansion of the function:

$$\overline{y}(x) = \cos ax + \sin ax = \frac{\sin a\pi}{a\pi} + 2\frac{(-1)^{n+1}}{n^2 - a^2}\left[\frac{a\sin a\pi}{\pi}\cos nx + n\frac{\sin a\pi}{\pi}\sin nx\right]$$

Therefore, the general solution to this equation is $y(x)=C_1cos(ax)+C_2sin(ax)$.

**Example 2.5: $y''+ay'+by=0$**

By replacing in the equation the unknown function and its derivatives by their expansions in SFS:

$$\sum_{n=1}^{\infty}\left[-A_n n^2 + 2(-1)^n \frac{y'_e(\pi)}{\pi} + aB_n n + 2(-1)^n a\frac{y_o(\pi)}{\pi} + bA_n\right]\cos(nx) +$$

$$+ \sum_{n=1}^{\infty}\left[-B_n n^2 - 2(-1)^n n\frac{y_o(\pi)}{\pi} - aA_n n + bB_n\right]\sin(nx) + \frac{y'_e(\pi)}{\pi} + a\frac{y_o(\pi)}{\pi} + by_0 = 0$$

All the coefficients of this expansion are zero:

$$y_0 = -\frac{y'_e(\pi) + ay_o(\pi)}{b\pi}, \ B_n = A_n\frac{n^2 - b}{an} - \frac{2(-1)^n}{an}\left[\frac{y'_e(\pi)}{\pi} + \frac{ay_o(\pi)}{\pi}\right], \ A_n = -B_n\frac{n^2 - b}{an} - \frac{2(-1)^n}{a}\frac{y_o(\pi)}{\pi}$$

$$\rightarrow A_n = -A_n\left(\frac{n^2 - b}{an}\right)^2 + 2(-1)^n \frac{y'_e(\pi) + ay_o(\pi)}{an\pi}\frac{n^2 - b}{an} - 2(-1)^n\frac{y_o(\pi)}{a\pi}$$

$$A_n\left(n^4 + a^2 n^2 - 2bn^2 + b^2\right)\pi = 2(-1)^n\left[y'_e(\pi)\left(n^2 - b\right) + a\left(2n^2 - b\right)y_o(\pi)\right]$$

$$\rightarrow A_n = \frac{2(-1)^n\left[y'_e(\pi)\left(n^2 - b\right) - aby_o(\pi)\right]}{\pi\left[n^2 + \frac{1}{2}\left(a^2 - 2b - a\Delta\right)\right]\left[n^2 + \frac{1}{2}\left(a^2 - 2b + a\Delta\right)\right]} = \frac{2(-1)^n\left[y'_e(\pi)\left(n^2 - b\right) - aby_o(\pi)\right]}{\pi\left[n^2 + \left(\frac{-a+\Delta}{2}\right)^2\right]\left[n^2 + \left(\frac{-a-\Delta}{2}\right)^2\right]},$$

where $\Delta = \sqrt{a^2 - 4b}$. We will note $\lambda_1 = \frac{-a + \sqrt{a^2 - 4b}}{2} = \frac{-a+\Delta}{2}$ and $\lambda_2 = \frac{-a-\Delta}{2}$

$$B_n = -\left(B_n\frac{n^2 - b}{an}\right)^2 - \frac{2(-1)^n}{a}\frac{y_o(\pi)}{\pi}\frac{n^2 - b}{an} - \frac{2(-1)^n}{an}\left[\frac{y'_e(\pi)}{\pi} + \frac{ay_o(\pi)}{\pi}\right]$$

$$B_n\left(n^4 + a^2 n^2 - 2bn^2 + b^2\right)\pi = 2(-1)^n n\left[y_o(\pi)\left(n^2 + a^2 - b\right) + a y'_e(\pi)\right]$$

$$\rightarrow B_n = \frac{2(-1)^n n\left[y_o(\pi)\left(n^2 + a^2 - b\right) + a y'_e(\pi)\right]}{\pi\left[n^2 + \left(\dfrac{-a+\Delta}{2}\right)^2\right]\left[n^2 + \left(\dfrac{-a-\Delta}{2}\right)^2\right]} = \frac{2(-1)^n n\left[y_o(\pi)\left(n^2 + a^2 - b\right) + a y'_e(\pi)\right]}{\pi\left(n^2 + \lambda_1^2\right)\left(n^2 + \lambda_2^2\right)}$$

For $y_{o.1}(\pi) = C_1 \sinh\dfrac{-a+\Delta}{2}\pi = C_1 \sinh(\lambda_1\pi)$ and $y'_{e.1}(\pi) = C_1\lambda_1 \sinh(\lambda_1\pi) = \lambda_1 y_{o.1}(\pi)$:

$$y_{0.1} = -\frac{C_1\lambda_1 y_{o.1}(\pi) + C_1 a y_{o.1}(\pi)}{b\pi} = -C_1\frac{y_{o.1}(\pi)}{\pi}\frac{-a+\Delta+2a}{2b} = -C_1\frac{y_{o.1}(\pi)}{\pi}\frac{4b}{2b(a-\Delta)} = C_1\frac{\sinh(\lambda_1\pi)}{\pi\lambda_1}$$

$$A_{n.1} = \frac{2(-1)^n C_1 \sinh(\lambda_1\pi)\left[\lambda_1\left(n^2 - b\right) - a\lambda_1\lambda_2\right]}{\pi\left(n^2 + \lambda_1^2\right)\left(n^2 + \lambda_2^2\right)} = 2(-1)^n C_1\frac{\sinh(\lambda_1\pi)}{\pi}\frac{\lambda_1}{n^2 + \lambda_1^2}$$

$$B_{n.1} = \frac{2(-1)^n C_1 n \sinh(\lambda_1\pi)\left(n^2 + a^2 - b + a\lambda_1\right)}{\pi\left(n^2 + \lambda_1^2\right)\left(n^2 + \lambda_2^2\right)} = 2(-1)^n C_1\frac{\sinh(\lambda_1\pi)}{\pi}\frac{n}{n^2 + \lambda_1^2},$$

And for $y_{o.2}(\pi) = C_2 \sinh\dfrac{-a-\Delta}{2}\pi = C_2 \sinh(\lambda_2\pi)$ et $y'_{e.2}(\pi) = C_2\lambda_2 \sinh(\lambda_2\pi) = \lambda_2 y_{o.2}(\pi)$:

$$y_{0.2} = C_2\frac{\sinh(\lambda_2\pi)}{\pi\lambda_2}, \quad A_{n.2} = 2(-1)^n C_2\frac{\sinh(\lambda_2\pi)}{\pi}\frac{\lambda_2}{n^2 + \lambda_2^2}, \quad B_{n.2} = 2(-1)^n C_2\frac{\sinh(\lambda_2\pi)}{\pi}\frac{n}{n^2 + \lambda_2^2}$$

The two solutions are independent. Therefore, the general solution of the equation is:

- for $a^2 > 4b$: $y(x) = C_1\exp\left(\dfrac{-a+\sqrt{a^2-4b}}{2}x\right) + C_2\exp\left(\dfrac{-a-\sqrt{a^2-4b}}{2}x\right)$

- for $a^2 < 4b$: $y(x) = \exp\left(\dfrac{-a}{2}x\right)\left[C_1\sin\left(\dfrac{\sqrt{a^2-4b}}{2}x\right) + C_2\cos\left(\dfrac{\sqrt{a^2-4b}}{2}x\right)\right]$

- for $a^2 < 4b$: $\lambda_1 = \lambda_2 = \dfrac{-a}{2}$. $y_1(x) = C_1\exp\left(\dfrac{-a}{2}x\right)$ is a solution to the equation. We can look for a second solution, in the form $y_2(x) = u(x)y_1(x)$. By the method "reduction of order" [4] we find for $u(x) = u_0 + \sum\limits_{n=1}^{\infty}C_n\cos\dfrac{n\pi x}{L} + \sum\limits_{n=1}^{\infty}D_n\sin\dfrac{n\pi x}{L}$ the condition:

$$\overline{u''} = \sum_{n=1}^{\infty}\left[\left(-C_n n^2 + 2(-1)^n\frac{u'_e(\pi)}{\pi}\right)\cos(nx) - \left(D_n n^2 + 2(-1)^n n\frac{u_o(\pi)}{\pi}\right)\sin(nx)\right] + \frac{u'_e(\pi)}{\pi} = 0$$

$$\rightarrow \frac{u'_e(\pi)}{\pi} = 0, \ -C_n n^2 + 2(-1)^n\frac{u'_e(\pi)}{\pi} = 0 \rightarrow C_n = 0, \ D_n = 2(-1)^n\frac{u_o(\pi)}{n\pi} \rightarrow u(x) = C_1 x + C_2$$

We note that for all ODEs with constant coefficients, whatever the degree of the equation, the expansion coefficients of the unknown function are provided by an equation to determine the mean value and a pair of equations for each harmonic of the expansion.

**Example 2.6: $y'=x$**

Even in the case of non-homogeneous ODEs with constant coefficients, there are no significantly more difficult equations to solve. For example, if in this equation, on the interval $[-L, L]$, we make the replacements:

$$\overline{y}(x) = y_0 + \sum_{n=1}^{\infty}A_n\cos nx + \sum_{n=1}^{\infty}B_n\sin nx \ \text{ et } \ \overline{x} = 2\sum_{n=1}^{\infty}\frac{(-1)^{n+1}}{n}\sin nx, \text{ we have:}$$

$$\sum_{n=1}^{\infty}\left(-A_n n\sin nx\right) + \sum_{n=1}^{\infty}n\left[B_n + 2(-1)^n\frac{y_o(\pi)}{n\pi}\right]\cos nx + \frac{y_o(\pi)}{\pi} = 2\sum_{n=1}^{\infty}\frac{(-1)^{n+1}}{n}\sin nx$$

$\rightarrow y_0 = K$ (arbitrary), $y_o(\pi) = 0$, $\rightarrow B_n + 2(-1)^n \dfrac{y_o(\pi)}{n\pi} = 0 \rightarrow B_n = 0$ and $A_n = -2\dfrac{(-1)^{n+1}}{n^2}$

We recognize the expansion: $\dfrac{x^2}{2} = \dfrac{1}{6} - 2\sum_{n=1}^{\infty}\dfrac{(-1)^{n+1}}{n^2}\cos nx \rightarrow y(x) = y_0 + \dfrac{x^2}{2} - \dfrac{1}{6} = \dfrac{x^2}{2} + C$

**Example 2.7: $y' = x^2$.**

$\sum_{n=1}^{\infty}(-A_n n \sin nx) + \sum_{n=1}^{\infty}\left[ B_n n + 2(-1)^n \dfrac{y_o(\pi)}{\pi}\right]\cos nx + \dfrac{y_o(\pi)}{\pi} = \dfrac{\pi^2}{3} - 4\sum_{n=1}^{\infty}\dfrac{(-1)^{n+1}}{n^2}\cos nx$

$\rightarrow y_0 = K$ (arbitrary), $\rightarrow y_o(\pi) = \dfrac{\pi^3}{3}$, $A_n = 0$, $B_n = \dfrac{2(-1)^{n+1}}{n}\dfrac{\pi^2}{3} - 4\sum_{n=1}^{\infty}\dfrac{(-1)^{n+1}}{n^3}$,

$\rightarrow y(x) = y_0 + \left[\dfrac{2(-1)^{n+1}}{n}\dfrac{\pi^2}{3} - 4\sum_{n=1}^{\infty}\dfrac{(-1)^{n+1}}{n^3}\right]\cos nx = \dfrac{x^3}{3} + C$

### 2.2. Linear ODEs with variable coefficients

A linear differential equation of order $m$ with variable, non-homogeneous coefficients, of the form $f_m(x)y^{(m)} + ... + f_2(x)y'' + f_1(x)y' + f_0(x)y = g(x)$ can be solved, on the interval $[-L, L]$, if $y^{(i)}(x)$, $f_i(x)$ and $g(x)$, $i = 0, 1, 2, ..., m$, are the square-integrable functions, by determining the coefficients $A_n$ and $B_n$, $n = 1, 2, ..., \infty$, of the expansion in SFS of the unknown function:

$$y(x) = y_0 + \sum_{n=1}^{\infty}A_n \cos \omega_n x + \sum_{n=1}^{\infty}B_n \sin \omega_n x$$

if the values of the function and its derivatives are known at the points situated at the ends of the interval. If we know only their value at one of these extremes, the relation (2.1a) allows us to find the other. If $L$ is taken as a parameter, estimation can be made, if $L \rightarrow \infty$, for all the values on the real axis for which the conditions of square-integrability are fulfilled. From the relations (2.1) - (2.3) of the previous subsection, by successive derivations, we will find the expressions of the series expansions for the higher order derivatives (if their values at the ends of the definition interval are known). By using relations similar to the relations (1.9) and (1.10) previously determined (relations for the coefficients of the product of the two functions), as well as relations deduced from the integration formula by parts one can obtain the expressions (dependent on $A_n$ and $B_n$) for the coefficients $P_n$ and $Q_n$ of the expansion in series for the terms of the form $f_i(x)y^{(i)}$, which are introduced into the basic equation, simultaneously with the expression of the expansion in series of the function $g(x)$. By simplifying the resulting relation (grouping the terms which have one of the elements of the expansion base as a common factor, including here also the unitary function $I$), we obtain the expression of the expansion in SFS of an identically zero function. Consequently, all the coefficients of this expansion (algebraic expressions in which the coefficients $A_n$ and $B_n$ appear) are zero and give rise to relations (algebraic equations) which allow the calculation of the unknowns $A_n$, $B_n$ et $f_0$. The method is similar to that in which these terms are determined by replacing the function $y(x)$ by its expansion in Taylor series.

The coefficients of order $N$, $P_N$ and $Q_N$, of product $f_i(x)y^{(i)}$ expansions are numerical series which can also contain infinite quantities of decreasing or alternating decreasing terms. When these sums can be calculated (they can be reduced to a finite sum of terms), we obtain exact equations, and by solving them, we obtain exact values for the expansion coefficients of the unknown function $y(x)$. For this, it is necessary that the expansion coefficients in SFS of the functions $f_i(x)$ can be calculated with exactness (the Euler integrals that correspond to them can be calculated with a finite number of terms). On the contrary, we can limit the number of terms in these series to the first $N$. Thus, $2N+1$ approximate equations will be obtained, with $2N+1$ unknown variables, by the resolution of which are obtained approximate values for $y_0$, $A_n$ et $B_n$. For the function $y(x)$ we obtain a value approximated by the Fourier sum $S_N$. The higher the number $N$ of terms is, the better the approximation.

In some cases, the solutions of the equations for some of the coefficients can be deduced by comparing the terms of the equations with the terms of the series expansions of certain known functions, which facilitates the search for the global solution.

A frequent situation is that of equations with polynomial coefficients. If the boundary conditions are given (the values of the unknown function y and its derivatives at the ends of the interval considered), they are sufficient to write the equations which determine the coefficients $A_n$ și $B_n$. If $y(x)$, $y^{(m)}(x)$, $f_i(x)$ and $g(x)$ are square-integrables on the interval $[-L,L]$, the relations (1.12) - (1.17) and those which can be derived from it, provide relations of calculation for all the coefficients of the expansion in SFS of the terms of these equations. Here, for example, are the relations calculated for the terms of the first and second order equations, having first and second degree polynomial coefficients. If

$$\bar{y}(x) = y_0 + \sum_{n=1}^{\infty} A_n \cos\frac{n\pi x}{L} + \sum_{n=1}^{\infty} B_n \sin\frac{n\pi x}{L} \text{ , then:}$$

$$(yx)^{SFS} = \frac{L}{\pi}\sum_{n=1}^{\infty}\frac{(-1)^{n+1}}{n}B_n + \frac{L}{\pi}\left[-\frac{B_n}{2n} + \sum_{\substack{l=1\\l\neq n}}^{\infty}(-1)^{l+n}\frac{2lB_l}{n^2-l^2}\right]\cos\frac{n\pi x}{L} +$$

$$+\frac{L}{\pi}\left[y_0\frac{2(-1)^{n+1}}{n} - \frac{A_n}{2n} + 2n(-1)^{n+1}\sum_{\substack{l=1\\l\neq n}}^{\infty}\frac{(-1)^l A_l}{n^2-l^2}\right]\sin\frac{n\pi x}{L}$$

$$(yx^2)^{SFS} = L^2\left[\frac{y_0}{3} + \sum_{n=1}^{\infty}\frac{2(-1)^n A_n}{n^2\pi^2}\right] + \frac{L^2}{\pi^2}\left[\frac{4y_0(-1)^n}{n^2} + \frac{A_n(2\pi^2n^2+3)}{6n^2} + 4\sum_{\substack{l=1\\l\neq n}}^{\infty}A_l(-1)^{n+l}\frac{(n^2+l^2)}{(n^2-l^2)^2}\right]\cos\frac{n\pi x}{L} +$$

$$+\frac{L^2}{\pi^2}\left[\frac{B_n(2\pi^2n^2+3)}{6n^2} + 8n\sum_{\substack{l=1\\l\neq n}}^{\infty}(-1)^{n+l}\frac{lB_l}{(n^2-l^2)^2}\right]\sin\frac{n\pi x}{L}$$

$$(y'x)^{SFS} = \sum_{n=1}^{\infty}(-1)^n A_n + \left[\frac{A_n}{2} - \sum_{\substack{l=1\\l\neq n}}^{\infty}(-1)^{l+n}\frac{2l^2A_l}{n^2-l^2}\right]\cos\frac{n\pi x}{L} + \left[-\frac{B_n}{2} + 2n(-1)^{n+1}\sum_{\substack{l=1\\l\neq n}}^{\infty}\frac{lB_l(-1)^l}{n^2-l^2}\right]\sin\frac{n\pi x}{L}$$

$$(y'x^2)^{SFS} = \frac{y_o(L)L}{3} + \frac{2L}{\pi^2}\sum_{n=1}^{\infty}\left[\frac{(-1)^n\pi B_n}{n} + \frac{2y_o(L)}{n^2}\right] +$$

$$+\sum_{n=1}^{\infty}\left\{\frac{(-1)^n y_o(L)(4\pi^2n^2+30)}{6n^2L} + \frac{B_n(2\pi^2n^2+3)n\pi}{6n^2L} + 4(-1)^n\sum_{\substack{l=1\\l\neq n}}^{\infty}\left[\frac{\pi}{L}B_l(-1)^l\frac{(n^2+l^2)}{(n^2-l^2)^2} + 2\frac{y_o(L)}{L}\frac{(n^2+l^2)}{(n^2-l^2)^2}\right]\right\}\cos\frac{n\pi x}{L} +$$

$$+\sum_{n=1}^{\infty}\frac{L}{\pi}\left[\frac{-A_n(2\pi^2n^2-3)}{6n} - 8n\sum_{\substack{l=1\\l\neq n}}^{\infty}(-1)^{n+l}\frac{l^2A_l}{(n^2-l^2)^2}\right]\sin\frac{n\pi x}{L}$$

$$(y''x)^{SFS} = \sum_{n=1}^{\infty}\left[(-1)^n\frac{n\pi}{L}B_n + 2\frac{y_o(L)}{L}\right] + \sum_{n=1}^{\infty}\left\{\frac{B_n n\pi}{2L} + (-1)^n\frac{y_o(L)}{L} + (-1)^n\sum_{\substack{l=1\\l\neq n}}^{\infty}\frac{(-1)^{l+1}2l^2}{n^2-l^2}\left[lB_l\frac{\pi}{L} - 2\frac{y_o(L)}{L}\right]\right\}\cos\frac{n\pi x}{L} +$$

$$+\sum_{n=1}^{\infty}\left\{\frac{3(-1)^{n+1}y'_x(L)}{n\pi} + \frac{A_n}{2}\frac{n\pi}{L} + 2n(-1)^{n+1}\sum_{\substack{l=1\\l\neq n}}^{\infty}\left[\frac{(-1)^{l+1}A_l l^2}{n^2-l^2}\frac{L}{L} + 2\frac{y'_x(L)}{\pi}\frac{1}{n^2-l^2}\right]\right\}\sin\frac{n\pi x}{L}$$

$$(y''x^2)^{SFS} = y'_e(L)L - 2\sum_{i=1}^{\infty}(-1)^i A_i +$$

$$+ \sum_{n=1}^{\infty}\left\{2(-1)^n \frac{y'_e(L)}{L}\frac{L^2\left(2\pi^2 n^2+15\right)}{6\pi^2 n^2} - \frac{A_n\left(2\pi^2 n^2+3\right)}{6} + (-1)^n\sum_{\substack{l=1\\l\neq n}}^{\infty}\left[A_l(-1)^{l+1}\frac{4l^2\left(n^2+l^2\right)}{\left(n^2-l^2\right)^2} + \frac{8L^2}{\pi^2}\frac{y'_e(L)}{L}\frac{\left(n^2+l^2\right)}{\left(n^2-l^2\right)^2}\right]\right\}\cos\frac{n\pi x}{L} +$$

$$+ \sum_{n=1}^{\infty}\left\{\frac{\left(2\pi^2 n^2-3\right)}{6\pi n}\left[-\pi n B_n - 2(-1)^n y_o(L)\right] - 8n\sum_{\substack{l=1\\l\neq n}}^{\infty}\left[\frac{l^3 B_l(-1)^{n+l}}{\left(l^2-n^2\right)^2}\right] - \frac{16nL}{\pi}\frac{y_o(L)}{L}\sum_{\substack{l=1\\l\neq n}}^{\infty}\frac{(-1)^l l^2}{\left(l^2-n^2\right)^2}\right\}\sin\frac{n\pi x}{L}$$

For example, a non homogeneous Euler type equation $y''x^2 + ay'x + by = g(x)$, after having used these relations and after the expansion in SFS of the function

$$g(x) = g_0 + \sum_{n=1}^{\infty}P_n\cos\frac{n\pi x}{L} + \sum_{n=1}^{\infty}Q_n\sin\frac{n\pi x}{L}, \text{ leads to the following algebraic equations:}$$

$$y'_e(L)L - 2\sum_{n=1}^{\infty}(-1)^n A_n + a\sum_{n=1}^{\infty}(-1)^n A_n + by_0 = g_0 \;\rightarrow\; y_0 = \frac{1}{b}\left[g_0 - y'_e(L)L + (2-a)\sum_{n=1}^{\infty}(-1)^n A_n\right],$$

then, for each $n \in \boldsymbol{N}^+$:

$$2(-1)^n\frac{y'_e(L)}{L}\frac{L^2\left(2\pi^2 n^2+15\right)}{6\pi^2 n^2} - \frac{A_n\left(2\pi^2 n^2+3\right)}{6} + (-1)^n\sum_{\substack{l=1\\l\neq n}}^{\infty}\left[A_l(-1)^{l+1}\frac{4l^2\left(n^2+l^2\right)}{\left(n^2-l^2\right)^2} + \frac{8L^2}{\pi^2}\frac{y'_e(L)}{L}\frac{\left(n^2+l^2\right)}{\left(n^2-l^2\right)^2}\right] +$$

$$+ a\frac{A_n}{2} - a\sum_{\substack{l=1\\l\neq n}}^{\infty}(-1)^{l+n}\frac{2l^2 A_l}{n^2-l^2} + bA_n = P_n$$

$$\frac{\left(2\pi^2 n^2-3\right)}{6\pi n}\left[-\pi n B_n - 2(-1)^n y_o(L)\right] - 8n\sum_{\substack{l=1\\l\neq n}}^{\infty}\left[\frac{l^3 B_l(-1)^{n+l}}{\left(l^2-n^2\right)^2}\right] - \frac{16nL}{\pi}\frac{y_o(L)}{L}\sum_{\substack{l=1\\l\neq n}}^{\infty}\frac{(-1)^l l^2}{\left(l^2-n^2\right)^2} +$$

$$- a\frac{B_n}{2} + 2an(-1)^{n+1}\sum_{\substack{l=1\\l\neq n}}^{\infty}\frac{lB_l(-1)^l}{n^2-l^2} + bB_n = Q_n$$

### 2.3. Nonlinear ODEs

The determination of the solution by calculating the coefficients of its expansion in sinusoidal Fourier series can also be successfully applied to non-linear ODEs, with polynomial non-linearities, equations whose terms contain only natural powers of the unknown, its derivatives and their products. The coefficients of these terms can be constant, or variable. From expansion

$$\overline{y}(x) = y_0 + \sum_{n=1}^{\infty}A_n\cos\frac{n\pi x}{L} + \sum_{n=1}^{\infty}B_n\sin\frac{n\pi x}{L},$$

the relations of type (2.2) - (2.3), deduced for calculation the Fourier coefficients of the derivatives of any order and of the type $(1.9 - (1.10)$, for calculation the Fourier coefficients of the product of two certain functions, are sufficient to convert the given equation to a system of $2N+1$ algebraic equations($N\to\infty$), with $2N+1$ unknowns: $y_0, A_1, A_2, ..., A_N, B_1, B_2, ..., B_N$. For exemple:

$$(y^2)^{SFS} = y_0^2 + \frac{1}{2}\sum_{l=1}^{\infty}\left(A_l^2+B_l^2\right) + \left[2y_0 A_n + \frac{1}{2}\sum_{l=1}^{\infty}A_l\left(A_{n+l}+A_{|n-l|}\right) + \frac{1}{2}\sum_{l=1}^{\infty}B_l\left(B_{n+l}-\text{sgn}(l-n)B_{|n-l|}\right)\right]\cos\frac{n\pi x}{L} +$$

$$+ \left[2y_0 B_n + \frac{1}{2}\sum_{l=1}^{\infty}A_l\left(B_{n+l}-\text{sgn}(l-n)B_{|n-l|}\right) - \frac{1}{2}\sum_{l=1}^{\infty}B_l\left(A_{n+l}-A_{|n-l|}\right)\right]\sin\frac{n\pi x}{L}$$

$$\left(yy'\right)^{SFS} = y_0\,\frac{y_o(L)}{L} + 2y_0 A_n \cos\frac{n\pi x}{L} +$$

$$+ \sum_{n=1}^{\infty}\sum_{l=1}^{\infty}\frac{1}{2}\left\{A_l\left[\omega_{n+l}B_{n+l}+\omega_{|n-l|}B_{|n-l|}+4(-1)^{l+n}\frac{y_o(L)}{L}\right]-B_l\left[\omega_{n+l}A_{n+l}-\mathrm{sgn}(l-n)\omega_{|n-l|}A_{|n-l|}\right]\right\}\cos\frac{n\pi x}{L} -$$

$$+ \sum_{n=1}^{\infty}\left\{y_0\omega_n A_n + \frac{y_o(L)}{L}B_n - \frac{1}{2}\sum_{l=1}^{\infty}A_l\left[\omega_{n+l}A_{n+l}-\mathrm{sgn}(l-n)\omega_{|n-l|}A_{|n-l|}\right]-\frac{1}{2}\sum_{l=1}^{\infty}B_l\left[\omega_{n+l}B_{n+l}-\omega_{|n-l|}B_{|n-l|}\right]\right\}\sin\frac{n\pi x}{L}$$

$$\left(y'^2\right)^{SFS} = \frac{y_o(L)^2}{L^2} + \frac{1}{2}\sum_{l=1}^{\infty}\left[\omega_l^2 A_l^2 + \omega_l^2 B_l^2 + 4(-1)^l B_l\frac{y_o(L)}{L}+4\frac{y_o(L)^2}{L^2}\right]+2y_0\sum_{n=1}^{\infty}\left[\omega_n B_n + 2(-1)^n\frac{y_o(L)}{L}\right]\cos\frac{n\pi x}{L}+$$

$$+ \sum_{n=1}^{\infty}\sum_{l=1}^{\infty}\frac{1}{2}\left\{\left[\omega_l B_l + 2(-1)^l\frac{y_o(L)}{L}\right]\left[-\omega_{n+l}B_{n+l}-\omega_{|n-l|}B_{|n-l|}+4(-1)^{n+l}\frac{y_o(L)}{L}\right]+\omega_l A_l\left[\omega_{n+l}A_{n+l}-\mathrm{sgn}(l-n)\omega_{|n-l|}A_{|n-l|}\right]\right\}\cos\frac{n\pi x}{L} +$$

$$+ \sum_{n=1}^{\infty}\left\{-2y_0\omega_n A_n + \frac{1}{2}\sum_{l=1}^{\infty}\left[\omega_l B_l + 2(-1)^l\frac{y_o(L)}{L}\right]\left(-\omega_{n+l}A_{n+l}+\mathrm{sgn}(l-n)\omega_{|n-l|}A_{|n-l|}\right)-\omega_l A_l\left[\omega_{n+l}B_{n+l}-\omega_{|n-l|}B_{|n-l|}\right]\right\}\sin\frac{n\pi x}{L}$$

It can be seen that to obtain the values of the Fourier coefficients of the unknown *y(x)*, it is necessary to solve some algebraic equations of the second order or more, if the degree of nonlinearity increases. It can also be noticed that if, in one of the types of equations analyzed, the substitutions made to solve the equation are made using the partial sums $S_N$ (with *2N+1* terms) instead of the whole expansion (with a infinite number of terms), an approximate ODE is obtained. The larger is *N,* the better is the approximation. To solve it, a system of *2N+1* algebraic equations must be solved, a problem for which there are many electronic calculation programs. Therefore, it is possible to develop a simple, universal and accurate algorithm to solve these types of ODS.

### 3. Methods for solving differential equations

The examples discussed in the previous section highlighted that by replacing certain ODE terms with their SFS / SFN expansions, the initial equation remains, or may become, linear. In the general case, the unknown function *y(x)* is replaced by an infinite sum:

$$\hat{y}(x) = y_0\cdot 1 + \sum_{n=-\infty}^{\infty}C_n R_n(x), \text{ sau } \hat{y}(x) = y_0\cdot 1 + \sum_{n=1}^{\infty}A_n P_n(x) + \sum_{n=1}^{\infty}B_n Q_n(x)$$

which can be written as well:

$$y(x) = y_0 + \sum_{n=-\infty}^{\infty}y_n(x), \text{ respectively } y(x) = y_0 + \sum_{n=1}^{\infty}y_{en}(x) + \sum_{n=1}^{\infty}y_{on}(x)$$

If all the terms of the given equation can be replaced by combinations of the functions $R_n(x)$, respectively $P_n(x)$ and $Q_n(x)$, an superposition of independent equations is obtained which, due to the linearity of the equation, is reduced to a system of algebraic equations (the functions $y_n(x)$, respectively $y_{en}(x)$ and $y_{on}(x)$ are perfectly determined by determining the coefficients $C_n$, respectively $A_n$ and $B_n$).

The examples analyzed (2.1-2.7) also highlighted the fact that this method makes it possible to determine, in many cases, an exact solution, expressed in a closed form. In many other situations (for example, for variable coefficient ODEs), exact solutions are obtained, expressed in the form of infinite sums of terms. To be useful in practice, these forms can be approximated by neglecting the least significant terms: the **approximation of the solution**. The approximate solution is determined at all points in the interval, without interpolation, as in many other common approximate methods.

In many other situations, to solve the equation, the approximation of the expansion in series must be done **before** arriving at the exact solution: the impossibility of solving the system of algebraic equations wiht infinite dimensions, resulted by the ODE processing, the impossibility of determining by analytical calculations the Fourier coefficients, etc. We can consider that we used an **approximation of the method**. In these situations, we also obtain a solution expressed in a closed form (a finite sum of the first $2N+1$ harmonics). In this case too, the values of the solution are perfectly determined over the entire definition interval.

There are also situations in which the replacement of certain terms of ODEs is done by other expressions which approximate them, on the basis of empirical, methodological or other criteria. These are methods of **approximating the equation**. In most cases, this approximation is made in a finite number of points, the value of the solution at the intermediate points being determined by interpolation.

### 4. Solutions of ODEs, determined by calculating the coefficients of its expansion in non-sinusoidal Fourier series

In the paper [1], we proposed a generalization of the expansion in SFS of the periodic functions defined on an interval $[-L, L]$, generalization by which we replaced the periodic sinusoidal base (the usual name is trigonometric base) with a non-sinusoidal basis, comprising the following functions: the unit function *1*, the fundamental periodic quasi-harmonics *g(x)*-even and *h(x)*-odd, with the period *2L*, with null mean values over the interval of definition and the secondary quasi-harmonics, defined on $[-L, L]$, $g_n(x)=g(nx)$ and $h_n(x)=h(nx)$, with the period *2L/n*, for $n \in \mathbf{Z}^+$. The fundamental quasi-harmonics *g(x)* and *h(x)*, defined over the interval $[-L, L]$, can be any function admitting expansions in sinusoidal series, extended on the real axis by successive translations: for all $x_R \in R$ we define the function $K(x_R)=E[(x_R -x_1)/2L]$, such that for each $x_R \in R$ and each $x \in [-L, L]$, there are the relations $x_R=x+KT$ and $g(x_R)=g(x_R-KT)=g(x)$. $E(x)=\lfloor x \rfloor=k$, is the *floor function* (*k* is the nearest integer less than or equal to *x*, namely $E(x_R) \le x < E(x)+1$. The secondary quasi-harmonics are obtained by the dilation of these fundamentals: $g_n(x)=g_n[x-2L \cdot E[(x+L)n/2L]$. The function $g_n(x)$ receives on the interval $[2L(2k-1)/n, 2L(2k+1)/n]$, the same values with those that *g(x)* receives on the interval $[-L, L]$. We have introduced the simplified notation for this function:
$g_n(x)=G[-L/n<g(nx)>L/n]_n$, où $n \in \mathbf{N}^+$.

The coefficients $A_n$ and $B_n$ of the expansion in SFN of the function *f(x)* are obtained by means of the algebraic relationships between the Fourier coefficients of the expansions in SFS of the functions *f(x), g(x)* and *h(x)*.

Article [1] highlighted a wide range of possible applications for this new type of series expansion, highlighting the diversity of possible solutions for each of them, illustrated by the different ways of approximating the functions. Concerning the resolution of differential equations, in the previous section, we analyzed a new methods, replacing the unknown function with its expansion in SFS and solving the resulting algebraic equations. The method is also applicable to expansions in SFN, provided that through the derivation/ integration operations any item of another base appears (therefore, only for the bases composed of circular, hyperbolic, exponential functions) . By using this method, the advantages offered by the possibility of using orthogonalized bases can also be capitalized, although this implies additional computational volume.

### 4.1. Methods of solving differential equations by approximating the solution

In some ODEs, it is not necessary to replace all the terms of the differential equation by their expansions in sinusoidal/non-sinusoidal Fourier series. The equation that results by replacing only certain terms, even when it is not linear, can lead to an equation easier to solve than the initial equation, using classical integration/derivation procedures (the method is frequently applied in current practice by techniques of linearization of nonlinear equations). We will illustrate these statements by solving a very simple inhomogeneous equation:

$$y' = x \tag{4.1}$$

whith the condition $y(0) = y_0$

We will initially select a series of Fourier series expansions [1], which will be useful:

$$\bar{x} = X\left[-L < x > L\right] = \frac{2L}{\pi} \sum_{n=1}^{\infty} \frac{(-1)^{n+1}}{n} \sin\frac{n\pi x}{L} \tag{4.2}$$

$$\left(x^2\right)^{SFS} = X_e^2\left[-L > x^2 < L\right] = \frac{L^2}{3} - \frac{4L^2}{\pi^2} \sum_{n=1}^{\infty} \frac{(-1)^{n+1}}{n^2} \cos\frac{n\pi x}{L} \tag{4.3}$$

$$\left(1_o\right)^{SFS} = J\left[-L > -1 < 0 > 1 < L\right] = \frac{4}{\pi} \sum_{n=1}^{\infty} \frac{1}{2n-1} \sin\frac{(2n-1)\pi x}{L} \tag{4.4}$$

$$\left(x^{\prime\prime}\right)^{SFS} = X_e\left[-L > -x - \frac{L}{2} < 0 > x - \frac{L}{2} < L\right] = -\frac{4L}{\pi^2} \sum_{n=1}^{\infty} \frac{1}{(2n-1)^2} \cos\frac{(2n-1)\pi x}{L} \tag{4.5}$$

From the first two relations, we can write the expansions in non-sinusoidal series:

$$\hat{x} = X\left[-L < x > L\right] = \frac{L}{2}\left(J_1 - \sum_{n=1}^{\infty} \frac{1}{2^n} J_{2^n}\right) =$$

$$= \frac{L}{2} J\left[-L > -1 < 0 > 1 < L\right] - \frac{L}{2} \sum_{n=1}^{\infty} \frac{1}{2^n}\left[-\frac{L}{2^n} > -1 < 0 > 1 < \frac{L}{2^n}\right]_{2^n} \tag{4.6}$$

$$\left(x^2\right)^{SFN} = X_e^2\left[-L > x^2 < L\right] = \frac{L^2}{3} + L \cdot X_{e1} - \sum_{n=1}^{\infty} \frac{L}{4^n} \cdot \left[X_e\right]_{2^n} =$$

$$= \frac{L^2}{3} + L\left[-L > -x - \frac{L}{2} < 0 > x - \frac{L}{2} < L\right]_1 - \sum_{n=1}^{\infty} \frac{L}{4^n}\left[-\frac{1}{2^n} > -2^n x - \frac{L}{2} < 0 > 2^n x - \frac{L}{2} < \frac{1}{2^n}\right]_{2^n} \tag{4.7}$$

For the interval $[-L, L]$ the inhomogeneous term of equation (4.1) can be developed in a sinusoidal series, or in a non-sinusoidal series (for the last variant the possibilities being multiple), so that using the reduced notation, we can write:

$$\left(y'\right)^{SFS} = \frac{2L}{\pi} \sum_{n=1}^{\infty} \frac{(-1)^{n+1}}{n}\left[-\frac{L}{n} < \sin\frac{n\pi x}{L} > \frac{L}{n}\right]_n \text{, respectively} \tag{4.8}$$

$$\left(y'\right)^{SFN} = \frac{L}{2}\left[-L < -1 > 0 < 1 > L\right] - \frac{L}{2} \sum_{n=1}^{\infty} \frac{1}{2^n}\left[-\frac{L}{2^n} < -1 > 0 < 1 > \frac{L}{2^n}\right]_{2^n} \tag{4.9}$$

By integration:

$$y + C = -\frac{2L^2}{\pi^2} \sum_{n=1}^{\infty} \frac{(-1)^{n+1}}{n^2}\left[-\frac{L}{n} < \cos\frac{n\pi x}{L} > \frac{L}{n}\right]_n = \frac{1}{2}\left(X_e^2\left[-L > x^2 < L\right] - \frac{L^2}{3}\right) = \frac{x^2}{2} - \frac{L^2}{6},$$

$$y + C = \frac{L}{2}\left[-L < -x > 0 < x > L\right] - \frac{L}{2} \sum_{n=1}^{\infty} \frac{1}{2^n}\left[-\frac{L}{2^n} < -x > 0 < x > \frac{L}{2^n}\right]_{2^n} =$$

$$= \frac{L}{2}\left(\left[-L < -x - \frac{L}{2} > 0 < x - \frac{L}{2} > L\right] + \frac{L}{2}\right) - \frac{L}{2} \sum_{n=1}^{\infty}\left(\frac{1}{4^n}\left\{\left[-\frac{L}{2^n} < -2^n x - \frac{L}{2} > 0 < 2^n x - \frac{L}{2} > \frac{L}{2^n}\right]_{2^n} + \frac{L}{2}\right\}\right) =$$

$$= \frac{1}{2}\left(X_{e1} - \sum_{n=1}^{\infty} \frac{1}{4^n}\left[X_e\right]_{2^n}\right) + \frac{L}{2}\left(\frac{L}{2} - \sum_{n=1}^{\infty} \frac{1}{4^n}\frac{L}{2}\right) = \frac{1}{2}\left(x^2 - \frac{L^2}{3}\right) + \frac{L^2}{4}\left(1 - \frac{1}{3}\right) = \frac{x^2}{2} - \frac{L^2}{6} + \frac{L^2}{6} = \frac{x^2}{2}$$

In both variants, the particular solution is $y - y_0 = x^2/2$.

Note that the derivative of the even function $y_e(x)$ is an odd function which can be developed in a series of odd functions, with zero mean value, with the base $h(x)$:

$$y_e'(x) = \sum_{n=1}^{\infty} A_n h_n(x)$$

The integration of this relation leads to a relation for the expansion of the even function $y_e(x)$, in a base $g(x)$ of even functions, for which the mean value $g_0$ can be non-null:

$$y_e(x) = y_0 + \sum_{n=1}^{\infty} A_n (g - g_0)_n \quad , \text{ where } \quad y_0 = g_0 \sum_{n=1}^{\infty} A_n \tag{4.10}$$

We also note that in the interval $[-L, L]$, each of equations (4.8) and (4.9) is obtained by the linear superposition of the following subequations:

$$y'_n = \frac{2L}{\pi} \frac{(-1)^{n+1}}{n} \sin \frac{n\pi x}{L}, \text{ respectively } y'_0 = \frac{L}{2} J_1 \text{ and } y'_n = \frac{L}{2} \frac{1}{2^n} J_{2n}, \text{ for } n=1, 2, ...,\infty$$

By solving these partial equations and adding the solutions, we get the same results as those obtained by integrating the whole equation. In the present case, the derivatives $y'_n(x)$ of the partial solutions $y_n(x)$ are even the harmonics (or quasi-harmonics) of the expansion in a known series of the function $y'$.

Consequently, the method can be applied to all equations in which the expansion in series of one term leads to a linear equation represented as an infinite superposition of partial equations for which the solution can be found by classical solving methods.

### 4.2. Linearization of the differential equation of the simple gravity pendulum, by introducing an infinite sum of ramp functions

We will apply the method of linearization of nonlinear equations to the equation of the **simple gravity pendulum**:

$$\frac{d^2\theta}{dt^2} + \omega_0^2 \sin\theta = 0 \quad , \quad \theta \in [-\pi, \pi] \tag{4.11}$$

In a classical approach [13], the equation is solved by the method of approximation of the equation, replacing $sin\theta \approx \theta$, acceptable for $\theta \in [-\theta_1, \theta_1]$, if $\theta_1 \approx 0$. Here, we will try to find an exact solution, valid for the whole interval. For this, we will use the non-sinusoidal series expansion of the term $sin\theta$, the basis of the expansion being the following odd function:
$g_1(\theta) = X_o[-\pi < -\theta - \pi > -\pi/2 < \theta > \pi/2 < -\theta + \pi > \pi]$.
The Fourier series expansion of the function g1 ($\theta$) leads us to:

$$\overline{g}_1(\theta) = \frac{4}{\pi} \sum_{n=1}^{\infty} \frac{(-1)^{n+1} \sin(2n-1)\theta}{(2n-1)^2} \qquad \text{whose coefficients are:}$$

$d_1=4/\pi$, $d_2=0$, $d_3=-4/9\pi$, $d_4=0$, $d_5=4/25\pi$, $d_6=0$, $d_7=-4/49\pi$, $d_8=0$, $d_9=4/81\pi$, $d_{10}=0$, $d_{11}=-4/121\pi$, $d_{12}=0$, ... For the function $f(\theta)=sin\theta$, we can write:

$$\sin\theta = \hat{f}(\theta) = \sum_{n=1}^{\infty} B_n g_n(\theta), \text{ in which:}$$

$$B_1 = \frac{\pi}{4}, B_2 = 0, B_3 = \frac{\pi}{36}, B_4 = 0, B_5 = -\frac{\pi}{100}, B_6 = 0, B_7 = \frac{\pi}{196}, B_8 = 0, B_9 = 0, ...$$

$$B_{2n-1} = \frac{\pi}{4} \cdot \frac{(-1)^n}{(2n-1)^2}, B_{2n} = 0, \text{ but } B_{n^2} = 0, \text{ for } n=2, 3, .., \infty. \text{ We get the linear equation:}$$

$$\frac{d^2\theta}{dt^2} + \omega_0^2 \sum_{n=1}^{\infty} B_n g_n(\theta) = 0$$

which, for $\theta = \sum_{n=1}^{\infty} \theta_n$, is the linear superposition of partial equations:

$$\frac{d^2\theta_n}{dt^2} + \omega_0^2 B_n g_n(\theta_n) = 0, \text{ où } n=1, 2, 3, ...,\infty \tag{4.12}$$

If $\theta(0)=\theta_0$, and $v_0=0$, these initial conditions are valid for all the partial equations:

$\theta_n(0)=\theta_0$, et $v_{0n}=0$.

If $v_0 \neq 0$, the given problem is replaced by an equivalent problem: we consider the initial angular speed $v_0$ imposed on the pendulum, as coming from the transformation into kinetic energy $E_c = 1/2 \cdot mv^2$ of the part $E_p = mg(h_e - h_0)$ of the potential energy of the pendulum of length $l$, located at an equivalent height $h_e$ (which corresponds to an equivalent angle $\theta_e$). Thus, the initial problem is replaced by that in which the initial angular speed $v_0$ is zero and the initial position is:

$$\theta_e = 2k\pi \pm \arccos\left(\cos\theta_0 - \frac{lv_0^2}{2g}\right) \tag{4.13}$$

Values $k \neq 0$ appear when $h_e > 2l$. For these cases, we can calculate $v_\pi$, the speed at which the pendulum crosses the position $h = \pm\pi$. If $\theta_e$ and $\theta_0$ belong to the interval $(-\pi, \pi)$, the equivalent partial solutions $\theta_n(t)$ will all have the initial conditions:

$$\theta_n(0) = \theta_e \text{ and } \frac{d\theta_n}{dt}\bigg|_{t=0} = 0 \text{ , and for } h_e > 2l: \ \theta_n(0) = \pi \text{ and } \frac{d\theta_n}{dt}\bigg|_{t=0} = v_\pi \tag{4.14}$$

The resolution of these $n$ partial equations leads us to find the $n$ partial responses of the system to these $n$ periodic forces acting on the pendulum. The sum of these responses is the solution to the equation for the equivalent problem. This solution-trajectory $\theta(t)=\Sigma\theta_n(t)$ passes through the position $\theta_0$ at time $t_e$. The initial condition for the initial problem becomes $\theta_0=\theta(t-t_e)=\Sigma\theta_n(t-t_e)$. The positions $\theta_n(t_e)$ and the speeds $\theta_n'(t_e)$ of the particular trajectories are achieved by the equivalent pendulum when it passes through the position $\theta_0$. In this way, the equivalent potential energy of the pendulum (corresponding to the initial speed $v_0$) is distributed over the partial components of the total force system and transformed by each equivalent pendulum according to its own specificity, into kinetic energy, resulting in an initial angular speed equivalent for each partial equation.

Because of the laws of conservation of energy, the trajectory of the pendulum, for each partial system of forces $F_n$, is symmetrical compared to the equilibrium position $\theta_{mn}$ $(\theta_0)$ (for $n=1$). It passes through the positions $\theta_m - \theta_i$ with the same speeds $v_i$ with which it passes through the positions $\theta_m + \theta_i$ and reaches the maximum position $-\theta_e$ also at zero speed. If $\theta_e > \pi$, the equivalent speed $v_\pi$ that corresponds to the position $\theta = \pi$ is equal to that of the position $\theta = -\pi$. For the same reasons, the partial trajectories with initial velocities $\theta_e$ negative are symmetrical with respect to the axis $\theta = 0$ with those whose initial velocity $\theta_e$ is positive.

From (4.12), for $n=1$, we obtain the partial equation

$$\frac{d^2\theta_1}{dt^2} + \omega_0^2 \frac{\pi}{4} g_1(\theta) = 0 \text{ , where } g_1(\theta)=X_o[-\pi<-\theta-\pi>-\pi/2<\theta>\pi/2<-\theta+\pi>\pi],$$

which is resolved successively, for each subinterval in which the force exerted on the pendulum has a certain law of variation. For $-\pi/2 < \theta_e < \pi/2$, on the subinterval $[-\pi/2$, $\pi/2]$ the force increases uniformly from a negative value to a positive value, the equation is harmonic and the solution is :

$$\theta_1(t) = \theta_e \cos\left(\sqrt{\pi/4}\omega_0 t\right) \tag{4.15}$$

For $-\pi < \theta_e < -\pi/2$, the partial equation takes different forms for different subintervals. On the subinterval $[\theta_e, -\pi/2]$ the equation takes the form:

$$\frac{d^2\theta_1}{dt^2} - \omega_0^2 \frac{\pi}{4} \cdot (\theta_1 + \pi) = 0 \text{ ,} \tag{4.16}$$

with the initial conditions: $\theta_{1.0}=\theta_1(0) = \theta_e$ and $v_{1.0} = \dfrac{d\theta_1}{dt}\bigg|_{t=0} = 0$ and with the solution:

$$\theta_1(t) = (\theta_e + \pi)ch\left(\sqrt{\pi/4} \cdot \omega_0 t\right) - \pi \tag{4.17}$$

$$v_1(t) = (\theta_e + \pi)\omega_0\sqrt{\pi/4}\,sh\left(\sqrt{\pi/4}\cdot\omega_0 t\right) \qquad (4.18)$$

The equations (4.16) and (4.18) allow us to find the moment $t_{1.1}$ of switching of the force field, i.e. the moment when $\theta_I = -\pi/2$, then the angular velocity of the pendulum from this moment: $v_{1.1} = v_1(t_1) = (\theta_e + \pi)\omega_0\sqrt{\pi/4}\,sh\left(\sqrt{\pi/4}\cdot\omega_0 t_1\right)$. These values become the initial conditions of the harmonic equation in the subinterval $[-\pi/2 \, , \, \pi/2]$, available for $t > t_{1.1}$. The solution is:

$$\theta_1(t - t_{1.1}) = (\pi/2)\cos\left(\sqrt{\pi/4}\,\omega_0(t - t_{1.1})\right) + (v_{1.1}/\omega_0\sqrt{\pi/4})\sin\left(\sqrt{\pi/4}\,\omega_0(t - t_{1.1})\right) \qquad (4.19)$$

Now, from $\theta_I(t) = 0$, we can calculate the time $t_{1.2}$ of the passage through the equilibrium point.

Further on, the evolution of the pendulum is symmetrical with respect to this point. The pendulum reaches the position $\theta_I(t) = -\theta_e$ after the time $2t_{1.2}$ and will continue with a symmetrical trajectory compared to the axis $t = 2t_{1.2}$, to arrive after a total time $t_t = 4t_{1.2}$ again in the position $\theta_I(t) = \theta_e$ and to continue a periodic trajectory with the period $T = 4t_{1.2}$.

Figure 2.a. graphically shows the solutions for the first harmonic in the case of zero initial velocities, for initial positions greater than and less than $-\pi/2$, and Figure 2.b., for non-zero initial velocities.

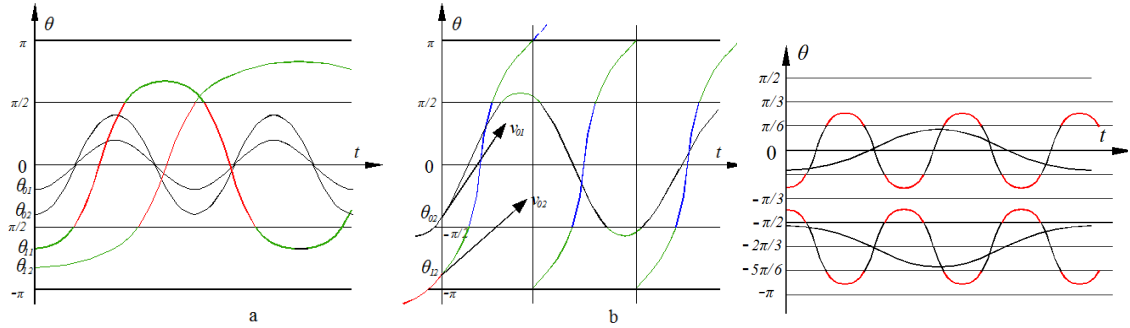

*Fig.2. **a**: first order partial solutions for $v_0 = 0$, $-\pi/2 < \theta_{01}$, $\theta_{02} < 0$ and $-\pi < \theta_{11}$, $\theta_{12} < -\pi/2$.*
*  **b**: for $v_0 \neq 0$  **c**: third order partial solutions for $v_0 = 0$ şi $\theta_e < 0$*

For $n = 3$ we have the equation $\quad \dfrac{d^2\theta_3}{dt^2} + \omega_0^2 \dfrac{\pi}{4}\dfrac{1}{9}g_3(\theta) = 0, \quad$ where

$g_3(\theta) = X_0[-\pi/3 < 3(-\theta - \pi/3) > -\pi/6 < 3\theta > \pi/6 < 3(-\theta + \pi/3) > \pi/3]_3 = X_0[-\pi < -3\theta - 3\pi > -5\pi/6 < 3\theta + 2\pi > -\pi/2 < -3\theta - \pi > -\pi/6 < 3\theta > \pi/6 < -3\theta + \pi > \pi/2 < 3\theta - 2\pi > 5\pi/6 < -3\theta + 3\pi > \pi]$

The initial conditions of the equation are given by (4.14).

For $v_0 = 0$, if $-\pi/6 < \theta_0 < \pi/6 \rightarrow g_3(\theta) = 3\theta$, and the solution is: $\theta_3(t) = \theta_0\cos\left(\sqrt{\pi/12}\,\omega_0 t\right)$,

For $v_0 \neq 0$, if $|\theta_e| = \sqrt{\theta_0^2 + \left(v_0/\omega_0\sqrt{\pi/12}\right)^2} \leq \pi/6$, so

$$\theta_3(t) = \theta_0\cos\left(\sqrt{\dfrac{\pi}{12}}\omega_0 t\right) + \dfrac{v_0}{\omega_0\sqrt{\pi/12}}\sin\left(\sqrt{\dfrac{\pi}{12}}\cdot\omega_0 t\right) = \theta_e\cos\left[\sqrt{\dfrac{\pi}{12}}\cdot\omega_0 t + arctg\left(-v_0/\theta_0\omega_0\sqrt{\dfrac{\pi}{12}}\right)\right]$$

For the other subintervals in which $\theta_e$ can be located, the partial solution is:

$\theta_3(t) = (\theta_e - \theta_m)\cos\left(\sqrt{\dfrac{\pi}{12}}\omega_0 t\right) - \theta_m,\quad$ or $\quad \theta_3(t) = (\theta_e - \theta_m)ch\left(\sqrt{\dfrac{\pi}{12}}\omega_0 t\right) - \theta_m,\quad$ as the slope of the function $g_3(\theta)$ is positive, respectively negative. Here, $\theta_m$ is the midpoint of this subinterval and it always has the same sign as $\theta_e$.

Consequently, the cosine $\theta_3(t)$ (as well as those of higher order) is an oscillating movement with respect to the point $\theta_m$, with the amplitude $\theta_e - \theta_m$ (it never arrives at the point of equilibrium $\theta = 0$ nor in $\theta_0$, if it is outside the subinterval). The hyperbolic cosine $\theta_3(t)$ (as well as the higher-order quasi-cosine) describes a divergent movement, moving away from $\theta_m$ (which is an unstable equilibrium point).

Figure 2.c presents some of these partial solutions, for different negative values of $\theta_e$ ($v_0 = 0$).

For n = 5, the partial equation is $\dfrac{d^2\theta_5}{dt^2} - \omega_0^2 \dfrac{\pi}{4} \dfrac{1}{25} g_5(\theta) = 0$, where

$g_5(\theta) = X_o[-\pi/5 < 5(-\theta-\pi) > -\pi/10 < 5\theta > \pi/10 < 5(-\theta+\pi) > \pi/5]_5$

It is solved in the same way as the previous equations. It should be noted that due to the negative coefficient of this quasi-harmonic, the position $\theta = 0$ is an unstable equilibrium. For $-\pi/10 < \theta_e < \pi/10$, the trajectory of the pendulum is divergent.

The general solution of the equivalent equation is the sum of all the partial solutions: $\theta(t) = \Sigma \theta n(t)$. It is also the exact solution of the given equation, expressed as an infinite sum of partial solutions, each of these solutions having different expressions on different subintervals. For this solution to be practically useful, we will only retain the first $N$ partial solutions (the approximation of the solution). The same result is achieved by approximating the method, if from the expansion in SFN of the function $\sin\theta$ we retain the first $N$ terms.

For a certain value $t = t_e$, the sum $\theta(t) = \Sigma \theta_n(t)$ is zero: $\Sigma \theta_n(t_e) = 0$. As shown in Fig. 3 (for simplicity we have chosen the case $\theta_0 = 0$), at this time, each partial solution has a value $\theta_n(t_e) = \theta_{n0}$ and a speed $v_n(t_e) = v_{n0}$. These are the initial conditions for the components of the initial equation.

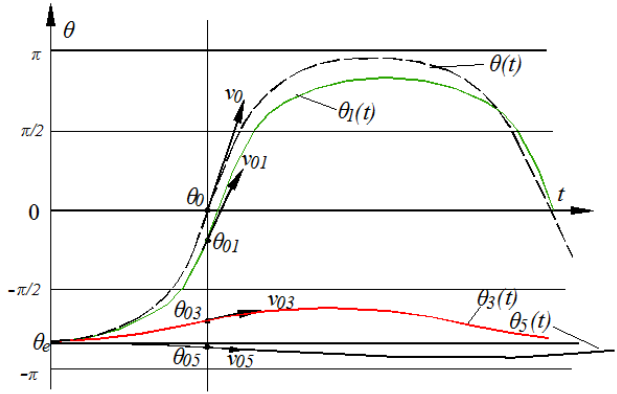

*Fig.3. The general solution and the partial solutions of order 1, 3 and 5 for $v_0 \neq 0$*

### 4.3. Linearization of the differential equation of the gravitational pendulum, by introducing an infinite sum of rectangular functions

The expansion in non-sinusoidal Fourier series that we chose for the linearization of equation (4.11) is the closest to the traditional approach, but it is only one of the many possibilities offered by the method that we propose. Another interesting solution is the replacement in the equation of the function $\sin\theta$ with its expansion in infinite series of periodic rectangular pulses:

$$\frac{d^2\theta}{dt^2} + \omega_0^2 \sum_{n=1}^{\infty} B_n \Pi_n(\theta) = 0 \;, \tag{4.20}$$

where the variable $\theta$ can have values in the interval $[-\pi, \pi]$,
$\Pi_1(\theta) = \Pi[-\pi(-1)0(1)\pi]$, et $\Pi_n(\theta) = \Pi[-\pi/n(-1)0(1)\pi/n]_n$.

As in the previous approach, the initial conditions of the equivalent problem are:

$$\theta_0 = \theta(0) = \theta_e \text{ and } v_0 = \frac{d\theta}{dt}\Big|_{t=0} = 0 \ , \tag{4.21}$$

where $\theta_e$ is the equivalent angle, corresponding to the equivalent height $h_e$.

Let be $B_n$ the expansion coefficients of $\sin\theta = \hat{f}(\theta) = \sum_{n=1}^{\infty} B_n \Pi_n(\theta)$

Because $\Pi_1 = \frac{4}{\pi} \sum_{n=1}^{\infty} \frac{\sin(2n-1)\theta}{2n-1}$, the expansion coefficients of the function $sin\theta$ are:

$$B_1 = \frac{\pi}{4}, B_2 = 0, B_3 = -\frac{\pi}{12}, \ B_4 = 0, B_5 = -\frac{\pi}{20}, B_6 = 0, B_7 = -\frac{\pi}{28}, B_8 = 0, \ B_9 = 0, B_{10} = 0, \ ...$$

$$B_{2n-1} = \frac{\pi}{4} \cdot \frac{-1}{2n-1} \ , \ \ B_{2n} = 0, \text{ but } \ B_{2^2} = 0, \qquad \text{for } n=2, 3, .., \infty$$

All the coefficients, with the exception of that of the fundamental quasi-harmonics, are negative. The linear equation (4.20) is the superposition of an infinite number of equations of the type:

$$\frac{d^2\theta_n}{dt^2} = -\omega_0^2 B_n \Pi_n(\theta) \ , \text{ où } n=1, 2, 3, ...,\infty, \text{ with the initial conditions (4.21). The solutions of}$$

these equations, for all subintervals $[t_k, \ t_{k+1}]$ in which $\Pi_n(t)=ct$ are:

$$\theta_n(t) = -\omega_0\theta^2 B_n sgn\Pi_{nk}(t-t_k)^2/2 + v_{0k}(t-t_k) + \theta_{0k}, \text{ for } k= 0, 1, 2, 3, ..., n, ... \tag{4.22}$$

where $sgn\Pi_{nk}$ is the sign of the functions $\Pi_n(\theta)$ for the subinterval $[\theta_k, \ \theta_{k+1}]$, $t_k$ are the moments of the commutations of the function $\Pi_n(\theta)$, $\theta_{0k}$ and $v_{0k}$ are the initial conditions for the available equation on this subinterval (the final values of the previous subinterval). The pendulum speed for each system of forces is uniformly accelerated:

$$v_n(t) = -sgn(\Pi_n)B_n\omega_0^2(t-t_k) + v_{0k}$$

For $n=1$, the point $\theta=0$ is a point of stable equilibrium. Whatever the initial position of the pendulum, it tends to reach a stable equilibrium position and oscillates around this position, describing a quasi-sinusoid built on the basis of a polynomial of the second degree. This curve will have the inflection points located on the axis $\theta=0$. The oscillations will have the amplitude $\theta_e$ and the period of oscillation $T_1 = 4\sqrt{8|\theta_e|/\pi\omega_0^2}$ (dependent on $\theta_e$). A pendulum that starts from the position $\theta_{0I}=\theta_e$ with the speed $v_{0I}=0$, arrives after the time $t=T_I/4$ in the position $\theta_I=0$ (where the force which acts on it changes direction) with the speed $v_1 = -\omega_0\sqrt{|\theta_e|\pi/2}$. Due to the inertia, the pendulum continues to move towards the position $\theta_I=-\theta_e$, where it arrives after the temp $t=T_I/4$, with the speed $v_I=0$. Under the action of the same system of forces, the pendulum continues its movement in the opposite direction and after another quarter of period, it finds the position $\theta_I=0$, this time with speed $v_1 = \omega_0\sqrt{|\theta_e|\pi/2}$. After a new change in direction of the force, after another quarter of a period, the pendulum returns to its position $\theta_I=\theta_e$, with the speed $v_I=0$. The oscillations continue with the $T_I$ period. With the notations $A=\pi\omega^2/8$ and $T = \sqrt{8|\theta_e|/\pi\omega_0^2}$, we have:

$$\theta_I(t) = X_I^2 = X^2[0*(\theta_e - At^2)*T*(-\theta_e + A(t-T)^2)*3T*(\theta_e - A(t-3T)^2)*4T]_1 \tag{4.23}$$

If the initial conditions change: $\theta_0=\theta_e$ and $v_I(0)\neq0$, the solution of the equation changes its period, its amplitude and its initial phase shift.

The partial solution of order $n$ of equation (4.20) at the initial conditions (4.21) is also a quasi-sinusoidal polynomial of second degree: (4.24)

$$\theta_n(t) = X_n^2 = X^2[0*(-\theta_n + At^2/n)*Tn^{1/2}*(\theta_n - A(t-T/n)^2)*3 \ n^{1/2}*(-\theta_n + A(t-3T/n)^2)*4T \ n^{1/2}]_n,$$

where $\theta_n = \theta_e - \theta_m$. This curve has all the inflection points located on the axis $\theta = \theta_m$. The oscillation has amplitude $\theta_e - \theta_m$ and the oscillation period $T_n = 4\sqrt{2|\theta_e - \theta_e|/B_n\omega_0^2}$.

In figure 4.a, we have represented, for the fundamental harmonic, three of these partial solutions, for different values of the equivalent initial position $\theta_e$: two for $\theta_e \in [-\pi, \pi]$ (red line) and one for $\theta_e \notin [-\pi, \pi]$ (dotted black line); with a dotted green line, we have represented the fictitious continuation of the quasi-sinusoid outside the interval $[-\pi, \pi]$).

The general solution of the equivalent equation is the sum of all the partial solutions: $\theta(t) = \Sigma\theta_n(t)$. An approximation of the solution is obtained by adding the first $N$ partial solutions, where $N$ must be large enough to obtain a satisfactory error. For a value $t = t_e$, the sum is zero: $\Sigma\theta_n(t_e) = 0$. As shown in Figure 4b, for the moment $t_e$, each partial solution provides a value $\theta_n(t_e) = \theta_{n0}$ and a speed $v_n(t_e) = v_{n0}$. These are the initial conditions of the components of the original equation. The solution of the equation is $\theta(t) = \Sigma\theta_n(t - t_e)$.

With the notations of (4.23) and (4.24) we can advance a new expression for the movement of the pendulum:
$\boldsymbol{\theta(t) = \Sigma X n^2}$, pour $t > t_e$.

*Fig.4. **a**: Partial solutions for $\theta_e < \pi$: order 1 (red line), order 3 (black line). For $\theta_e > \pi$: order 1 (dotted line), rang 3 (blu line) b: solution générale*

## 5. Conclusions

We have tried here, using a few simple examples, to prove that the method of solving differential equations by determining the coefficients of expansion in sinusoidal or non-sinusoidal series of the unknown function is a particularly solid method, applicable to all types of differential and integro-differential equations, linear and nonlinear, partial differential equations, systems of such equations, whatever their order and whatever the complexity of the coefficients.

## 6. Bibliographie

[1] Török A., Petrescu S., Feidt M., Séries de Fourier périodiques non sinusoïdales, https://hal.archives-ouvertes.fr/hal-02485085

[2] Nagle R. K., Saff E. B., Snider A. D., Fundamentals of Differential Equations and Boundary Value Problems, 8th ed, 2012, Pearson Education, Inc.

[3] O'Neil P. V., Advanced Engineering Mathematics, Seventh Edition, part 1, Cengage Learning, 2012, Publisher: Global Engineering: Christopher M. Shortt

[4]    Kreyszig E., Advanced Engineering Mathematics, 10-th edition, chapters 1, 2, 5, 6, 11, 12, John Wiley & Sons, Inc., 2011

[5]    Polyanin A. D., A. V. Manzhirov, Handbook of Mathematics for Engineers and Scientists, Chapman & Hall/CRC Press, 2007

[6]    Polyanin A. D., A. V. Manzhirov, Handbook of Exact Solutions for Ordinary Differential Equations, CRC Press, 2003

[7]    Zwillinger D., CRC Standard Mathematical Tables and Formulae, 31st Edition, 2003, by CRC Press Company

[8]    Spiegel M. R., Mathematical Handbook of Formulas and Tables, Schaum's Outline Series, Schaum Publishing Co. 1st.ed. 1968, McGraw-Hill Book Company
https://archive.org/details/MathematicalHandbookOfFormulasAndTables

[9]    Dourmashkin, Classical Mechanics: MIT 8.01 Course Notes, Chapter 23 Simple Harmonic Motion, pp 8,
http://web.mit.edu/8.01t/www/materials/modules/chapter23.pdf

[10]   Jordan D. W., Smith P., Nonlinear Ordinary Differential Equations. An introduction for Scientists and Engineers,4-th edition, Keele University, Oxford University Press, 2007,

[11]   Olver P. J., Equations Introduction to Partial Differential Equations, Springer Cham Heidelberg New York Dordrecht London, Springer International Publishing Switzerland, 2014, 6056 ISSN - 2197 5604 (electronic)

[12]   Dourmashkin, Classical Mechanics: MIT 8.01 Course Notes, Chapter 23 Simple Harmonic Motion, pp 8,
http://web.mit.edu/8.01t/www/materials/modules/chapter23.pdf

[13]   Tolstov G. P., Silverman R. A., Fourier Series, Courier Corporation, 1976

[14]   Al-Gwaiz M.A., Sturm-Liouville Theory and its Applications, Springer Undergraduate Mathematics Series ISSN 1615-2085, Springer-Verlag London Limited 2008

# LE DÉVELOPPEMENT EN SÉRIES PÉRIODIQUES,
# METHODE DE RESOLUTION D'EQUATIONS DIFFERENTIELLES

Arpad Török[1], Stoian Petrescu[2], Michel Feidt[3]

[1] PhD student, The Polytechnic University of Bucharest, Department of Engineering Thermodynamics, 313, Splaiul Independentei, 060042 Bucharest, Romania, e-mail: arpi_torok@yahoo.com
[2] Prof. Dr. Eng., Polytechnic University of Bucharest, Department of Engineering Thermodynamics, România
[3] Prof. Dr. Eng., L.E.M.T.A., U.R.A. C.N.R.S. 7563, Université de Loraine Nancy 12, avenue de la Foret de Haye, 54516 Vandœuvre-lès-Nancy, France

**Résumé:** Le développement des fonctions de variables réelles en séries de Taylor et de Frobenius (séries entières lesquels sont constituées dans des bases nonorthogonales, nonpériodiques), en séries de Fourier sinusoïdales (des bases des fonctions orthogonales, périodiques), en séries de fonctions spéciales (des bases des fonctions orthogonales, nonpériodiques), etc est une procédé couramment utilisé pour résoudre une large gamme d'équations différentielles ordinaires (ODEs) et d'équations aux dérivées partielles (PDEs). Dans cet article, basé sur une analyse approfondie des propriétés des séries de Fourier périodiques sinusoïdales (SFS), nous serons en mesure d'appliquer cette procédure à une catégorie beaucoup élargie d'ODEs (toutes les équations linéaires, homogènes et non homogènes à coefficients constants, une large catégorie d'équations linéaires et non linéaires à coefficients variables, systèmes d'ODEs, équations intégro-différentielles, etc.). Nous allons également étendre cette procédure et l'utiliser pour résoudre certains ODEs, sur des bases périodiques non orthogonales, représentées par des séries de Fourier périodiques non sinusoïdales (SFN).

**Mots clefs**: séries de Fourier sinusoïdales, séries de Fourier non sinusoïdales, bases orthogonales, équations différentielles, approximation des solutions, le pendule gravitationnel

## 1. Introduction

Les méthodes de résolution des ODE, proposées ici, utilisent des résultats récemment obtenus [1] dans le domaine de l'analyse fonctionnelle, concernant le développement de fonctions variables réelles, définies sur un intervalle $[-L, L]$, en séries infinies de fonctions périodiques sur le même intervalle, formant des bases orthogonales, mais aussi non orthogonales). Selon l'analyse harmonique (Fourier), toute fonction $f(x)$, périodique sur l'intervalle $[-L, L]$, qui satisfait aux conditions de Dirichlet, peut être développée en une somme infinie, connue dans la littérature sous le nom de série trigonométrique (pour laquelle, pour des raisons mises en évidence dans l'article [1], nous avons utilisé le nom de série sinusoïdale). Cette série est formée par les composantes d'une base biortogonale complète, composée par la fonction unité $1$, les harmoniques fondamentales $sin(\pi x/L)$-paire et $cos(\pi x/L)$-impaire, avec la période $2L$ et les harmoniques secondaires $sin(n\pi x/L)$ et $cos(n\pi x/L)$, avec la période $2L/n$, pour $n \in \mathbf{N}^+$. Les coefficients de cette développement (coefficients de Fourier) peuvent être calculés à l'aide d'intégrales définies (formules d'Euler). Ce papier généralise cette affirmation en montrant que la fonction $f(x)$ peut également être développée en séries périodiques non sinusoïdales (SFN), consistant en la somme infinie des composantes pondérées d'une base *complète*, non orthogonale: la fonction unitaire $1$, les quasi-harmoniques fondamentales $g(x)$-paire et $h(x)$-impaire, périodiques avec la période $2L$, avec valeur moyenne nulle sur l'intervalle de définition et les quasi-harmoniques secondaires, définies sur $[-L, L]$, $g_n(x)=g(nx)$ et $h_n(x)=h(nx)$, avec la période $2L/n$, pour $n \in \mathbf{N}^+$. Les quasi-harmoniques fondamentales $g(x)$ et $h(x)$ peuvent être toute fonction réelle, de variable réelle, qui admette sur l'intervalle $[-L, L]$ des développements en série sinusoïdale. On obtient les coefficients $A_n$ et $B_n$ du développement en SFN de la fonction $f(x)$, à l'aide des relations algébriques entre les coefficients du développement en SFS des fonctions $g(x)$ et $h(x)$.

Ainsi, toute fonction $f(x):[x_1, x_2] \rightarrow \mathbf{R}$ $T$-périodique ($T=x_2-x_1$), de l'espace $L^2$ (c'est-à-dire, de carré intégrable), peut être représentée par la somme :

$$\hat{f}(x) = f_0 + \sum_{n=1}^{\infty} A_n [\bar{g}_n(x) - g_0] + \sum_{n=1}^{\infty} B_n \bar{h}_n(x), \text{ où :}$$

$f_0 = \dfrac{1}{T} \int_{x_1}^{x_2} f(x)dx$ et $g_0 = \dfrac{1}{T} \int_{x_1}^{x_2} g(x)dx$ sont les valeurs moyennes sur l'intervalle $[x_1, x_2]$,

$\bar{g}_n(x)$ et $\bar{h}_n(x)$ ($n=1, 2, 3, ..., \infty$) sont des **Fourier-fonctions** (des fonctions continues, justement les développements en SFS des fonctions $g(nx)$ et $h(nx)$, definis sur $[-T/n, T/n]$.

À partir de ce résultat général, dans le travail [1] sont décrits aussi quelques situations particulières :

- $\hat{f}(x) = f_0 + \sum_{n=-\infty}^{\infty} C_n [\bar{g}_n(x) - g_0]$, où $g(x)$ a tous les deux composants: le pair et l'impaire

- $\hat{f}(x) = f_0 + \sum_{n=1}^{\infty} A_n \mathbf{S}[g(x)]_L + \sum_{n=1}^{\infty} B_n \mathbf{C}[g(x+L/2)]_L$ où la fonction $g(x)$ est une fonction

de $L^2$ définie sur l'intervalle $[0, L/2]$, $\mathbf{S}[g(x)]_L$ et $\mathbf{C}[g(x+L/2)]_L$ sont les fonctions dénommées quasi−sinusoïdes, dérivés d'elle

- $\tilde{f}(x) = f_0 + \sum_{n=1}^{\infty} [A_n^0 \Phi_n(x) + B_n^0 \Psi_n(x)]$ où $\Phi_n(x)$ et $\Psi_n(x)$ sont des fonctions orthogonales, générées par les Fourier-fonctions $\bar{g}(x) - g_0$ et $\bar{h}(x)$, par un procédé d'orthogonalisation.

## 1.1. Propriétés des développements en séries de Fourier périodiques (le cas général)

Le développement en SFN d'une fonction $f(x)$ est obtenu de son développement en SFS, par une redistribution et un regroupement des termes du développement, de manière à obtenir les développements en SFS des fonctions périodiques $g(x)$, $h(x)$ et des quasi-harmoniques $g_n(x)$, $h_n(x)$ correspondants. Par conséquent, ce développement bénéficie également des propriétés de convergence, de dérivabilité et d'intégrabilité, similaires à celles de SFS [1-9]. Donc :

- soit $f(x)$ une fonction $2L$-périodique, continue dans l'intervalle $[-L, L]$. Son

développement de Fourier $f(x) = f_0 + \sum_{n=1}^{\infty} A_n g_n(x) + \sum_{n=1}^{\infty} B_n h_n(x)$, où $g_0=0$, sinusoïdal ou non,

convergent ou non, peut être intégré terme par terme, entre toutes limites d'intégration:

$$F(x) = \int_0^x f(x)dx = C + f_0 x + \sum_{n=1}^{\infty} A_n G_n(x)dx + \sum_{n=1}^{\infty} B_n H_n(x)dx, \qquad (1.1)$$

où $G(x)$ et $H(x)$ sont les primitives des $g(x)$, respectivement $h(x)$, et $C$ est une constante d'intégration lequel dépende des coefficients du développement. Après les remplacements

$x = \sum_{n=1}^{\infty} C_n G_n(x)$ (les fonctions $x$ et $G(x)$ sont impaires) et $C = \sum_{n=1}^{\infty} B_n H_n'(0) + F(0)$, résulte un développement en SFN de la primitive $F(x)$, en une base, le plus souvent différente de celle de la fonction $f(x)$.

- soit $f(x)$ une fonction $2L$-périodique, continue dans l'intervalle $[-L, L]$, avec $f(-L)=f(L)$ et avec la dérivée $f'(x)$ lisse par portions dans cet intervalle. Le développement de Fourier, sinusoïdal ou non, de la fonction $f'(x)$, peut être obtenu en dérivant terme par terme le développement de Fourier de la fonction $f(x)$. La série obtenue converge ponctuellement vers $f'(x)$ en tous les points de continuité et vers $[f'(x) + f'(-x)]/2$ en ceux de discontinuité.

$$\text{Si } f(x) = f_0 + \sum_{n=1}^{\infty} A_n g_n(x) + \sum_{n=1}^{\infty} B_n h_n(x) \text{ , alors: } \hat{f}'(x) = \sum_{n=1}^{\infty} A_n g_n'(x) + \sum_{n=1}^{\infty} B_n h_n'(x) \qquad (1.2)$$

Dans ce cas aussi, la base du développement en série non sinusoïdal de la dérivée diffère de celle de la fonction $f(x)$.

La condition $f(-L)=f(L)$ fait que le nombre des problèmes dans lesquels la formule (1.2) peut être utile est assez petit, mais elle peut être évitée si le saut à partir du point $x=L$ (ainsi que tout autre saut du composant impair) est compensé par un saut dans la direction opposée (obtenu en soustrayant une autre fonction impaire qui fait un saut identique dans le même point). Dans l'exemple suivant, par ce procédé, la composante impaire $f_o$ de la fonction $f(x)$ est décomposée en une somme de la fonction différenciable $f_{os}$ (pour laquelle $f(-L)=f(L)$) et la fonction de rampe $f_r=x \cdot f_o(L)/L \rightarrow f_{os}=f_o-f_r$. Par conséquent:

$$\frac{d}{dx} f_o(x) = \frac{d}{dx}\left[ f_{os}(x) + \frac{f_o(L)}{L} x \right] = \frac{d}{dx} f_{os}(x) + \frac{f_o(L)}{L},$$

relation qui nous permet de trouver une expression pour le développement de la dérivée $f'(x)$ de la fonction $f(x)$, pour toutes les catégories des fonctions qui satisfont les autres conditions. Le développement en SFN de la fonction $f_{os}=f_o-f_r$ nous permet le calcul des coefficients du développement en série de la dérivée.

Le cas général, celui du développement des fonctions dans des bases périodiques non sinusoïdales, met en évidence le fait que l'élément $I=1$ de la base a un caractère particulier. Il fait partie de toutes les bases périodiques, il ne change pas lorsque les autres composantes de la base (de valeur moyenne nulle sur l'intervalle de définition) changent après l'intégration, ou après la dérivation. Son coefficient est calculé par une intégrale définie et non par des relations algébriques. $I=1$ est une fonction paire, mais pour $f_0=0$, il reçoit simultanément un caractère impair aussi. Les dérivées et les primitives de toutes les fonctions paires (y compris la fonction $I \cdot f_0=f_0$) sont des fonctions impaires et, inversement, celles des fonctions impaires (y compris la fonction $I \cdot 0=0$) sont des fonctions paires. En dérivant toute fonction paire $f_0$ on obtient la fonction impaire $0$ et vice versa, en intégrant la fonction impaire $0$, on obtient toute fonction paire $f_0$ (la constante d'intégration)

## 1.2. Propriétés des développements en séries de Fourier sinusoïdales

Comme dans le cas général, dans le cas des développements en SFS, sur l'intervalle $[-L, L]$, pour $\omega_n=n\omega_0=n\pi/L$, parce que $\bar{x}=2\sum_{n=1}^{\infty}\dfrac{(-1)^{n+1}}{\omega_n}\sin(\omega_n x)$, nous pouvons écrire [9]:

$$\text{Pour } \bar{f}(x) = f_0 + \bar{f}_e + \bar{f}_o = f_0 + \sum_{n=1}^{\infty} a_n \cos(\omega_n x) + \sum_{n=1}^{\infty} b_n \sin(\omega_n x), \qquad (1.3)$$

on a: $f_o(L) = -f_o(-L)$, $f_e(L)=f_e(-L)$, $f_e(L)=\dfrac{f(L)+f(-L)}{2}$, $f_o(L)=\dfrac{f(L)-f(-L)}{2}$,

$f(0) = f_0 + \sum_{n=1}^{\infty} a_n$ ,

$$\bar{f}'(x) = \overline{\Phi}(x) = \Phi_0 + \overline{\Phi}_e + \overline{\Phi}_o = \Phi_0 \cdot 1 + \sum_{n=1}^{\infty} \alpha_n \cdot \cos(\omega_n x) + \sum_{n=1}^{\infty} \beta_n \cdot \sin(\omega_n x), \text{ dans lequel:} \qquad (1.4)$$

$$\Phi_0 = \frac{f_o(L)}{L} = \lim_{n \to \infty} 2\left[(-1)^{n+1} \omega_n b_n\right], \quad \alpha_n = \omega_n\left( b_n + 2(-1)^n \frac{f_o(L)}{\omega_n L} \right), \quad \beta_n = -a_n \omega_n, \quad \Phi(0) = \Phi_0 + \sum_{n=1}^{\infty} \alpha_n$$

$$\left(\int_0^x f(x)dx\right)^{SFS} = \overline{F}(x) - F(0) = F_{00} + \overline{F}_e + \overline{F}_o = F_{00} \cdot 1 + \sum_{n=1}^{\infty} A_n \cdot \cos(\omega_n x) + \sum_{n=1}^{\infty} B_n \cdot \sin(\omega_n x), \text{ dans lequel:}$$

$$F_{00} = \sum_{n=1}^{\infty} \frac{b_n}{\omega_n}, \quad A_n = -\frac{b_n}{\omega_n}, \quad B_n = \frac{a_n}{\omega_n} + \frac{2(-1)^{n+1}}{\omega_n} \frac{f_0(L)}{L}, \quad F(0) = F_0 + \sum_{n=1}^{\infty} A_n = F_0 - F_{00} \qquad (1.5)$$

Également sont valables les relations: $\overline{f}(x) = \int_0^x \overline{\Phi}(x) + \dfrac{f_o(L)}{L} + \sum_{n=1}^{\infty} a_n$ , $f_0 = \dfrac{F_o(L)}{L}$ et

$$\int_a^b \overline{f}(x)dx = f_0(b-a) + \sum_{n=1}^{\infty} \frac{a_n(\sin \omega_n b - \sin \omega_n a) - b_n(\cos \omega_n b - \cos \omega_n a)}{n\omega_0} \qquad (1.6)$$

Ainsi, nous avons obtenu une série de relations pour calculer les valeurs moyennes et les coefficients des développement des fonctions dérivées et primitives du premier rang, à partir des valeurs des coefficients du développement en SFS de la fonction *f(x)* et des valeurs de la fonction dans les points limites *f(–L)* et *f(L)*. Ces relations nous permettent également de calculer, pas à pas, les expressions des dérivées et desprimitives de rang supérieur, après avoir calculé les valeurs de ces fonctions à la limite de l'intervalle, les valeurs moyennes et les valeurs des coefficients de leurs développements. Les nouvelles relations peuvent être utilisées pour résoudre les équations différentielles et intégro-différentielles de rang supérieur, pour calculer des intégrales définies ou indéfinies, etc. A titre d'exemple, voici la dérivation et l'intégration du développement de la fonction *f(x)=eˣ* pour l'intervalle [−π, π] [8, 9]:

$$\overline{f}(x) = e^x = \frac{\sinh \pi}{\pi} + \sum_{n=1}^{\infty} \left[ \frac{2\sinh \pi}{\pi} \cdot \frac{(-1)^n}{1+n^2} \cos(nx) - \frac{2n \cdot \sinh \pi}{\pi} \cdot \frac{(-1)^n}{1+n^2} \sin(nx) \right].$$

$$\Phi_0 = \frac{\sinh \pi}{\pi}, \alpha_n = n\left( -\frac{\sinh \pi}{\pi} \cdot \frac{2n(-1)^n}{1+n^2} + 2(-1)^n \frac{\sinh \pi}{n \cdot \pi} \right) = \frac{2\sinh \pi}{\pi} \cdot \frac{(-1)^n}{1+n^2}, \beta_n = -n\frac{2\sinh \pi}{\pi} \cdot \frac{(-1)^n}{1+n^2},$$

$$F_{00} = \sum_{n=1}^{\infty} \frac{b_n}{\omega_n} = \frac{2n \cdot \sinh \pi}{n \cdot \pi} \sum_{n=1}^{\infty} \frac{(-1)^{n+1}}{1+n^2} = \frac{2\sinh \pi}{\pi} \cdot \frac{1}{2}\left(1 - \frac{\pi}{\sinh \pi}\right) = \frac{\sinh \pi}{\pi} - 1, \quad F_0 = \frac{\sinh \pi}{\pi} \qquad (1.7)$$

$$A_n = \frac{-1}{n}\left( -\frac{2n \cdot \sinh \pi}{\pi} \cdot \frac{(-1)^n}{1+n^2} \right), \quad B_n = \frac{1}{n}\left[ \frac{2\sinh \pi}{\pi} \cdot \frac{(-1)^n}{1+n^2} - 2(-1)^n \frac{\sinh \pi}{\pi} \right] = -\frac{2n \cdot \sinh \pi}{\pi} \cdot \frac{(-1)^n}{1+n^2}$$

Donc: $\overline{\Phi}(x) = \left[ \left( e^x \right)' \right]^{SFS} = e^x$ et $\overline{F}(x) = \left( \int_0^x e^x dx \right)^{SFS} + F(0) = e^x$

Pour obtenir (1.7), on start de $\overline{f}(0) = 1 = \dfrac{\sinh \pi}{\pi}\left[ 1 + 2\sum_{n=1}^{\infty} \dfrac{(-1)^n}{1+n^2} \right] \rightarrow \sum_{n=1}^{\infty} \dfrac{(-1)^n}{1+n^2} = \dfrac{1}{2}\left( \dfrac{\pi}{\sinh \pi} - 1 \right)$

On peut le remarquer que pour *f₀=0* et *f₀(L)=0* (sans discontinuités de la composante impaire), à la fois l'intégration et la dérivation se font terme par terme: $\alpha_n = b_n \omega_n$, $\beta_n = -a_n \omega_n$, $A_n = b_n / \omega_n$, $B_n = -a_n / \omega_n$. Si *f₀(x)* a des discontinuités aux extrémités de l'intervalle [−L, L] (ou à l'intérieur), elles provoqueront, lors de la dérivation, l'apparition d'une valeur moyenne $\Phi_0 \neq 0$ et une modification en conséquence du composant $\Phi_e$. Lors de l'intégration, l'effet des discontinuités est transmis au composant impair *Fₒ*. Les dérivés de toutes les fonctions *f(x)+C* ont la même expression, et le retour à la fonction initiale, par l'intégration, est assuré par la relation $F(0) = F_0 + \sum_{n=1}^{\infty} A_n$ , lequel pour *C=0* conduit à *F(0)=0*.

Sur l'axe réel, la fonction périodique impaire *f(x)*, discontinue aux extrémités de l'intervalle [−L, L] (Fig.1A), est la somme entre la fonction continue *fₑ(x)* (Fig.1D) et la fonction «échelle» *f_H(x)* (une succession des « échelons » Heaviside négatifs) (Fig.1E). Les dérivés de ces fonctions sont les fonctions périodiques *f'(x)=*$\Xi_{2\pi}(x+L) \overset{def}{=} \sum_{n=-\infty}^{\infty} \delta_{2n\pi}(x+L)$ (peigne de Dirac, Fig.1B), dont le développement en SFS est $\overline{\Xi} = \frac{1}{\pi}\sum_{n=1}^{\infty}(-1)^n\cos(n\omega_0 x)$ [6], multiplié par le coefficient $f_o(L)/L$ (Fig.1C). Les développements en SFS de ces deux dérivées sont donc divergents, mais leur somme est convergente.

Si la fonction $f(x)$ est paire, la dérivation des sauts à l'intérieur de l'intervalle génère deux peigne de Dirac de signe opposé, qui s'annulent réciproquement, et si $f(x)$ a des discontinuités finies dans un nombre fini des points dans l'intervalle [- $L$, $L$], cela se reflète dans la position et l'amplitude des impulsions Dirac correspondantes.

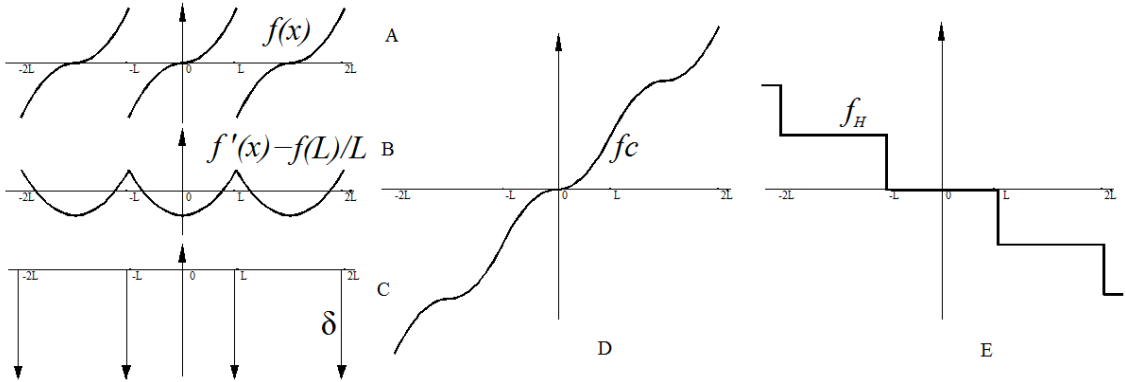

*Fig.1. La dérivation du développement en SFS d'une fonction périodique avec des discontinuités*

### 1.3. Le développement en série de Fourier du produit de deux quelconques fonctions périodiques

Sont connus les relations de calcul pour déterminer les coefficients du développement en SFS du produit $p(x)q(x)$ de deux quelconques fonctions périodiques, de carré intégrable (qui fait que leur produit est aussi une fonction de ce type), définis sur l'intervalle [$-L$, $L$] [9]. Le calcul de ces coefficients est possible si les coefficients de Fourier du développement en série de chacune des deux fonctions sont connus: on doit calculerles les intégrales définies du type $\int_{-L}^{L} p_e \cos\frac{n\pi x}{L} dx$ et $\int_{-L}^{L} p_o \sin\frac{n\pi x}{L} dx$. Si une certaine approximation est permis, le calcul peut être effectué par des méthodes numériques.

Dans le cas des ODEs, un ou tous les deux termes du produit $p(x)q(x)$, peuvent être même l'une des fonctions $y(x)$, $y'(x)$, $y''(x)$, $\int y(x)dx$, etc, avec des coefficients de Fourier inconnus. Le remplacement de ces produits, lorsqu'ils apparaissent dans une équation différentielle, avec leur développement en série de Fourier, ou avec les sommes de Fourier $S_N$ qui les rapproche, conduit à un système de $2N+1$ ($N\to\infty$) d'équations algébriques. On les résolvant, peuvent être trouvées les coefficients $f_0$, $a_1$, ..., $a_N$, $b_1$, ..., $b_N$ de la série de Fourier du développement de la fonction $y(x)$. Pour

$y(x) = y_0 + y_e + y_o = y_0 + \sum_{n=1}^{\infty} A_n \cos\frac{n\pi x}{L} + \sum_{n=1}^{\infty} B_n \sin\frac{n\pi x}{L}$ et

$p(x) = p_0 + p_e + p_o = p_0 + \sum_{n=1}^{\infty} C_n \cos\frac{n\pi x}{L} + \sum_{n=1}^{\infty} D_n \sin\frac{n\pi x}{L}$

(où $y_{e0}=p_{e0}=0$), les coefficients de développement en série du produit $yp$ sont:

$$(yp)^{SFS} = y_0 p_0 + \left(y_0 \overline{p}_e + p_0 \overline{y}_e + \overline{y}_e \overline{p}_e + \overline{y}_o \overline{p}_o\right) + \left(y_0 \overline{p}_o + p_0 \overline{y}_o + \overline{y}_e \overline{p}_o + \overline{y}_o \overline{p}_e\right) =$$

où: (1.8)

$$= \overline{P}(x) = P_0 + \overline{P}_e + \overline{Q}_o = P_0 + \sum_{n=1}^{\infty} P_n \cos\frac{n\pi x}{L} + \sum_{n=1}^{\infty} Q_n \sin\frac{n\pi x}{L}$$

$$P_0 = \frac{1}{2L}\int_{-L}^{L}(yp)^{SFS}\,dx = y_0 p_0 + \frac{1}{2L}\int_{-L}^{L}\overline{y}_e \overline{p}_e\,dx + \frac{1}{2L}\int_{-L}^{L}\overline{y}_o \overline{p}_o\,dx = y_0 p_0 + \frac{1}{2}\left(\sum_{l=1}^{\infty} A_l C_l + \sum_{l=1}^{\infty} B_l D_l\right) \quad (1.9)$$

$$P_n = \frac{1}{L}\int_{-L}^{L}\left(y_0 \overline{p}_e + p_0 \overline{y}_e + \overline{y}_e \overline{p}_e + \overline{y}_o \overline{p}_o\right)\cos\frac{n\pi x}{L}\,dx = \frac{y_0}{L}\int_{-L}^{L}p_e \cos\frac{n\pi x}{L}\,dx + \frac{p_0}{L}\int_{-L}^{L}y_e \cos\frac{n\pi x}{L}\,dx +$$

$$+ \frac{1}{L}\int_{-L}^{L}\left(\sum_{l=1}^{\infty} A_l \cos\frac{l\pi x}{L}\right)\left(\sum_{m=1}^{\infty} C_m \cos\frac{m\pi x}{L}\right)\cos\frac{n\pi x}{L}\,dx + \frac{1}{L}\int_{-L}^{L}\left(\sum_{l=1}^{\infty} B_l \sin\frac{l\pi x}{L}\right)\left(\sum_{m=1}^{\infty} D_m \sin\frac{m\pi x}{L}\right)\cos\frac{n\pi x}{L}\,dx$$

$$\rightarrow P_n = y_0 C_n + p_0 A_n + \frac{1}{2}\sum_{l=1}^{\infty} A_l\left(C_{n+l} + C_{|n-l|}\right) + \frac{1}{2}\sum_{l=1}^{\infty} B_l\left[D_{n+l} + \mathrm{sgn}(l-n)D_{|n-l|}\right] \quad (1.10a)$$

$$Q_n = \frac{1}{L}\int_{-L}^{L}\left(y_0 \overline{p}_o + p_0 \overline{y}_o + \overline{y}_e \overline{p}_o + \overline{y}_o \overline{p}_e\right)\sin\frac{n\pi x}{L}\,dx = \frac{y_0}{L}\int_{-L}^{L}p_o \sin\frac{n\pi x}{L}\,dx + \frac{p_0}{L}\int_{-L}^{L}y_o \sin\frac{n\pi x}{L}\,dx +$$

$$+ \frac{1}{L}\int_{-L}^{L}\left(\sum_{l=1}^{\infty} A_l \cos\frac{l\pi x}{L}\right)\left(\sum_{m=1}^{\infty} D_m \sin\frac{m\pi x}{L}\right)\sin\frac{n\pi x}{L}\,dx + \frac{1}{L}\int_{-L}^{L}\left(\sum_{l=1}^{\infty} B_l \sin\frac{l\pi x}{L}\right)\left(\sum_{m=1}^{\infty} C_m \cos\frac{m\pi x}{L}\right)\sin\frac{n\pi x}{L}\,dx$$

$$\rightarrow Q_n = y_0 D_n + p_0 B_n + \frac{1}{2}\sum_{l=1}^{\infty} A_l\left[D_{n+l} - \mathrm{sgn}(l-n)D_{|n-l|}\right] - \frac{1}{2}\sum_{l=1}^{\infty} B_l\left(C_{n+l} - C_{|n-l|}\right) \quad (1.10b)$$

Dans ces relations, nous avons considéré: $A_0 = B_0 = C_0 = D_0 = 0$ (pour $l=n$). Les relations sont les mêmes que dans [9], dans lequel $A_0 = 2y_0$ et $C_0 = 2p_0$, mais ici, nous avons choisi une autre méthode pour atteindre le résultat, en utilisant des calculs intermédiaires tels que:

$$\int_{-L}^{L}\sum_{l=1}^{\infty}\left(A_l \cos\frac{l\pi x}{L}\right)\sum_{m=1}^{\infty}\left(C_m \cos\frac{m\pi x}{L}\right)\cos\frac{n\pi x}{L}\,dx = \int_{-L}^{L}\sum_{l=1}^{\infty}\left(A_l \cos\frac{l\pi x}{L}\right)\sum_{m=1}^{\infty}\frac{C_m}{2}\left[\cos\frac{(m-n)\pi x}{L} + \cos\frac{(m+n)\pi x}{L}\right]dx =$$

$$= \int_{-L}^{L}\sum_{l=1}^{\infty}\left(A_l \cos\frac{l\pi x}{L}\right)\left[\sum_{m=1}^{n-1}\frac{C_m}{2}\cos\frac{(m-n)\pi x}{L} + \frac{C_m}{2} + \sum_{m=n+1}^{\infty}\frac{C_m}{2}\cos\frac{(m-n)\pi x}{L} + \sum_{m=1}^{\infty}\frac{C_m}{2}\cos\frac{(m+n)\pi x}{L}\right]dx =$$

$$= \frac{L}{2}\left(\sum_{l=1}^{n-1} A_l C_{n-l} + \sum_{l=1}^{\infty} A_l C_{l+n} + \sum_{l=1}^{\infty} A_l C_{l-n}\right) = \frac{L}{2}\sum_{l=1}^{\infty} A_l\left(C_{l+n} + C_{|l-n|}\right)$$

De toute évidence, les formules suivantes sont également valables:

$$P_n = p_0 A_n + y_0 C_n + \frac{1}{2}\sum_{l=1}^{\infty} C_l\left(A_{n+l} + A_{|n-l|}\right) + \frac{1}{2}\sum_{l=1}^{\infty} D_l\left(B_{n+l} + \mathrm{sgn}(l-n)B_{|n-l|}\right) \quad (1.11)$$

$$Q_n = p_0 B_n + y_0 D_n + \frac{1}{2}\sum_{l=1}^{\infty} C_l\left(B_{n+l} - \mathrm{sgn}(l-n)B_{|n-l|}\right) - \frac{1}{2}\sum_{l=1}^{\infty} D_l\left(A_{n+l} - A_{|n-l|}\right)$$

## 1.4. Le calcul de certaines séries numériques infinies

Les coefficients $P_n$ et $Q_n$ des relations 1.9 et 1.10 contiennent, chacun, les expressions des séries numériques infinies, dont la somme peut être déterminée (si les coefficients $A_n$, $B_n$, $C_n$, et $Dn$ sont connus) par un volume, parfois important, de calculs. Entre chaque de ces séries infinies et le coefficient correspondant du développement en SFS du produit $yp$, il y a une relation directe, de sorte que la somme de série peut être déterminée en calculant le coefficient respectif. Par conséquent, si les calculs pour déterminer les coefficients de Fourier du produit $yp$ sont moins volumineux, cette méthode de calcul devient une solution préférable pour le calcul de cette quantité. Pour illustrer cela, nous comparerons ces volumes de calcul pour les produits entre la fonction:

$$\overline{y}(x) = y_0 + \sum_{n=1}^{\infty} A_n \cos\frac{n\pi x}{L} + \sum_{n=1}^{\infty} B_n \sin\frac{n\pi x}{L} = \overline{x} = \sum_{n=1}^{\infty} 2L\frac{(-1)^{n+1}}{n\pi}\sin\frac{n\pi x}{L} \quad \text{et les fonctions:}$$

$$\overline{p}_1(x) = p_{1.0} + \sum_{n=1}^{\infty} C_{1.n}\cos\frac{n\pi x}{L} + \sum_{n=1}^{\infty} D_{1.n}\sin\frac{n\pi x}{L} = \overline{x} = \sum_{n=1}^{\infty} 2L\frac{(-1)^{n+1}}{n\pi}\sin\frac{n\pi x}{L} \quad \text{et}$$

$$\overline{p}_2(x) = p_{2,0} + \sum_{n=1}^{\infty} C_{2,n} \cos\frac{n\pi x}{L} + \sum_{n=1}^{\infty} D_{2,n} \sin\frac{n\pi x}{L} = \left(x^2\right)^{SFS} = \frac{L^2}{3} + \sum_{n=1}^{\infty} \frac{4L^2}{\pi^2}\frac{(-1)^n}{n^2}\cos\frac{n\pi x}{L}$$

avec les coefficients: $y_0{=}0$, $A_n{=}0$, $B_n = 2L\frac{(-1)^{n+1}}{n\pi}$; $p_{1,0}{=}0$, $C_{1,n}{=}0$, $D_{1,n} = 2\frac{L}{\pi}\frac{(-1)^{n+1}}{n}$; $p_{2,0} = \frac{L^2}{3}$,

$C_{2,n} = \frac{4L^2(-1)^n}{n^2\pi^2}$, $D_{2,n}{=}0$. On découle les coefficients des produits:

$$\overline{P}_1 = P_{1,0} + \sum_{n=1}^{\infty} P_{1,n}\cos\frac{n\pi x}{L} + \sum_{n=1}^{\infty} Q_{1,n}\sin\frac{n\pi x}{L} = \left(yp_1\right)^{SFS} = \left(x^2\right)^{SFS} = \frac{L^2}{3} + \frac{4L^2}{\pi^2}\frac{(-1)^n}{n^2}\cos\frac{n\pi x}{L},$$

$$P_2 = P_{2,0} + \sum_{n=1}^{\infty} P_{2,n}\cos\frac{n\pi x}{L} + \sum_{n=1}^{\infty} Q_{2,n}\sin\frac{n\pi x}{L} = \left(yp_2\right)^{SFS} = \left(x^3\right)^{SFS} = (-1)^{n+1}\left[\frac{2L^3}{n\pi} - \frac{12L^3}{\pi^3 n^3}\right]\sin\frac{n\pi x}{L}$$

Nous allons garder les significations précédemment attribuées: $A_0{=}B_0{=}C_0{=}D_0{=}0$. Conform avec 1.9−1.10:

$$P_{1,0} = \frac{1}{2}\sum_{n=1}^{\infty} B_n \frac{2L(-1)^{n+1}}{n\pi} \quad \text{Pour } y(x){=}x{:} \ \rightarrow P_{1,0} = 2\sum_{n=1}^{\infty}\frac{L^2}{\pi^2}\frac{(-1)^{n+1}}{n}\frac{(-1)^{n+1}}{n} = \frac{2L^2}{\pi^2}\sum_{n=1}^{\infty}\frac{1}{n^2} = p_{2,0} \tag{1.12}$$

$$P_{1,n} = \frac{B_n}{2}\frac{2L}{\pi}\frac{(-1)^{2n+1}}{2n} + \frac{2L}{2\pi}\sum_{\substack{l=1\\l\neq n}}^{\infty} B_l\left[\frac{(-1)^{n+l+1}}{n+l} + \mathrm{sgn}(l-n)\frac{(-1)^{|n-l|+1}}{|n-l|}\right] = \frac{L}{\pi}\left[-\frac{B_n}{2n} + \sum_{\substack{l=1\\l\neq n}}^{\infty}(-1)^{n+l}\frac{2lB_l}{n^2-l^2}\right] \tag{1.13}$$

$$\rightarrow P_{1,n} = \frac{L}{\pi}\left[\frac{L(-1)^n}{n^2\pi} + \frac{4L}{\pi}\sum_{\substack{l=1\\l\neq n}}^{\infty}\frac{(-1)^n}{l^2-n^2}\right] = (-1)^n\frac{L^2}{\pi^2}\left(\frac{1}{n^2} + 4\frac{3}{4n^2}\right) = 4\frac{(-1)^n L^2}{n^2\pi^2} = C_{2,n}$$

$$Q_{1,n} = \frac{L}{\pi}\left\{y_0\frac{2(-1)^{n+1}}{n} - \frac{A_n}{2n} + \sum_{\substack{l=1\\l\neq n}}^{\infty} A_l\left[\frac{(-1)^{n+l+1}}{n+l} - \mathrm{sgn}(l-n)\frac{(-1)^{|n-l|+1}}{|n-l|}\right]\right\} \tag{1.14}$$

$$\rightarrow Q_{1,n} = \frac{L}{\pi}\left\{y_0\frac{2(-1)^{n+1}}{n} - \frac{A_n}{2n} + 2n\sum_{\substack{l=1\\l\neq n}}^{\infty}(-1)^{n+l+1}\frac{A_l}{n^2-l^2}\right\} = 0$$

$$P_{2,0} = y_0\frac{L^2}{3} + \frac{2L^2}{\pi^2}\sum_{n=1}^{\infty}\frac{(-1)^n}{n^2}A_n \quad \text{Pour } y(x){=}x{:} \ \rightarrow P_{2,0}{=}0 \tag{1.15}$$

$$P_{2,n} = y_0\frac{4L^2}{\pi^2}\frac{(-1)^n}{n^2} + \frac{A_n L^2}{3} + \frac{2L^2}{\pi^2}\left\{\frac{A_n}{4n^2} + \sum_{\substack{l=1\\l\neq n}}^{\infty} A_l\left[\frac{(-1)^{n+l}}{(n+l)^2} + \frac{(-1)^{|n-l|}}{(n-l)^2}\right]\right\} \rightarrow P_{n,2}{=}0 \tag{1.16}$$

$$Q_{2,n} = \frac{L^2 B_n}{3} - \frac{2L^2}{\pi^2}\left\{\frac{(-1)^{2n}}{4n^2}B_n + \sum_{\substack{l=1\\l\neq n}}^{\infty} B_l\left[\frac{(-1)^{n+l}}{(n+l)^2} - \frac{(-1)^{|n-l|}}{(n-l)^2}\right]\right\} \tag{1.17}$$

$$\rightarrow Q_{2,n} = \frac{L^2}{3}\frac{2L(-1)^{n+1}}{n\pi} - \frac{2L^2}{\pi^2}\left[\frac{(-1)^{2n}}{4n^2}\frac{2L(-1)^{n+1}}{n\pi} + \sum_{\substack{l=1\\l\neq n}}^{\infty}(-1)^{2l+n+1}\frac{2L}{l\pi}\frac{4nl}{(n^2-l^2)^2}\right] =$$

$$= \frac{2(-1)^{n+1}L^3}{3n\pi} - (-1)^{n+1}\frac{4L^3}{\pi^3}\left[\frac{1}{4n^3} + 4n\left(\frac{11}{16n^4} - \frac{\pi^2}{12n^2}\right)\right] = (-1)^{n+1}\left(\frac{2L^3}{n\pi} - \frac{12L^3}{n^3\pi^3}\right)$$

Pour obtenir ces résultats, nous avons collecté des données de la littérature spécialisée [10]:

$$\sum_{n=1}^{\infty}\frac{1}{n^2} = \frac{\pi^2}{6}, \ \sum_{\substack{l=1\\l\neq n}}^{\infty}\frac{1}{l^2-n^2} = \frac{3}{4n^2}, \ \text{et pour } \sum_{\substack{l=1\\l\neq n}}^{\infty}\frac{1}{(l^2-n^2)^2} = \frac{1}{4n^2}\left(\frac{11}{4n^2} - \frac{\pi^2}{3}\right) \ \text{nous sommes partis de:}$$

$$\sum_{l=0}^{\infty} \frac{1}{(l^2 - a^2)^2} = \frac{1}{2a^4} + \frac{\pi}{4a^3}ctg(a\pi) + \frac{\pi^2}{4a^2}\cos ec(a\pi), \quad \text{où} \quad a \neq 0, \ 1, \ 2,..., \ \infty. \text{ Pour } \forall n \in \mathbf{N}, \text{ si}$$

$a = n + \varepsilon \to n$, c'est-à dire $\varepsilon \in \mathbf{R} \to 0$ et $\sin(a\pi) = \sin(n\pi + \varepsilon\pi) = \sin(\varepsilon\pi) \to 0$. Alors,

$$\sum_{l=0}^{\infty} \frac{1}{(l^2 - a^2)^2} = \sum_{\substack{l=1 \\ l \neq n}}^{\infty} \frac{1}{(l^2 - a^2)^2} + \frac{1}{a^4} + \frac{1}{4an}\left[\frac{1}{(n-a)^2} - \frac{1}{(n+a)^2}\right]$$

$$\lim_{\varepsilon \to 0} \frac{1}{4n(n \pm \varepsilon)}\frac{1}{\varepsilon^2} + \frac{\pi}{4(n \pm \varepsilon)^3}\cot(\pm \varepsilon\pi) + \frac{\pi^2}{4(n \pm \varepsilon)^2}\csc(\pm \varepsilon\pi) = \frac{\pi^2}{12n^2} - \frac{1}{4n^4}$$

## 2. La solution d'une ODE, déterminée en calculant les coefficients de son développement en série de Fourier sinusoïdale

Pour résoudre les équations différentielles linéaires avec des coefficients $a_i(x)$ variables, lorsque ces coefficients sont des fonctions analytiques au point $x=x_0$, on peut utiliser la méthode du développement en série de Taylor: $y(x) = \sum_{n=1}^{\infty} a_n(x - x_0)^n$, si $x_0$ est un point ordinaire, ou la méthode de Frobenius: $y(x) = \sum_{n=1}^{\infty} a_n(x - x_0)^{n+r}$, si $x_0$ est un point singulier régulier [4, 5]. La méthode, qui permet de trouver une solution valide sur son disque de convergence de rayon R, est basée sur la propriété de dérivabilité terme par terme de ces séries et permet de trouver des relations de récurrence pour les coefficients de la série. C'est possible parce que la base sur laquelle la fonction–solution générale est développée est la même que la base sur laquelle sont développées ses dérivés aussi. Pour trouver les solutions particulières il faut savoir les valeurs de la fonction et de ses dérivés à un point de l'intervalle $(-R, R)$. Par cette méthode, la résolution d'équations différentielles est transformée en un problème d'algèbre, qui implique le plus souvent de trouver des formules de récurrence.

Étant donné que les SFS et les SFN aussi, jouissent de la propriété de la dérivabilité terme par terme, une méthode similaire à la méthode décrite peut être appliquée en utilisant ces nouveaux types de développements en série, dans les équations impliquant des fonctions périodiques de carré intégrable, définies sur un intervalle $[-L, L]$, si nous connaissons les conditions aux limites correspondantes.

La méthode de résolution des équations différentielles en déterminant les coefficients du développement en séries sinusoïdales ou non sinusoïdales de la fonction inconnue est une méthode particulièrement solide, applicable à tous les types d'équations différentielles et intégro-différentielles, linéaires et non linéaires, d'équations aux dérivées partielles, des systèmes de telles équations, quels que soient leur ordre et quels que soient la complexité des coefficients. Les conditions requises pour l'application de la méthode sont peu contraignantes et faciles à remplir, notamment pour les situations rencontrées en physique et en génie. La méthode peut être appliquée à des larges classes de telles équations, et encore plus, la même équation peut être résolue en utilisant plusieurs types de développement. La méthode peut également être facilement étendue aux fonctions du domaine complexe C.

Comme pour la méthode de développement en série entière, les opérations de dérivation et d'intégration effectuées sur ces développements en série sont transformées en opérations algébriques effectuées sur les coefficients du développement. Par conséquent, la résolution d'ODE se transforme en résolution d'équations algébriques. De toutes les bases utilisées comme support des développements, les plus avantageuses sont les sinusoïdes (*1, cos(nx), sin(nx)*), en raison de la facilité avec laquelle les coefficients du développement peuvent être calculés (les formules d'Euler), en raison de la vitesse de convergence élevée (ce qui les rend idéales pour les méthodes d'approximation numérique) et la relative facilité avec laquelle les équations algébriques résultantes peuvent être résolues (le degré de difficulté est dicté par l'ordre de l'équation et par le degré de non-linéarité de l'équation différentielle). Contrairement aux bases des puissances entières positives qui, par des dérivations répétées, "patinent" le long des éléments de la base (conduisant à la génération des chaînes récurrentes), les bases sinusoïdales "oscillent" entre les mêmes éléments, pairs et impairs, de la base, conduisant à des équations algébriques plus simples.

## 2.1. ODEs linéaires à coefficients *constants*

Comme nous l'avons trouvé dans la section précédente, sur l'intervalle $[-L, L]$, pour

$$\overline{y}(x) = y_0 + \overline{y}_e + \overline{y}_o = y_0 + \sum_{n=1}^{\infty} A_n \cos \omega_n x + \sum_{n=1}^{\infty} B_n \sin \omega_n x \tag{2.1}$$

$$\rightarrow \overline{y}'(x) = \sum_{n=1}^{\infty} \left[ B_n \omega_n + 2(-1)^n \frac{y(L) - y(-L)}{2L} \right] \cos \omega_n x - \sum_{n=1}^{\infty} (A_n \omega_n \sin \omega_n x) + \frac{y(L) - y(-L)}{2L} \tag{2.2}$$

$$\rightarrow \overline{y}'' = \sum_{n=1}^{\infty} \left[ \left( -A_n \omega_n^2 + 2(-1)^n \frac{y_e'(L)}{L} \right) \cos \omega_n x - \left( B_n \omega_n^2 + 2(-1)^n \omega_n \frac{y_o(L)}{L} \right) \sin \omega_n x \right] + \frac{y_e'(L)}{L} \tag{2.3}$$

Aussi, $y_o(L) = L \lim_{n \to \infty} 2 \left[ (-1)^{n+1} \omega_n b_n \right] = 1/2 \cdot \left[ y(L) - y(-L) \right]$ \hfill (2.1a)

Plus loin, par des dérivations successives (à chaque dérivation, les coefficients de développement sont corrigés en tenant compte de l'existence des points de discontinuité) on retrouve les expressions des développements en série pour les dérivées d'ordre supérieur, lorsque sont connues leurs valeurs aux points à la limite de l'intervalle de définition. Par conséquent, la solution particulière *y(x)* de **toute** équation différentielle ordinaire linéaire avec des coefficients constants, homogène ou inhomogène, quel que soit l'ordre de l'équation, si elle est définie sur n'importe quel intervalle $[-L, L]$, ou équivalent, peut être déterminée si *y(x)* est de carré intégrable, en calculant les coefficients de son développement en série sinusoïdale, si les valeurs de la fonction et de ses dérivées sont connues aux points des extrémités de l'intervalle (conditions aux limites). Les solutions des équations peuvent être trouvées pour d'autres types de conditions aux limites aussi.

Pour illustrer cela, nous appliquerons la méthode dans le cas d'équations linéaires homogènes à coefficients constants, sur l'intervalle $[-\pi, \pi]$, pour quelques équations simples:

**Exemple 2.1: *y' = a*,** \hfill (2.4)

avec la condition à la limite *y(0)=C*. Pour l'intervalle choisi, on peut écrire:

$$\sum_{n=1}^{\infty} \left[ B_n n + 2(-1)^n \frac{y(\pi) - y(-\pi)}{2\pi} \right] \cos nx + \sum_{n=1}^{\infty} (-A_n n \sin nx) + \frac{y(\pi) - y(-\pi)}{2\pi} = a \tag{2.5}$$

$$\rightarrow A_n=0, \quad B_n = \frac{y(\pi) - y(-\pi)}{2\pi} \frac{2(-1)^{n+1}}{n}, \quad \frac{y(\pi) - y(-\pi)}{2\pi} = \frac{y_o(\pi)}{\pi} = a$$

Pour la valeur moyenne $y_0$, il n'y a aucun conditionnement, donc ça peut prendre n'importe quelle valeur *K*. Donc:

$$y(x) = a \frac{2(-1)^{n+1}}{n} \sin nx + K = ax + K \text{. Pour } y=0 \rightarrow y(0)=K=C.$$

Pour *a=0 → y(x)=K*

**Exemple 2.2: *y' = ay*,** \hfill (2.6)

avec la condition à la limite *y(π)=C*. A partir des relations (2.2) et (2.6), on obtient:

$$ay_0 = \frac{y_o(\pi)}{\pi}, \quad aB_n = -nA_n \rightarrow B_n = -nA_n/a \text{, pour } n=1, 2, 3, ...\infty \text{ et}$$

$$aA_n = -\frac{n^2}{a}A_n + 2(-1)^n \frac{y_o(\pi)}{\pi} \rightarrow A_n = 2\frac{(-1)^n a}{n^2+a^2}\frac{y_o(\pi)}{\pi} \rightarrow B_n = -2n\frac{(-1)^n}{n^2+a^2}\frac{y_o(\pi)}{\pi} \text{ et donc:}$$

$$\bar{y}(x) = 2\frac{y_o(\pi)}{\pi}\left[\frac{1}{2a} + \sum_{n=1}^{\infty}(-1)^n \frac{a\cos nx - n\sin nx}{n^2+a^2}\right] \tag{2.7}$$

On reconnaît ici, pour $y_o(\pi)=sinha\pi$, le développement en série de Fourier de la fonction $y(x)=e^{ax}$, pour lequel $y(\pi)=e^{a\pi}$, et pour $y_o(\pi)=Ksinha\pi$, le développement en série de Fourier de la fonction $y(x)=K{\cdot}e^{ax}$ (la solution générale), pour lequel $y(\pi)=Ke^{a\pi}$, et $y(-\pi)=Ke^{-a\pi}$.
En conclusion, pour $y(\pi)=C \rightarrow K=Ce^{-a\pi}$, la solution particulière étant $y(x)=Ce^{a(x-\pi)}$.

Pendant la résolution de l'équation, est apparue la contrainte $aB_n=-nA_n$. Cela signifie que $A_n=0$ impose $B_n=0$, donc la solution ne peut pas avoir une seule composante (elle ne peut pas être seulement paire ou impaire). Ce fait est imposé même par la relation d'égalité (2.4). Donc $f_o(\pi)/\pi{\neq}0$, par conséquent, toujours $y_o{\neq}0$.

**Exemple 2.3:** $y'' = a^2 y$, avec $y(\pi)=K_1$ et $y'(\pi)=K_2$. Parce que:

$$\bar{y}'(x) = \bar{y}'_e + \bar{y}'_o = \sum_{n=1}^{\infty}\left(-A_n n\sin nx\right) + \sum_{n=1}^{\infty}\left[B_n n + 2(-1)^n \frac{y_o(\pi)}{\pi}\right]\cos nx + \frac{y_o(\pi)}{\pi} \tag{2.8}$$

$$\bar{y}'' = \sum_{n=1}^{\infty}\left[\left(-A_n n^2 + 2(-1)^n \frac{y'_e(\pi)}{\pi}\right)\cos(nx) - \left(B_n n^2 + 2(-1)^n n\frac{y_o(\pi)}{\pi}\right)\sin(nx)\right] + \frac{y'_e(\pi)}{\pi} \tag{2.9}$$

et par l'égalité avec $a^2 y$ on découvre que:

$$y_0 = \frac{y'_e(\pi)}{a^2\pi}, \quad A_n = 2\frac{(-1)^n}{n^2+a^2}\frac{y'_e(\pi)}{\pi} \text{ et } B_n = -2n\frac{(-1)^n}{n^2+a^2}\frac{y_o(\pi)}{\pi}, \text{ donc}$$

$$\bar{y}(x) = \frac{y'_e(\pi)}{a^2\pi} + 2\sum_{n=1}^{\infty}\frac{(-1)^n}{n^2+a^2}\left[\frac{y'_e(\pi)}{\pi}\cos nx - n\frac{y_o(\pi)}{\pi}\sin nx\right] \tag{2.10}$$

Pour $y_o(\pi)=sinh(a\pi)$ et $y'_e(\pi)=a{\cdot}sinh(a\pi)$, on reconnaît le développement en série de la fonction: $\bar{y}(x) = \cosh(ax) + \sinh(ax) = e^{ax}$. Mais, parce qu'il n'y a pas des contraintes entre les coefficients des deux composants (pairs et impairs), chacun d'eux peut être une solution: $y_1=cosh(ax)$ et $y_2=sinh(ax)$ (pour $f_o(\pi)=0$, respectivement $f'_e(\pi)=0$). Par conséquent, la solution générale de l'équation est $y(x)=C_1cosh(ax)+C_2sinh(ax)$. La solution particulière résulte de: $y(\pi)=C_1cosh(a\pi)+C_2sinh(a\pi)=K_1$ et $y'(\pi)=aC_1sinh(a\pi)+aC_2cosh(a\pi)=K_2$

**Exemple 2.4:** $y'' = -a^2 y$

De la relation (2.10), en remplaçant $a^2$ avec $-a^2$, on a:

$$\bar{y}(x) = -\frac{y'_e(\pi)}{a^2\pi} + 2\frac{(-1)^n}{n^2-a^2}\left[\frac{y'_e(\pi)}{\pi}\cos nx - n\frac{y_o(\pi)}{\pi}\sin nx\right],$$

dans lequel, pour $y_o(\pi)=sin(a\pi)$ et $y'_e(\pi)=-a{\cdot}sin(a\pi)$ nous reconnaissons, pour $a{\notin}Z$, le développement en série de la fonction:

$$\bar{y}(x) = \cos ax + \sin ax = \frac{\sin a\pi}{a\pi} + 2\frac{(-1)^{n+1}}{n^2-a^2}\left[\frac{a\sin a\pi}{\pi}\cos nx + n\frac{\sin a\pi}{\pi}\sin nx\right]$$

Par conséquent, la solution générale de cette équation est $y(x)=C_1cos(ax)+C_2sin(ax)$.

**Exemple 2.5:** $y''+ay'+by=0$
En remplaçant dans l'équation la fonction inconnue et ses dérivées par leurs développements en SFS:

$$\sum_{n=1}^{\infty}\left[-A_n n^2 + 2(-1)^n \frac{y'_e(\pi)}{\pi} + aB_n n + 2(-1)^n a\frac{y_o(\pi)}{\pi} + bA_n\right]\cos(nx) +$$

$$+ \sum_{n=1}^{\infty}\left[-B_n n^2 - 2(-1)^n n\frac{y_o(\pi)}{\pi} - aA_n n + bB_n\right]\sin(nx) + \frac{y'_e(\pi)}{\pi} + a\frac{y_o(\pi)}{\pi} + by_0 = 0$$

Tous les coefficients de cette development sont nuls:

$$y_0 = -\frac{y'_e(\pi) + a y_o(\pi)}{b\pi} \quad , \quad B_n = A_n \frac{n^2 - b}{an} - \frac{2(-1)^n}{an}\left[\frac{y'_e(\pi)}{\pi} + \frac{a y_o(\pi)}{\pi}\right], \quad A_n = -B_n \frac{n^2 - b}{an} - \frac{2(-1)^n}{a}\frac{y_o(\pi)}{\pi}$$

$$\to A_n = -A_n\left(\frac{n^2 - b}{an}\right)^2 + 2(-1)^n \frac{y'_e(\pi) + a y_o(\pi)}{an\pi}\frac{n^2 - b}{an} - 2(-1)^n \frac{y_o(\pi)}{a\pi}$$

$$A_n(n^4 + a^2 n^2 - 2bn^2 + b^2)\pi = 2(-1)^n\left[y'_e(\pi)(n^2 - b) + a(2n^2 - b)y_o(\pi)\right]$$

$$\to A_n = \frac{2(-1)^n\left[y'_e(\pi)(n^2 - b) - aby_o(\pi)\right]}{\pi\left[n^2 + \frac{1}{2}(a^2 - 2b - a\Delta)\right]\left[n^2 + \frac{1}{2}(a^2 - 2b + a\Delta)\right]} = \frac{2(-1)^n\left[y'_e(\pi)(n^2 - b) - aby_o(\pi)\right]}{\pi\left[n^2 + \left(\frac{-a+\Delta}{2}\right)^2\right]\left[n^2 + \left(\frac{-a-\Delta}{2}\right)^2\right]},$$

où $\Delta = \sqrt{a^2 - 4b}$ . Nous noterons $\lambda_1 = \frac{-a + \sqrt{a^2 - 4b}}{2} = \frac{-a + \Delta}{2}$ et $\lambda_2 = \frac{-a - \Delta}{2}$

$$B_n = -\left(B_n \frac{n^2 - b}{an}\right)^2 - \frac{2(-1)^n}{a}\frac{y_o(\pi)}{\pi}\frac{n^2 - b}{an} - \frac{2(-1)^n}{an}\left[\frac{y'_e(\pi)}{\pi} + \frac{a y_o(\pi)}{\pi}\right]$$

$$B_n(n^4 + a^2 n^2 - 2bn^2 + b^2)\pi = 2(-1)^n n\left[y_o(\pi)(n^2 + a^2 - b) + a n y'_e(\pi)\right]$$

$$\to B_n = \frac{2(-1)^n n\left[y_o(\pi)(n^2 + a^2 - b) + a y'_e(\pi)\right]}{\pi\left[n^2 + \left(\frac{-a+\Delta}{2}\right)^2\right]\left[n^2 + \left(\frac{-a-\Delta}{2}\right)^2\right]} = \frac{2(-1)^n n\left[y_o(\pi)(n^2 + a^2 - b) + a y'_e(\pi)\right]}{\pi(n^2 + \lambda_1^2)(n^2 + \lambda_2^2)}$$

Pour $y_{o1}(\pi) = C_1 \sinh\frac{-a+\Delta}{2}\pi = C_1 \sinh(\lambda_1\pi)$ et $y'_{e1}(\pi) = C_1\lambda_1 \sinh(\lambda_1\pi) = \lambda_1 y_{o1}(\pi)$:

$$y_{0.1} = -\frac{C_1\lambda_1 y_{o1}(\pi) + C_1 a y_{o1}(\pi)}{b\pi} = -C_1\frac{y_{o1}(\pi)}{\pi}\frac{-a+\Delta+2a}{2b} = -C_1\frac{y_{o1}(\pi)}{\pi}\frac{4b}{2b(a-\Delta)} = C_1\frac{\sinh(\lambda_1\pi)}{\pi\lambda_1}$$

$$A_{n.1} = \frac{2(-1)^n C_1 \sinh(\lambda_1\pi)\left[\lambda_1(n^2 - b) - a\lambda_1\lambda_2\right]}{\pi(n^2 + \lambda_1^2)(n^2 + \lambda_2^2)} = 2(-1)^n C_1\frac{\sinh(\lambda_1\pi)}{\pi}\frac{\lambda_1}{n^2 + \lambda_1^2}$$

$$B_{n.1} = \frac{2(-1)^n C_1 n \sinh(\lambda_1\pi)(n^2 + a^2 - b + a\lambda_1)}{\pi(n^2 + \lambda_1^2)(n^2 + \lambda_2^2)} = 2(-1)^n C_1\frac{\sinh(\lambda_1\pi)}{\pi}\frac{n}{n^2 + \lambda_1^2},$$

Et pour $y_{o2}(\pi) = C_2 \sinh\frac{-a-\Delta}{2}\pi = C_2 \sinh(\lambda_2\pi)$ et $y'_{e2}(\pi) = C_2\lambda_2 \sinh(\lambda_2\pi) = \lambda_2 y_{o2}(\pi)$:

$$y_{0.2} = C_2\frac{\sinh(\lambda_2\pi)}{\pi\lambda_2}, \quad A_{n.2} = 2(-1)^n C_2\frac{\sinh(\lambda_2\pi)}{\pi}\frac{\lambda_2}{n^2 + \lambda_2^2}, \quad B_{n.2} = 2(-1)^n C_2\frac{\sinh(\lambda_2\pi)}{\pi}\frac{n}{n^2 + \lambda_2^2}$$

Les deux solutions sont indépendants. Par conséquent, la solution générale de l'équation est:

- pour $a^2 > 4b$ : $y(x) = C_1 \exp\left(\frac{-a + \sqrt{a^2 - 4b}}{2}x\right) + C_2 \exp\left(\frac{-a - \sqrt{a^2 - 4b}}{2}x\right)$

- pour $a^2 < 4b$ : $y(x) = \exp\left(\frac{-a}{2}x\right)\left[C_1 \sin\left(\frac{\sqrt{a^2 - 4b}}{2}x\right) + C_2 \cos\left(\frac{\sqrt{a^2 - 4b}}{2}x\right)\right]$

- pour $a^2 < 4b$ : $\lambda_1 = \lambda_2 = \frac{-a}{2}$. $y_1(x) = C_1 \exp\left(\frac{-a}{2}x\right)$ est une solution de l'équation. On peut chercher une deuxième solution, sous la forme $y_2(x) = u(x)y_1(x)$. Par la méthode de "la réduction de l'ordre" [4] on trouve pour $u(x) = u_0 + \sum_{n=1}^{\infty} C_n \cos\frac{n\pi x}{L} + \sum_{n=1}^{\infty} D_n \sin\frac{n\pi x}{L}$ la condition:

$$\overline{u}'' = \sum_{n=1}^{\infty}\left[\left(-C_n n^2 + 2(-1)^n \frac{u'_e(\pi)}{\pi}\right)\cos(nx) - \left(D_n n^2 + 2(-1)^n n \frac{u_o(\pi)}{\pi}\right)\sin(nx)\right] + \frac{u'_e(\pi)}{\pi} = 0$$

$$\to \frac{u'_e(\pi)}{\pi} = 0, \; -C_n n^2 + 2(-1)^n \frac{u'_e(\pi)}{\pi} = 0 \to C_n{=}0, \; D_n = 2(-1)^n \frac{u_o(\pi)}{n\pi} \to u(x) = C_1 x + C_2$$

On note que pour tous les ODE à coefficients constants, quel que soit le degré de l'équation, les coefficients de développement de la fonction inconnue sont fournis par une équation pour déterminer la valeur moyenne et une paire d'équations pour chaque harmonique du développement.

**Exemple 2.6: $y'=x$**

Même dans le cas d'ODEs non homogènes à coefficients constants, il n'y a pas d'équations significativement plus difficiles à résoudre. Par exemple, si dans cette équation, sur l'intervalle $[-L, L]$, nous faisons les remplacements:

$$\overline{y}(x) = y_0 + \sum_{n=1}^{\infty} A_n \cos nx + \sum_{n=1}^{\infty} B_n \sin nx \; \text{ et } \; \overline{x} = 2\sum_{n=1}^{\infty} \frac{(-1)^{n+1}}{n} \sin nx \text{ , on a:}$$

$$\sum_{n=1}^{\infty}(-A_n n \sin nx) + \sum_{n=1}^{\infty} n\left[B_n + 2(-1)^n \frac{y_o(\pi)}{n\pi}\right]\cos nx + \frac{y_o(\pi)}{\pi} = 2\sum_{n=1}^{\infty} \frac{(-1)^{n+1}}{n} \sin nx$$

$$\to y_0{=}K \text{ (arbitraire)}, \, y_o(\pi){=}0, \, \to \, B_n + 2(-1)^n \frac{y_o(\pi)}{n\pi} = 0 \to B_n{=}0 \text{ et } A_n{=}-2\frac{(-1)^{n+1}}{n^2}$$

On reconnaît le développement: $\dfrac{x^2}{2} = \dfrac{1}{6} - 2\sum_{n=1}^{\infty} \dfrac{(-1)^{n+1}}{n^2}\cos nx \; \to \; y(x) = y_0 + \dfrac{x^2}{2} - \dfrac{1}{6} = \dfrac{x^2}{2} + C$

**Exemple 2.7: $y'=x^2$.**

$$\sum_{n=1}^{\infty}(-A_n n \sin nx) + \sum_{n=1}^{\infty}\left[B_n n + 2(-1)^n \frac{y_o(\pi)}{\pi}\right]\cos nx + \frac{y_o(\pi)}{\pi} = \frac{\pi^2}{3} - 4\sum_{n=1}^{\infty} \frac{(-1)^{n+1}}{n^2}\cos nx$$

$$\to y_0{=}K \text{ (arbitraire)}, \; \to \; y_o(\pi) = \frac{\pi^3}{3}, \; A_n = 0, \; B_n = \frac{2(-1)^{n+1}}{n}\frac{\pi^2}{3} - 4\sum_{n=1}^{\infty} \frac{(-1)^{n+1}}{n^3},$$

$$\to y(x) = y_0 + \left[\frac{2(-1)^{n+1}}{n}\frac{\pi^2}{3} - 4\sum_{n=1}^{\infty} \frac{(-1)^{n+1}}{n^3}\right]\cos nx = \frac{x^3}{3} + C$$

## 2.2. ODE linéaires à coefficients *variables*

Une équation différentielle linéaire d'ordre $m$ avec des coefficients variables, non homogènes, de la forme

$f_m(x)y^{(m)}+...+f_2(x)y''+f_1(x)y'+f_0(x)y=g(x)$

peut être résolu, sur l'intervalle $[-L, L]$, si $y^{(i)}(x)$, $f_i(x)$ et $g(x)$, $i=0, 1, 2, ...,m$, sont des fonctions de carré intégrable, en déterminant les coefficients $A_n$ et $B_n$, $n=1, 2, ..., \infty$, du développement en SFS de la fonction inconnu $y(x)$:

$$y(x) = y_0 + \sum_{n=1}^{\infty} A_n \cos \omega_n x + \sum_{n=1}^{\infty} B_n \sin \omega_n x$$

si les valeurs de la fonction et de ses dérivées sont connues aux points situés aux extrémités de l'intervalle. Si nous ne connaissons pas leur valeur qu'à l'une de ces extrémités, la relation (2.1a) nous permet également de trouver l'autre. Si $L$ est pris comme paramètre, des estimations peuvent être faites, si $L{\to}\infty$, pour toutes les valeurs sur l'axe réel pour lesquelles les conditions d'intégrabilité au carré sont remplies. À partir des relations (2.1)–(2.3) de la sous-section précédente, par dérivations successives, nous trouverons les expressions des développements en série pour les dérivées d'ordre supérieur (si leurs valeurs aux extrémités de l'intervalle de définition sont connues). En utilisant des relations similaires aux relations (1.9) et (1.10) antérieurement déterminées (relations pour les coefficients du produit des deux fonctions), ainsi que des relations déduites de la formule d'intégration par parties on peut obtenir des expressions (dépendantes de $A_n$ et $B_n$) pour les coefficients $P_n$ et $Q_n$ du développement en série pour les termes de la forme $f_i(x)y^{(i)}$, lesquelles sont introduits dans l'équation de base, simultanément avec l'expression du développement en série de la fonction $g(x)$. En simplifiant la relation résultante (regroupant les termes qui ont comme facteur commun l'un des éléments de la base du développement, y compris ici aussi la fonction unitaire *1*), on obtient l'expression du développement en SFS d'une fonction de valeur identiquement nulle. Par conséquent, tous les coefficients de ce développement (expressions algébriques dans lesquelles apparaît les coefficients $A_n$ et $B_n$) sont nuls et donnent naissance à des relations (équations algébriques) qui permettent le calcul des inconnues $A_n$, $B_n$ et $f_0$. La méthode est similaire à celle dans laquelle ces termes sont déterminés en remplaçant la fonction *y(x)* par son développement en série de Taylor.

Les coefficients d'ordre $N$, $P_N$ et $Q_N$, des développements des produits $f_i(x)y^{(i)}$ sont des séries numériques qui peuvent également contenir des quantités infinies de termes décroissants ou alternants décroissants. Lorsque ces sommes peuvent être calculées (elles peuvent être réduites à une somme finie de termes), on obtienes des équations exactes, et en les résolvant, on obtienes des valeurs exactes pour les coefficients du développement de la fonction inconnue *y(x)*. Pour cela, il est nécessaire que les coefficients du développement en SFS des fonctions $f_i(x)$ puissent être calculés avec exactitude (les intégrales d'Euler que leur correspond puissent être calculées avec un nombre fini de termes). Au contraire, on peut limiter au premier $N$, le nombre de termes de ces séries. Ainsi, *2N+1* équations approximatives seront obtenues, avec *2N+1* variables inconnues, par la résolution desquelles sont obtenues des valeurs approximatives pour $y_0$, $A_n$ et $B_n$. Pour la fonction *y(x)* on obtienne une valeur approximée par la somme de Fourier $S_N$. Plus élevé est le nombre $N$ de termes, meilleure est l'approximation.

Dans certains cas, les solutions des équations pour certains des coefficients peuvent être déduites en comparant les termes des équations avec les termes des développements en série de certaines fonctions connues, ce qui facilite la recherche de la solution globale.

Une situation fréquent est celle des équations à coefficients polynomiaux. Si les conditions aux limites sont données (les valeurs de la fonction inconnue *y* et de ses dérivées aux extrémités de l'intervalle considérée), elles sont suffisantes pour écrire les équations qui déterminent les coefficients $A_n$ și $B_n$. Si *y(x), $y^{(m)}(x)$, $f_i(x)$* et *g(x)* sont des carré intégrable sur l'intervalle $[-L, L]$, les relations (1.12)–(1.17) et ceux qui peuvent en être dérivés, fournissent des relations de calcul pour tous les coefficients du développement en SFS des termes de ces équations. Voilà, par exemple, les relations calculées pour les termes des équations du premier et du second ordre, ayant des coefficients polynomiaux du premier et du deuxième degré. Si

$$\bar{y}(x) = y_0 + \sum_{n=1}^{\infty} A_n \cos\frac{n\pi x}{L} + \sum_{n=1}^{\infty} B_n \sin\frac{n\pi x}{L} \text{ , alors:}$$

$$(yx)^{SFS} = \frac{L}{\pi}\sum_{n=1}^{\infty}\frac{(-1)^{n+1}}{n}B_n + \frac{L}{\pi}\left[-\frac{B_n}{2n} + \sum_{\substack{l=1\\l\neq n}}^{\infty}(-1)^{l+n}\frac{2lB_l}{n^2-l^2}\right]\cos\frac{n\pi x}{L} +$$

$$+\frac{L}{\pi}\left[y_0\frac{2(-1)^{n+1}}{n} - \frac{A_n}{2n} + 2n(-1)^{n+1}\sum_{\substack{l=1\\l\neq n}}^{\infty}\frac{(-1)^l A_l}{n^2-l^2}\right]\sin\frac{n\pi x}{L}$$

$$\left(yx^2\right)^{SFS} = L^2\left[\frac{y_0}{3} + \sum_{n=1}^{\infty}\frac{2(-1)^n A_n}{n^2\pi^2}\right] + \frac{L^2}{\pi^2}\left[\frac{4y_0(-1)^n}{n^2} + \frac{A_n(2\pi^2 n^2+3)}{6n^2} + 4\sum_{\substack{l=1\\l\neq n}}^{\infty}A_l(-1)^{n+l}\frac{\left(n^2+l^2\right)}{\left(n^2-l^2\right)^2}\right]\cos\frac{n\pi x}{L} +$$

$$+ \frac{L^2}{\pi^2}\left[\frac{B_n(2\pi^2 n^2+3)}{6n^2} + 8n\sum_{\substack{l=1\\l\neq n}}^{\infty}(-1)^{n+l}\frac{lB_l}{\left(n^2-l^2\right)^2}\right]\sin\frac{n\pi x}{L}$$

$$\left(y'x\right)^{SFS} = \sum_{n=1}^{\infty}(-1)^n A_n + \left[\frac{A_n}{2} - \sum_{\substack{l=1\\l\neq n}}^{\infty}(-1)^{l+n}\frac{2l^2 A_l}{n^2-l^2}\right]\cos\frac{n\pi x}{L} + \left[-\frac{B_n}{2} + 2n(-1)^{n+1}\sum_{\substack{l=1\\l\neq n}}^{\infty}\frac{lB_l(-1)^l}{n^2-l^2}\right]\sin\frac{n\pi x}{L}$$

$$\left(y'x^2\right)^{SFS} = \frac{y_o(L)L}{3} + \frac{2L}{\pi^2}\sum_{n=1}^{\infty}\left[\frac{(-1)^n \pi B_n}{n} + \frac{2y_o(L)}{n^2}\right] +$$

$$+ \sum_{n=1}^{\infty}\left\{\frac{(-1)^n y_o(L)(4\pi^2 n^2+30)}{6n^2 L} + \frac{B_n(2\pi^2 n^2+3)n\pi}{6n^2 L} + 4(-1)^n\sum_{\substack{l=1\\l\neq n}}^{\infty}\left[\frac{\pi}{L}B_l(-1)^l l\frac{\left(n^2+l^2\right)}{\left(n^2-l^2\right)^2} + 2\frac{y_o(L)}{L}\frac{\left(n^2+l^2\right)}{\left(n^2-l^2\right)^2}\right]\right\}\cos\frac{n\pi x}{L} +$$

$$+ \sum_{n=1}^{\infty}\frac{L}{\pi}\left[\frac{-A_n(2\pi^2 n^2-3)}{6n} - 8n\sum_{\substack{l=1\\l\neq n}}^{\infty}(-1)^{n+l}\frac{l^2 A_l}{\left(n^2-l^2\right)^2}\right]\sin\frac{n\pi x}{L}$$

$$\left(y''x\right)^{SFS} = \sum_{n=1}^{\infty}\left[(-1)^n\frac{n\pi}{L}B_n + 2\frac{y_o(L)}{L}\right] + \sum_{n=1}^{\infty}\left\{\frac{B_n n\pi}{2L} + (-1)^n\frac{y_o(L)}{L} + (-1)^n\sum_{\substack{l=1\\l\neq n}}^{\infty}\frac{(-1)^{l+1}2l^2}{n^2-l^2}\left[lB_l\frac{\pi}{L} - 2\frac{y_o(L)}{L}\right]\right\}\cos\frac{n\pi x}{L} +$$

$$+ \sum_{n=1}^{\infty}\left\{\frac{3(-1)^{n+1}y'_e(L)}{n\pi} + \frac{A_n}{2}\frac{n\pi}{L} + 2n(-1)^{n+1}\sum_{\substack{l=1\\l\neq n}}^{\infty}\left[\frac{(-1)^{l+1}A_l l^2}{n^2-l^2}\frac{\pi}{L} + 2\frac{y'_e(L)}{\pi}\frac{1}{n^2-l^2}\right]\right\}\sin\frac{n\pi x}{L}$$

$$\left(y''x^2\right)^{SFS} = y'_e(L)L - 2\sum_{l=1}^{\infty}(-1)^l A_l +$$

$$+ \sum_{n=1}^{\infty}\left\{2(-1)^n\frac{y'_e(L)}{L}\frac{L^2(2\pi^2 n^2+15)}{6\pi^2 n^2} - \frac{A_n(2\pi^2 n^2+3)}{6} + (-1)^n\sum_{\substack{l=1\\l\neq n}}^{\infty}\left[A_l(-1)^{l+1}\frac{4l^2(n^2+l^2)}{\left(n^2-l^2\right)^2} + \frac{8L^2}{\pi^2}\frac{y'_e(L)}{L}\frac{\left(n^2+l^2\right)}{\left(n^2-l^2\right)^2}\right]\right\}\cos\frac{n\pi x}{L} +$$

$$+ \sum_{n=1}^{\infty}\left\{\frac{(2\pi^2 n^2-3)}{6\pi n}\left[-\pi n B_n - 2(-1)^n y_o(L)\right] - 8n\sum_{\substack{l=1\\l\neq n}}^{\infty}\left[\frac{l^3 B_l(-1)^{n+l}}{\left(l^2-n^2\right)^2}\right] - \frac{16nL}{\pi}\frac{y_o(L)}{L}\sum_{\substack{l=1\\l\neq n}}^{\infty}\frac{(-1)^l l^2}{\left(l^2-n^2\right)^2}\right\}\sin\frac{n\pi x}{L}$$

Par exemple, une équation de type Euler, non homogène $y''x^2 + ay'x + by = g(x)$, après avoir utilisé ces relations et après le développement en SFS de la fonction

$g(x) = g_0 + \sum_{n=1}^{\infty}P_n\cos\frac{n\pi x}{L} + \sum_{n=1}^{\infty}Q_n\sin\frac{n\pi x}{L}$, conduit aux équations algébriques suivantes:

$$y'_e(L)L - 2\sum_{n=1}^{\infty}(-1)^n A_n + a\sum_{n=1}^{\infty}(-1)^n A_n + by_0 = g_0 \;\rightarrow\; y_0 = \frac{1}{b}\left[g_0 - y'_e(L)L + (2-a)\sum_{n=1}^{\infty}(-1)^n A_n\right],$$

puis, pour chaque $n \in \mathbf{N}^+$:

$$2(-1)^n\frac{y'_e(L)}{L}\frac{L^2(2\pi^2 n^2+15)}{6\pi^2 n^2} - \frac{A_n(2\pi^2 n^2+3)}{6} + (-1)^n\sum_{\substack{l=1\\l\neq n}}^{\infty}\left[A_l(-1)^{l+1}\frac{4l^2(n^2+l^2)}{\left(n^2-l^2\right)^2} + \frac{8L^2}{\pi^2}\frac{y'_e(L)}{L}\frac{\left(n^2+l^2\right)}{\left(n^2-l^2\right)^2}\right] +$$

$$+ a\frac{A_n}{2} - a\sum_{\substack{l=1\\l\neq n}}^{\infty}(-1)^{l+n}\frac{2l^2 A_l}{n^2-l^2} + bA_n = P_n$$

$$\frac{(2\pi^2 n^2-3)}{6\pi n}\left[-\pi n B_n - 2(-1)^n y_o(L)\right] - 8n\sum_{\substack{l=1\\l\neq n}}^{\infty}\left[\frac{l^3 B_l(-1)^{n+l}}{\left(l^2-n^2\right)^2}\right] - \frac{16nL}{\pi}\frac{y_o(L)}{L}\sum_{\substack{l=1\\l\neq n}}^{\infty}\frac{(-1)^l l^2}{\left(l^2-n^2\right)^2} +$$

$$- a\frac{B_n}{2} + 2an(-1)^{n+1}\sum_{\substack{l=1\\l\neq n}}^{\infty}\frac{lB_l(-1)^l}{n^2-l^2} + bB_n = Q_n$$

### 2.3. ODEs non linéaires

La détermination de la solution en calculant les coefficients de son développement en série de Fourier sinusoïdale peut être aussi appliquée avec succès aux ODEs non linéaires, avec des non-linéarités polynomiales, des équations dont les termes ne contiennent pas que des puissances naturelles de l'inconnu, ses dérivés et leurs produits. Les coefficients de ces termes peuvent être constants, ou variables. À partir du développement

$$\bar{y}(x) = y_0 + \sum_{n=1}^{\infty} A_n \cos\frac{n\pi x}{L} + \sum_{n=1}^{\infty} B_n \sin\frac{n\pi x}{L},$$

les relations de type (2.2)–(2.3), déduites pour le calcul des coefficients de Fourier des dérivées de tout ordre et de type (1.9–(1.10), pour le calcul des coefficients de Fourier du produit des deux certains fonctions, suffisent pour convertir l'équation donnée en un système de *2N+1* équations algébriques (*N→∞*), avec *2N+1* inconnues: *y₀, A₁, A₂, ..., A_N, B₁, B₂, ..., B_N*
Par exemple:

$$\left(y^2\right)^{SFS} = y_0^2 + \frac{1}{2}\sum_{l=1}^{\infty}\left(A_l^2 + B_l^2\right) + \left[2y_0 A_n + \frac{1}{2}\sum_{l=1}^{\infty} A_l\left(A_{n+l} + A_{|n-l|}\right) + \frac{1}{2}\sum_{l=1}^{\infty} B_l\left(B_{n+l} - \text{sgn}(l-n)B_{|n-l|}\right)\right]\cos\frac{n\pi x}{L} +$$

$$+ \left[2y_0 B_n + \frac{1}{2}\sum_{l=1}^{\infty} A_l\left(B_{n+l} - \text{sgn}(l-n)B_{|n-l|}\right) - \frac{1}{2}\sum_{l=1}^{\infty} B_l\left(A_{n+l} - A_{|n-l|}\right)\right]\sin\frac{n\pi x}{L}$$

$$\left(yy'\right)^{SFS} = y_0\frac{y_o(L)}{L} + 2y_0 A_n \cos\frac{n\pi x}{L} +$$

$$+ \sum_{n=1}^{\infty}\sum_{l=1}^{\infty}\frac{1}{2}\left\{A_l\left[\omega_{n+l}B_{n+l} + \omega_{|n-l|}B_{|n-l|} + 4(-1)^{l+n}\frac{y_o(L)}{L}\right] - B_l\left[\omega_{n+l}A_{n+l} - \text{sgn}(l-n)\omega_{|n-l|}A_{|n-l|}\right]\right\}\cos\frac{n\pi x}{L} -$$

$$- \sum_{n=1}^{\infty}\sum_{l=1}^{\infty}\left\{y_0\omega_n A_n + \frac{y_o(L)}{L}B_n - \frac{1}{2}\sum_{l=1}^{\infty} A_l\left[\omega_{n+l}A_{n+l} - \text{sgn}(l-n)\omega_{|n-l|}A_{|n-l|}\right] - \frac{1}{2}\sum_{l=1}^{\infty} B_l\left[\omega_{n+l}B_{n+l} - \omega_{|n-l|}B_{|n-l|}\right]\right\}\sin\frac{n\pi x}{L}$$

$$\left(y'^2\right)^{SFS} = \frac{y_o(L)^2}{L^2} + \frac{1}{2}\sum_{l=1}^{\infty}\left[\omega_l^2 A_l^2 + \omega_l^2 B_l^2 + 4(-1)^n B_l\frac{y_o(L)}{L} + 4\frac{y_o(L)^2}{L^2}\right] + 2y_0\sum_{n=1}^{\infty}\left[\omega_n B_n + 2(-1)^n\frac{y_o(L)}{L}\right]\cos\frac{n\pi x}{L} +$$

$$+ \sum_{n=1}^{\infty}\sum_{l=1}^{\infty}\frac{1}{2}\left\{\left[\omega_l B_l + 2(-1)^l\frac{y_o(L)}{L}\right]\left[-\omega_{n+l}B_{n+l} - \omega_{|n-l|}B_{|n-l|} + 4(-1)^{n+l}\frac{y_o(L)}{L}\right] + \omega_l A_l\left[\omega_{n+l}A_{n+l} - \text{sgn}(l-n)\omega_{|n-l|}A_{|n-l|}\right]\right\}\cos\frac{n\pi x}{L} +$$

$$+ \sum_{n=1}^{\infty}\left\{-2y_0\omega_n A_n + \frac{1}{2}\sum_{l=1}^{\infty}\left[\omega_l B_l + 2(-1)^l\frac{y_o(L)}{L}\right]\left(-\omega_{n+l}A_{n+l} + \text{sgn}(l-n)\omega_{|n-l|}A_{|n-l|}\right) - \omega_l A_l\left[\omega_{n+l}B_{n+l} - \omega_{|n-l|}B_{|n-l|}\right]\right\}\sin\frac{n\pi x}{L}$$

On peut voir que pour obtenir les valeurs des coefficients de Fourier de l'inconnu *y(x)*, il est nécessaire de résoudre certaines équations algébriques du second ordre ou plus, si le degré de non-linéarité augmente. On peut également remarquer que si, dans l'un des types d'équations analysés, les substitutions opérées pour résoudre l'équation sont faites en utilisant les sommes partiels *S_N* (avec *2N+1* termes) au lieu du développement entier (avec un nombre infini de termes), on obtienne une ODE approximative, l'approximation étant d'autant meilleure que le *N* est plus grand. Pour le résoudre, il faut résoudre un système de *2N+1* équations algébriques, problème pour lequel il existe des nombreux programmes de calcul électronique. Par conséquent, il est possible de développer un algorithme simple, universel et précis pour résoudre ces types d'ODEs.

### 3. Méthodes de résolution d'équations différentielles

Les exemples examinés dans la section précédente ont souligné qu'en remplaçant certains termes des ODE par leurs développements en SFS/SFN, l'équation initiale reste, ou peut devenir, linéaire. Dans le cas général, la fonction inconnue *y(x)* est remplacée par une somme infinie:

$$\hat{y}(x) = y_0 \cdot 1 + \sum_{n=-\infty}^{\infty} C_n R_n(x), \text{ sau } \hat{y}(x) = y_0 \cdot 1 + \sum_{n=1}^{\infty} A_n P_n(x) + \sum_{n=1}^{\infty} B_n Q_n(x)$$

qui peut être écrit aussi $y(x) = y_0 + \sum_{n=-\infty}^{\infty} y_n(x)$, respectivement $y(x) = y_0 + \sum_{n=1}^{\infty} y_{en}(x) + \sum_{n=1}^{\infty} y_{on}(x)$

Si tous les termes de l'équation donnée peuvent être remplacés par des combinaisons des fonctions $R_n(x)$, respectivement $P_n(x)$ et $Q_n(x)$, une superposition d'équations indépendantes est obtenue qui, en raison de la linéarité de l'équation, est réduite à un système d'équations algébriques (les fonctions $y_n(x)$, respectivement $y_{en}(x)$ et $y_{on}(x)$ sont parfaitement déterminées en déterminant les coefficients $C_n$, respectivement $A_n$ și $B_n$).

Les exemples analysés (2.1-2.7) ont également mis en évidence le fait que cette méthode permet de déterminer, dans de nombreux cas, une solution exacte, exprimée sous une forme fermée. Dans des nombreuses autres situations (par exemple, pour les ODEs à coefficients variables), on obtienne des solutions exactes, exprimées sous une forme des sommes infinies de termes. Pour être utiles en pratique, ces formes peuvent être approximées en négligeant les termes les moins significatifs: **l'approximation de la solution**. La solution approximative est déterminée à tous les points de l'intervalle, sans interpolation requise, comme dans des nombreuses autres méthodes approximatives courantes.

Dans de nombreuses autres situations, pour résoudre l'équation, l'approximation du développement en série doit se faire **en avant** d'arriver à la solution exacte: l'impossibilité de résoudre le système d'équations algébriques de dimensions infinies, résulté par le traitement d'ODE, l'impossibilité de déterminer par des calculs analytiques les coefficients de Fourier, etc. On peut considérer que nous avons eu recours à une **approximation de la méthode**. Dans ces situations, on obtienne aussi une solution exprimée sous une forme fermée (une somme finie des premières $2N+1$ harmoniques). Dans ce cas aussi, les valeurs de la solution sont parfaitement déterminées sur toute l'intervalle de définition.

Il existe également des situations dans lesquelles le remplacement des certains termes d'ODEs se fait par d'autres expressions qui les rapprochent, sur la base des critères empiriques, méthodologiques ou autres. Ce sont des méthodes d'**approximation de l'équation**. Dans la plupart des cas, cette approximation est faite en un nombre fini de points, la valeur de la solution aux points intermédiaires étant déterminée par interpolation.

## 4. Solutions des ODEs, déterminées en calculant les coefficients de son développement en séries de Fourier non sinusoïdales

Dans l'article [1], nous avons proposé une généralisation du développement en SFS des fonctions périodiques définies sur un intervalle $[-L, L]$, généralisation par laquelle nous avons remplacé la base sinusoïdale périodique (le nom usuel est base trigonométrique) avec une base non sinusoïdale, comprenant les fonctions suivantes: la fonction unité $I$, les quasi-harmoniques fondamentales $g(x)$-paire et $h(x)$-impaire, périodiques, avec la période $2L$, avec des valeurs moyenne nulles sur l'intervalle de définition et les quasi-harmoniques secondaires, définis sur $[-L, L]$, $g_n(x) = g(nx)$ et $h_n(x) = h(nx)$, avec la période $2L/n$, pour $n \in \mathbf{Z}^+$. Les quasi-harmoniques fondamentales $g(x)$ et $h(x)$, définies sur l'intervalle $[-L, L]$, peuvent être toutes fonctions admettant des développements en série sinusoïdale, étendus sur l'axe réel, par des translations successives: pour tout $x_R \in R$, on définit la fonction $K(x_R) = E[(x_R - x_1)/2L]$, tel que pour chaque $x_R \in R$ et chaque $x \in [-L, L]$, il y a les relations $x_R = x + KT$ et $g(x_R) = g(x_R - KT) = g(x)$. $E(x) = \lfloor x \rfloor = k$, est la fonction **partie entière** ($k$ est l'entier le plus proche inférieur ou égal à $x$, c'est à dire $E(x) \leq x < E(x)+1$. Dans la littérature anglaise, on utilise le nom de *floor function*). Les quasi-harmoniques secondaires sont obtenues par la dilatation de ces fondamentales: $g_n(x) = g_n[x - 2L \cdot E[(x+L)n/2L]$. La fonction $g_n(x)$ reçoit sur l'intervalle [*2L(2k–1)/n, 2L(2k+1)/n*], les mêmes valeurs que celles qu'il reçoit *g(x)* sur l'intervalle [*–L, L*]. Nous avons introduit pour cette fonction la notation simplifiée:

$g_n(x) = G[-L/n < g(nx) > L/n]_n$, où $n \in \mathbf{N}^+$.

Les coefficients $A_n$ et $B_n$ du développement en SFN de la fonction *f(x)* sont obtenus au moyen des relations algébriques entre les coefficients de Fourier des développements en SFS des fonctions *f(x)*, *g(x)* et *h(x)*.

L'article [1] mis en évidence un large éventail d'applications possibles de ce nouveau type de développement en série, en soulignant la diversité des solutions possibles pour chacun d'entre eux, illustrée par les différentes façons d'approximer les fonctions. Concernant la résolution des équations différentielles, dans la section précédente, nous avons analysé une nouvelle méthode de résolution, en remplaçant la fonction inconnue par son développement en SFS et en résolvant les équations algébriques résultantes. La méthode est également applicable aux développements en SFN, à condition qu'à travers les opérations de dérivation/ intégration aucun élément d'une autre base de développement n'apparaisse (donc, seulement pour les bases composées des fonctions circulaires, hyperboliques, exponentielles). En utilisant cette méthode, peuvent également être capitalisés les avantages offerts par la possibilité d'utiliser des bases orthogonalisées, bien que cela implique un volume de calcul supplémentaire.

## 4.1. Méthodes de résolution d'équations différentielles par l'approximation de la solution

Dans certains ODEs, il n'est pas nécessaire de remplacer tous les termes de l'équation différentielle par leurs développements en séries de Fourier sinusoïdales/non sinusoïdales. L'équation qui résulte par le remplacement seulement des certains termes, même lorsqu'elle n'est pas linéaire, peut conduire à une équation plus facile à résoudre que l'équation initiale, en utilisant des procédures classiques d'intégration/dérivation (la méthode est fréquemment appliquée dans la pratique actuelle par des techniques de linéarisation d'équations non linéaires). Nous illustrerons ces énoncés en résolvant une équation inhomogène très simple:

$$y' = x \tag{4.1}$$

avec la condition *y(0)=y₀*

Nous sélectionnerons pour commencer une série des développements en série de Fourier [1], qui seront utiles:

$$\bar{x} = X[-L < x > L] = \frac{2L}{\pi} \sum_{n=1}^{\infty} \frac{(-1)^{n+1}}{n} \sin\frac{n\pi x}{L} \tag{4.2}$$

$$(x^2)^{SFS} = X_e^2[-L > x^2 < L] = \frac{L^2}{3} - \frac{4L^2}{\pi^2} \sum_{n=1}^{\infty} \frac{(-1)^{n+1}}{n^2} \cos\frac{n\pi x}{L} \tag{4.3}$$

$$(1_o)^{SFS} = J[-L > -1 < 0 > 1 < L] = \frac{4}{\pi} \sum_{n=1}^{\infty} \frac{1}{2n-1} \sin\frac{(2n-1)\pi x}{L} \tag{4.4}$$

$$(x^{ir})^{SFS} = X_e\left[-L > -x - \frac{L}{2} < 0 > x - \frac{L}{2} < L\right] = -\frac{4L}{\pi^2} \sum_{n=1}^{\infty} \frac{1}{(2n-1)^2} \cos\frac{(2n-1)\pi x}{L} \tag{4.5}$$

A partir des deux premières relations, nous pouvons écrire les développements en séries non sinusoïdales:

$$\hat{x} = X[-L < x > L] = \frac{L}{2}\left(J_1 - \sum_{n=1}^{\infty} \frac{1}{2^n} J_{2^n}\right) =$$

$$= \frac{L}{2} J[-L > -1 < 0 > 1 < L] - \frac{L}{2} \sum_{n=1}^{\infty} \frac{1}{2^n}\left[-\frac{L}{2^n} > -1 < 0 > 1 < \frac{L}{2^n}\right]_{2^n} \tag{4.6}$$

$$\left(x^2\right)^{SFN} = X_e^2\left[-L > x^2 < L\right] = \frac{L^2}{3} + L \cdot X_{e1} - \sum_{n=1}^{\infty} \frac{L}{4^n} \cdot \left[X_e\right]_{2^n} =$$

$$= \frac{L^2}{3} + L\left[-L > -x - \frac{L}{2} < 0 > x - \frac{L}{2} < L\right]_1 - \sum_{n=1}^{\infty} \frac{L}{4^n}\left[-\frac{1}{2^n} > -2^n x - \frac{L}{2} < 0 > 2^n x - \frac{L}{2} < \frac{1}{2^n}\right]_{2^n} \qquad (4.7)$$

Pour l'intervalle $[-L, L]$ le terme inhomogène de l'équation (4.1) peut être développé dans une série sinusoïdale, ou dans une série non sinusoïdale (pour la dernière variante les possibilités étant multiples), ce qui fait qu'en utilisant la notation réduite, on peut écrire:

$$\left(y'\right)^{SFS} = \frac{2L}{\pi} \sum_{n=1}^{\infty} \frac{(-1)^{n+1}}{n}\left[-\frac{L}{n} < \sin\frac{n\pi x}{L} > \frac{L}{n}\right]_n, \text{ respectivement} \qquad (4.8)$$

$$\left(y'\right)^{SFN} = \frac{L}{2}\left[-L < -1 > 0 < 1 > L\right] - \frac{L}{2}\sum_{n=1}^{\infty} \frac{1}{2^n}\left[-\frac{L}{2^n} < -1 > 0 < 1 > \frac{L}{2^n}\right]_{2^n} \qquad (4.9)$$

Par intégration:

$$y + C = -\frac{2L^2}{\pi^2} \sum_{n=1}^{\infty} \frac{(-1)^{n+1}}{n^2}\left[-\frac{L}{n} < \cos\frac{n\pi x}{L} > \frac{L}{n}\right]_n = \frac{1}{2}\left(X_e^2\left[-L > x^2 < L\right] - \frac{L^2}{3}\right) = \frac{x^2}{2} - \frac{L^2}{6},$$

$$y + C = \frac{L}{2}\left[-L < -x > 0 < x > L\right] - \frac{L}{2}\sum_{n=1}^{\infty} \frac{1}{2^n}\left[-\frac{L}{2^n} < -x > 0 < x > \frac{L}{2^n}\right]_{2^n} =$$

$$= \frac{L}{2}\left(\left[-L < -x - \frac{L}{2} > 0 < x - \frac{L}{2} > L\right] + \frac{L}{2}\right) - \frac{L}{2}\sum_{n=1}^{\infty}\left(\frac{1}{4^n}\left\{\left[-\frac{L}{2^n} < -2^n x - \frac{L}{2} > 0 < 2^n x - \frac{L}{2} > \frac{L}{2^n}\right]_{2^n} + \frac{L}{2}\right\}\right) =$$

$$= \frac{1}{2}\left(X_{e1} - \sum_{n=1}^{\infty} \frac{1}{4^n}\left[X_e\right]_{2^n}\right) + \frac{L}{2}\left(\frac{L}{2} - \sum_{n=1}^{\infty}\frac{1}{4^n}\frac{L}{2}\right) = \frac{1}{2}\left(x^2 - \frac{L^2}{3}\right) + \frac{L^2}{4}\left(1 - \frac{1}{3}\right) = \frac{x^2}{2} - \frac{L^2}{6} + \frac{L^2}{6} = \frac{x^2}{2}$$

Dans tout les deux variantes, la solution particulière est $y - y_0 = x^2/2$.

On remarque que la dérivée de la fonction paire $y_e(x)$ est une fonction impaire qui peut être développée en série de fonctions impaires, avec la valeur moyenne nulle, avec la base $h(x)$:

$$y_e'(x) = \sum_{n=1}^{\infty} A_n h_n(x)$$

L'intégration de cette relation de dérivation conduit à une relation pour le développement de la fonction paire $y_e(x)$, dans une base $g(x)$ des fonctions paires, pour laquelle la valeur moyenne $g_0$ peut être non nulle:

$$y_e(x) = y_0 + \sum_{n=1}^{\infty} A_n(g - g_0)_n, \text{ où } y_0 = g_0\sum_{n=1}^{\infty} A_n \qquad (4.10)$$

On remarque également que dans l'intervalle $[-L, L]$, chacune des équations (4.8) et (4.9) est obtenu par la superposition linéaire des sous-équations suivantes:

$$y_n' = \frac{2L}{\pi}\frac{(-1)^{n+1}}{n}\sin\frac{n\pi x}{L}, \text{ respectivement } y_0' = \frac{L}{2}J_1 \text{ et } y_n' = \frac{L}{2}\frac{1}{2^n}J_{2^n}, \text{ pour } n=1, 2, ...,\infty$$

En résolvant ces équations partielles et en additionnant les solutions, nous obtenons les mêmes résultats que ceux obtenus en intégrant toute l'équation. Dans le cas présent, les dérivées $y'_n(x)$ des solutions partielles $y_n(x)$ sont même les harmoniques (ou quasi-harmoniques) du développement en une série connue de la fonction $y'$.

Par conséquent, la méthode peut être appliquée à toutes les équations dans lesquelles le développement en série d'un terme conduit à une équation linéaire représentée comme une superposition infinie d'équations partielles pour lesquelles la solution peut être trouvée par des méthodes de résolution classiques.

## 4.2. La linéarisation de l'équation différentielle du pendule gravitationnel, en introduisant une somme infinie de fonctions rampe

Nous appliquerons la méthode de linéarisation des équations non linéaires à l'équation du **pendule gravitationnel**:

$$\frac{d^2\theta}{dt^2} + \omega_0^2 \sin\theta = 0 \quad , \quad \theta \in \left[-\pi, \pi\right] \tag{4.11}$$

Dans une approche classique [13], l'équation est résolue par la méthode d'approximation de l'équation, en remplaçant $\sin\theta \approx \theta$, acceptable pour $\theta \in \left[-\theta_1, \theta_1\right]$, si $\theta_1 \approx 0$. Ici, nous allons essayer de trouver une solution exacte, valable pour tout l'intervalle $\left[-\pi, \pi\right]$. Pour cela, nous allons recourir au développement en série non sinusoïdale du terme $sin\theta$, la base du développement étant la fonction impaire
$g_1(\theta) = X_o[-\pi < -\theta - \pi > -\pi/2 < \theta > \pi/2 < -\theta + \pi > \pi]$.

Le développement en série de Fourier de la fonction $g_1(\theta)$ nous amène à:

$$\overline{g}_1(\theta) = \frac{4}{\pi} \sum_{n=1}^{\infty} \frac{(-1)^{n+1} \sin(2n-1)\theta}{(2n-1)^2} \qquad \text{dont les coefficients sont:}$$

$d_1 = 4/\pi$, $d_2 = 0$, $d_3 = -4/9\pi$, $d_4 = 0$, $d_5 = 4/25\pi$, $d_6 = 0$, $d_7 = -4/49\pi$, $d_8 = 0$, $d_9 = 4/81\pi$, $d_{10} = 0$, $d_{11} = -4/121\pi$, $d_{12} = 0$, ... Pour la fonction $f(\theta) = sin\theta$, nous pouvons écrire:

$$\sin\theta = \hat{f}(\theta) = \sum_{n=1}^{\infty} B_n g_n(\theta), \text{ dans lequel:}$$

$$B_1 = \frac{\pi}{4}, B_2 = 0, B_3 = \frac{\pi}{36}, B_4 = 0, B_5 = -\frac{\pi}{100}, B_6 = 0, B_7 = \frac{\pi}{196}, B_8 = 0, B_9 = 0, ...$$

$$B_{2n-1} = \frac{\pi}{4} \cdot \frac{(-1)^n}{(2n-1)^2}, B_{2n} = 0, \text{ mais } B_{n^2} = 0, \text{ pour } n = 2, 3, ..., \infty. \text{ On obtient l'équation linéaire:}$$

$$\frac{d^2\theta}{dt^2} + \omega_0^2 \sum_{n=1}^{\infty} B_n g_n(\theta) = 0$$

laquelle, pour $\theta = \sum_{n=1}^{\infty} \theta_n$, est la superposition linéaire des équations partielles:

$$\frac{d^2\theta_n}{dt^2} + \omega_0^2 B_n g_n(\theta_n) = 0 \text{ , où } n = 1, 2, 3, ..., \infty \tag{4.12}$$

Si $\theta(0) = \theta_0$, et $v_0 = 0$, ces conditions initiales sont valables pour toutes les équations partielles: $\theta_n(0) = \theta_0$, et $v_{0n} = 0$.

Si $v_0 \neq 0$, le problème donné est remplacé par un problème équivalent: nous considérons la vitesse angulaire initiale $v_0$ imposé au pendule, comme provenant de la transformation en énergie cinétique $E_c = 1/2 \cdot mv^2$ de la partie $E_p = mg(h_e - h_0)$ de l'énergie potentielle du pendule de longueur $l$, située à une hauteur équivalente $h_e$ (ce qui correspond à un angle équivalent $\theta_e$). Ainsi, le problème initial est remplacé par celui dans lequel la vitesse angulaire initiale $v_0$ est nulle et la position initiale est:

$$\theta_e = 2k\pi \pm \arccos\left(\cos\theta_0 - \frac{lv_0^2}{2g}\right) \tag{4.13}$$

Des valeurs $k \neq 0$ apparaissent quand $h_e > 2l$. Pour ces cas, nous pouvons calculer $v_\pi$, la vitesse à laquelle le pendule traverse la position $h = \pm\pi$. Si $\theta_e$ et $\theta_0$ appartiennent à l'intervalle $(-\pi, \pi)$, les solutions partielles équivalentes $\theta_n(t)$ ont tous, les conditions initiales:

$$\theta_n(0) = \theta_e \text{ et } \frac{d\theta_n}{dt}\bigg|_{t=0} = 0 \text{ , et pour } h_e > 2l: \theta_n(0) = \pi \text{ et } \frac{d\theta_n}{dt}\bigg|_{t=0} = v_\pi \tag{4.14}$$

La résolution de ces *n* équations partielles nous amène à trouver les *n* réponses partielles du système aux ces *n* forces périodiques agissant sur le pendule. La somme de ces réponses est la solution de l'équation pour le problème équivalent. Cette trajectoire-solution $\theta(t)=\Sigma\theta_n(t)$ passe par la position $\theta_0$ à l'instant $t_e$. La condition initiale au problème initial devient $\theta_0=\theta(t-t_e)=\Sigma\theta_n(t-t_e)$. Les positions $\theta_n(t_e)$ et les vitesses $\theta_n{}'(t_e)$ des trajectoires particulières sont enregistrées lorsque le pendule équivalent passe par la position $\theta_0$. De cette façon, l'énergie potentielle équivalente du pendule (correspondant à la vitesse initiale $v_0$) est distribuée sur les composants partiels du système de forces total et transformée par chaque pendule équivalent selon sa propre spécificité, en énergie cinétique, résultant une vitesse angulaire initiale équivalente pour chaque équation partielle

En raison des lois de conservation de l'énergie, la trajectoire du pendule, pour chaque système partiel de forces $F_n$ , est symétrique par rapport à une position d'équilibre $\theta_{mn}$ $(\theta_0)$ (pour *n=1*). Elle passe par les positions $\theta_m-\theta_i$ avec les mêmes vitesses $v_i$ avec lequel il passe à travers les positions $\theta_m+\theta_i$, et atteint la position maximale $-\theta_e$ également à vitesse nulle. Si $\theta_e{>}\pi$, la vitesse équivalente $v_\pi$ que correspond à la position $\theta{=}\pi$ est égal à celui de la position $\theta{=}-\pi$. Pour les mêmes raisons, les trajectoires partielles avec des vitesses initiales $\theta_e$ negative sont symétriques par rapport à l'axe $\theta{=}0$ avec ceux dont la vitesse initiale $\theta_e$ est positive.

De  (4.12), pour *n=1*, on obtient l'équation partielle

$$\frac{d^2\theta_1}{dt^2} + \omega_0^2\,\frac{\pi}{4}\,g_1(\theta) = 0 \quad\text{, où } g_1(\theta)=X_o[-\pi{<}-\theta{-}\pi{>}-\pi/2{<}\theta{>}\pi/2{<}-\theta{+}\pi{>}\pi],$$

qui se résout successivement, pour chaque sous-intervalle dans lequel la force exercée sur le pendule a une certaine loi de variation. Pour $-\pi/2{<}\theta_e{<}\pi/2$, sur le sous-intervalle $[-\pi/2,\ \pi/2]$ la force augmente uniformément d'une valeur négative à une valeur positive, l'équation est harmonique et la solution est:

$$\theta_1(t) = \theta_e\cos\!\left(\sqrt{\pi/4}\,\omega_0 t\right) \tag{4.15}$$

Pour $-\pi{<}\theta_e{<}-\pi/2$, l'équation partielle prend des formes différentes pour des sous-intervalles différents. Sur le sous-intervalle $[\theta_e,\ -\pi/2]$ l'équation prend la forme:

$$\frac{d^2\theta_1}{dt^2} - \omega_0^2\,\frac{\pi}{4}\cdot(\theta_1+\pi) = 0 \ , \tag{4.16}$$

avec les conditions initiales: $\theta_{1.0}=\theta_1(0) = \theta_e$  et  $v_{1.0}=\dfrac{d\theta_1}{dt}\Big|_{t=0} = 0$  et avec la solution:

$$\theta_1(t) = (\theta_e+\pi)ch\!\left(\sqrt{\pi/4}\cdot\omega_0 t\right) - \pi \tag{4.17}$$

$$v_1(t) = (\theta_e+\pi)\omega_0\sqrt{\pi/4}sh\!\left(\sqrt{\pi/4}\cdot\omega_0 t\right), \tag{4.18}$$

Les équations (4.16) et (4.18) permettent de trouver le moment $t_{1.1}$ de commutation de la champ des forces, c'est-à-dire le moment où $\theta_1{=}-\pi/2$, puis la vitesse angulaire du pendule à partir de ce moment: $v_{1.1}=v_1(t_1)=(\theta_e+\pi)\omega_0\sqrt{\pi/4}sh\!\left(\sqrt{\pi/4}\cdot\omega_0 t_1\right)$. Ces valeurs deviennent les conditions initiales de l'équation harmonique dans le sous-intervalle $[-\pi/2,\ \pi/2]$, valable pour $t{>}t_{1.1}$. La solution est:

$$\theta_1(t-t_{1.1}) = (\pi/2)\cos\!\left(\sqrt{\pi/4}\,\omega_0(t-t_{1.1})\right) + \left(v_{1.1}/\omega_0\sqrt{\pi/4}\right)\sin\!\left(\sqrt{\pi/4}\cdot\omega_0(t-t_{1.1})\right) \tag{4.19}$$

Maintenant, de $\theta_1(t){=}0$, on peut calculer l'instant $t_{1.2}$ du passage par le point d'équilibre. Plus loin, l'évolution du pendule est symétrique par rapport à ce point. Le pendule atteint la position $\theta_1(t){=}-\theta_e$ après le temp $2t_{1.2}$ et continuera une trajectoire symétrique par rapport à l'axe $t{=}2t_{1.2}$, pour arriver après un temps total $t_t{=}4t_{1.2}$ de nouveau en position $\theta_1(t){=}\theta_e$ et pour continuer une trajectoire périodique avec la période $T{=}4t_{1.2.}$

La figure 2.a. montre graphiquement les solutions pour la première harmonique dans le cas des vitesses initiales nulles, pour des positions initiales supérieures et inférieures à $-\pi/2$, et la figure 2.b., pour des vitesses initiales non nulles.

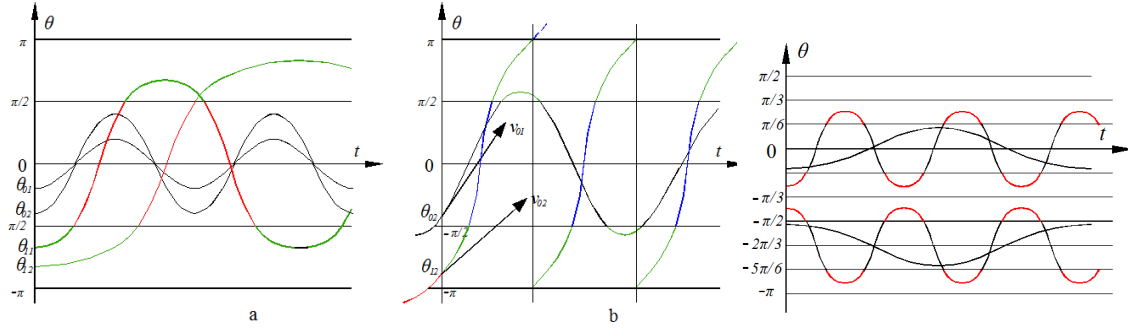

*Fig.2. **a**: des solutions partielles de premier ordre pour $v_0=0$, $-\pi/2 < \theta_{01}$, $\theta_{02} < 0$ et $-\pi < \theta_{11}$, $\theta_{12} < -\pi/2$. **b**: pour $v_0 \neq 0$ **c**: des solutions partielles de troisième ordre pour $v_0=0$ şi $\theta_e < 0$*

Pour $n=3$ on a l'équation $\quad \dfrac{d^2\theta_3}{dt^2} + \omega_0^2 \dfrac{\pi}{4}\dfrac{1}{9}g_3(\theta) = 0$, où

$g_3(\theta) = X_0[-\pi/3 < 3(-\theta - \pi/3) > -\pi/6 < 3\theta > \pi/6 < 3(-\theta + \pi/3) > \pi/3]_3 = X_0[-\pi < -3\theta - 3\pi > -5\pi/6 < 3\theta + 2\pi > -\pi/2 < -3\theta - \pi > -\pi/6 < 3\theta > \pi/6 < -3\theta + \pi > \pi/2 < 3\theta - 2\pi > 5\pi/6 < -3\theta + 3\pi > \pi]$

Les conditions initiales de l'équation sont données par (4.14).

Pour $v_0=0$, si $-\pi/6 < \theta_0 < \pi/6 \rightarrow g_3(\theta)=3\theta$, et la solution est: $\theta_3(t) = \theta_0 \cos\left(\sqrt{\pi/12}\,\omega_0 t\right)$,

Pour $v_0 \neq 0$, si $|\theta_e| = \sqrt{\theta_0^2 + \left(v_0 / \omega_0 \sqrt{\pi/12}\right)^2} \leq \pi/6$, alors

$$\theta_3(t) = \theta_0 \cos\left(\sqrt{\dfrac{\pi}{12}}\omega_0 t\right) + \dfrac{v_0}{\omega_0\sqrt{\pi/12}}\sin\left(\sqrt{\dfrac{\pi}{12}}\cdot\omega_0 t\right) = \theta_e \cos\left[\sqrt{\dfrac{\pi}{12}}\cdot\omega_0 t + arctg\left(-v_0/\theta_0\omega_0\sqrt{\dfrac{\pi}{12}}\right)\right]$$

Pour les autres sous-intervalles dans lesquels $\theta_e$ peut être situé, la solution partielle est:

$$\theta_3(t) = (\theta_e - \theta_m)\cos\left(\sqrt{\dfrac{\pi}{12}}\omega_0 t\right) - \theta_m, \quad \text{ou} \quad \theta_3(t) = (\theta_e - \theta_m)ch\left(\sqrt{\dfrac{\pi}{12}}\omega_0 t\right) - \theta_m, \quad \text{comme la pente de la}$$

fonction $g_3(\theta)$ est positif, respectivement négatif. Ici, $\theta_m$ est le point médian de ce sous-intervalle et il a toujours le même signe que $\theta_e$.

Par conséquent, la cosinusoïde $\theta_3(t)$ (ainsi que ceux d'ordre supérieur) est un mouvement oscillant vis à vis du point $\theta_m$, avec l'amplitude $\theta_e - \theta_m$ (elle n'arrive jamais au point d'équilibre $\theta=0$ ni dans $\theta_0$, s'il est en dehors du sous-intervalle). La cosinusoïde hyperbolique $\theta_3(t)$ (ainsi que les quasi-cosinusoïdes de rang supérieur) décrit un mouvement divergent, d'éloignément de $\theta_m$ (qui est un point d'équilibre instable).

La figure 2.c présente certaines de ces solutions partielles, pour différentes valeurs négatives de $\theta_e$ ($v_0=0$).

Pour $n=5$, l'équation partielle est $\dfrac{d^2\theta_5}{dt^2} - \omega_0^2 \dfrac{\pi}{4}\dfrac{1}{25}g_5(\theta) = 0$, où

$g_5(\theta) = X_0[-\pi/5 < 5(-\theta - \pi) > -\pi/10 < 5\theta > \pi/10 < 5(-\theta + \pi) > \pi/5]_5$

Elle est résolu de la même manière que les équations précédentes. Il convient de noter qu'en raison du coefficient négatif de cette quasi-harmonique, la position $\theta=0$ est une d'équilibre instable. Pour $-\pi/10 < \theta_e < \pi/10$, la trajectoire du pendule est divergente.

La solution générale de l'équation équivalente est la somme de toutes les solutions partielles: $\theta(t) = \Sigma\theta_n(t)$. C'est aussi la solution exacte de l'équation donnée, exprimée comme une somme infinie de solutions partielles, chacune de ces solutions ayant des expressions différentes sur des sous-intervalles différents. Pour que cette solution soit pratiquement utile, nous ne retiendrons que les $N$ premières solutions partielles (l'approximation de la solution). Le même résultat est atteint en approximant la méthode, si à partir du développement en SFN de la fonction $\sin\theta$ on retient les premiers $N$ termes.

Pour une certaine valeur $t=t_e$, la somme $\theta(t)=\Sigma\theta_n(t)$ est nulle: $\Sigma\theta_n(t_e)=0$. Comme le montre la Fig.3 (pour simplifier nous avons choisi le cas $\theta_0=0$), à ce stade, chaque solution partielle a une valeur $\theta_n(t_e)=\theta_{n0}$ et une vitesse $v_n(t_e)=v_{n0}$. Ce sont les conditions initiales pour les composantes de l'équation initiale.

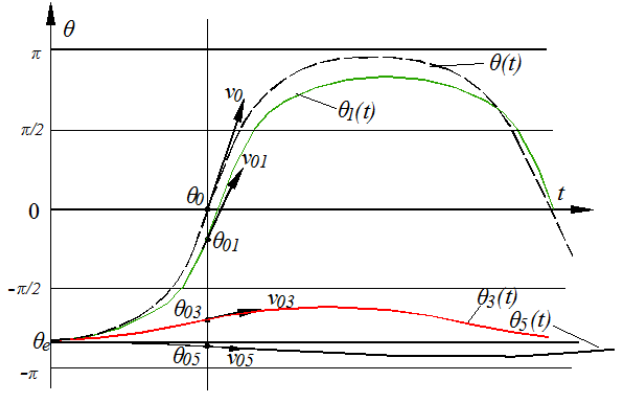

*Fig.3. La solution générale et les solutions partielles d'ordre 1, 3 et 5 pour $v_0\neq0$*

## 4.3. La linéarisation de l'équation différentielle du pendule gravitationnel, en introduisant une somme infinie de fonctions rectangulaires

Le développement en séries de Fourier non sinusoïdal que nous avons choisi pour la linéarisation de l'équation (4.11) est le plus proche de l'approche traditionnelle, mais ce n'est qu'une des nombreuses possibilités offertes par la méthode que nous proposons. Une autre solution intéressante est le remplacement dans l'équation de la fonction $\sin\theta$ avec son développement en série infinie d'impulsions rectangulaires périodiques:

$$\frac{d^2\theta}{dt^2} + \omega_0^2\sum_{n=1}^{\infty} B_n\Pi_n(\theta) = 0 \ , \tag{4.20}$$

où la variable $\theta$ peut avoir des valeurs dans l'intervalle $[-\pi, \pi]$,
$\Pi_1(\theta)= \Pi[-\pi(-1)0(1)\pi]$, et $\Pi_n(\theta)= \Pi[-\pi/n(-1)0(1)\pi/n]_n$ .
Comme dans l'approche précédente, les conditions initiales du problème équivalent sont:

$$\theta_0 =\theta(0) = \theta_e \ \text{ et } \ v_0 = \frac{d\theta}{dt}\bigg|_{t=0} = 0 \ , \tag{4.21}$$

où $\theta_e$ est l'angle équivalent, correspondant à la hauteur équivalente $h_e$.

Soit $B_n$ les coefficients du développement $\sin\theta = \hat{f}(\theta) = \sum_{n=1}^{\infty} B_n\Pi_n(\theta)$

Parce que $\Pi_1 = \frac{4}{\pi}\sum_{n=1}^{\infty}\frac{\sin(2n-1)\theta}{2n-1}$, les coefficients du développement de la fonction $\sin\theta$ sont:

$$B_1 = \frac{\pi}{4}, B_2 = 0, B_3 = -\frac{\pi}{12}, \ B_4 = 0, B_5 = -\frac{\pi}{20}, B_6 = 0, B_7 = -\frac{\pi}{28}, B_8 = 0, \ B_9 = 0, B_{10} = 0, \ldots$$

$$B_{2n-1} = \frac{\pi}{4}\cdot\frac{-1}{2n-1} \ , \ B_{2n} = 0, \ \text{mais} \ B_{n^2} = 0, \qquad \text{pour } n= 2, 3, .., \infty$$

Tous les coefficients, à l'exception de ce de la quasi-harmoniques fondamentales, sont négatifs. L'équation linéaire (4.20) est la superposition d'un nombre infini d'équations du type:

$$\frac{d^2\theta_n}{dt^2} = -\omega_0^2 B_n \Pi_n(\theta)$$ , où $n=1, 2, 3, ...,\infty$, avec les conditions initiales (4.21). Les solutions

de ces équations, pour tout les sous-intervalles $[t_k , t_{k+1}]$ dans lequel $\Pi_n(\theta)=ct$ sont:

$\theta_n(t) = -\omega_0^2 B_n sgn \Pi_{nk}(t-t_k)^2/2 + v_{0k}(t-t_k) + \theta_{0k}$, pour $k=0, 1, 2, 3, ...,n, ...$     (4.22)

où $sgn\Pi_{nk}$ est le signe des fonctions $\Pi_n(\theta)$ pour le sous-intervalle $[\theta_k , \theta_{k+1}]$, $t_k$ sont les moments des commutations de la fonction $\Pi_n(\theta)$, $\theta_{0k}$ et $v_{0k}$ sont les conditions initiales pour l'équation valable sur ce sous-intervalle (les valeurs finales du sous-intervalle précédent). La vitesse du pendule pour chaque système des forces est uniformément accélérée:

$$v_n(t) = -\text{sgn}(\Pi_n) B_n \omega_0^2 (t - t_k) + v_{0k}$$

Pour $n=1$, le point $\theta=0$ c'est un point d'équilibre stable. Quelle que soit la position initiale du pendule, il tend à atteindre une position d'équilibre stable et oscille autour de cette position, décrivant une quasi-sinusoïde construite sur la base d'un polynôme du deuxième degré. Cette courbe aura les points d'inflexion situés sur l'axe $\theta=0$. Les oscillations auront l'amplitude $\theta_e$ et la période d'oscillation $T_1 = 4\sqrt{8|\theta_e|/\pi\omega_0^2}$ (dependente de $\theta_e$). Un pendule que départs de la position $\theta_{01}=\theta_e$ avec la vitesse $v_{01}=0$, arrive après le temps $t=T_1/4$ en la position $\theta_1=0$ (où la force qui agit sur lui change de sens) avec la vitesse $v_1 = -\omega_0\sqrt{|\theta_e|\pi/2}$. En raison de l'inertie, le pendule continue de se déplacer vers la position $\theta=-\theta_e$, où il arrive après le temps $t=T_1/4$, avec la vitesse $v_1=0$. Sous l'action du même système de forces, le pendule continue son mouvement dans la direction opposée et après un autre quart de période, il retrouve la position $\theta_1=0$, cette fois avec la vitesse $v_1 = \omega_0\sqrt{|\theta_e|\pi/2}$. Après un nouveau changement de direction de la force, après un autre quart de période, le pendule retrouve sa position $\theta_1=\theta_e$, avec la vitesse $v_1=0$. Les oscillations continuent avec la période $T_1$. Avec les notations $A=\pi\omega^2/8$ et $T = \sqrt{8|\theta_e|/\pi\omega_0^2}$ , on a:

$\theta_1(t) = X_1^2 = X^2 [0*(\theta_e-At^2)*T*(-\theta_e + A(t-T)^2) *3T*(\theta_e-A(t-3T)^2) *4T]_1$     (4.23)

Si les conditions initiales changent: $\theta_0=\theta_e$ et $v_1(0)\neq0$, la solution de l'équation change sa période, son amplitude et son déphasage initial.

La solution partielle d'ordre $n$ de l'équation (4.20) aux conditions initiales (4.21) est également un polynôme quasi-sinusoïdal de deuxième degré:     (4.24)

$\theta_n(t) = X_n^2 = X^2 [0*(-\theta_n+At^2/n)*Tn^{1/2}*(\theta_n - A(t-T/n)^2) *3 n^{1/2}*(-\theta_n+A(t-3T/n)^2) *4T n^{1/2}]_n$,

où $\theta_n = \theta_e - \theta_m$ . Cette courbe a tous les points d'inflexion situés sur l'axe $\theta=\theta_m$. L'oscillation a l'amplitude $\theta_e - \theta_m$ et la période d'oscillation $T_n = 4\sqrt{2|\theta_e - \theta_l|/B_n\omega_0^2}$ .

Dans la figure 4.a, nous avons représenté, pour l'harmonique fondamentale, trois de ces solutions partielles, pour différentes valeurs de la position initiale équivalente $\theta_e$: deux pour $\theta_e \in [-\pi,\pi]$ (ligne rouge) et un pour $\theta_e \notin [-\pi, \pi]$ (ligne noire en pointillés); avec une ligne verte en pointillés, nous avons représenté la continuation fictive de la quasi-sinusoïde en dehors de l'intervalle $[-\pi, \pi]$).

La solution générale de l'équation équivalente est la somme de toutes les solutions partielles: $\theta(t)=\Sigma\theta_n(t)$. Une approximation de la solution est obtenue en additionnant les $N$ premières solutions partielles, où $N$ doit être suffisamment grand pour obtenir une erreur satisfaisante. Pour une valeur $t=t_e$, la somme est nulle: $\Sigma\theta_n(t_e)=0$. Comme le montre la figure 4b, pour le moment $t_e$, chaque solution partielle fournit une valeur $\theta_n(t_e)=\theta_{n0}$ et une vitesse $v_n(t_e)=v_{n0}$. Ce sont les conditions initiales des composantes de l'équation d'origine. La solution de l'équation est $\theta(t)=\Sigma\theta_n(t-t_e)$.

Avec les notations de (4.23) et (4.24) nous pouvons avancer une nouvelle expression pour le mouvement du pendule:

**$\theta(t)=\Sigma X_n^2$**, pour $t > t_e$ .

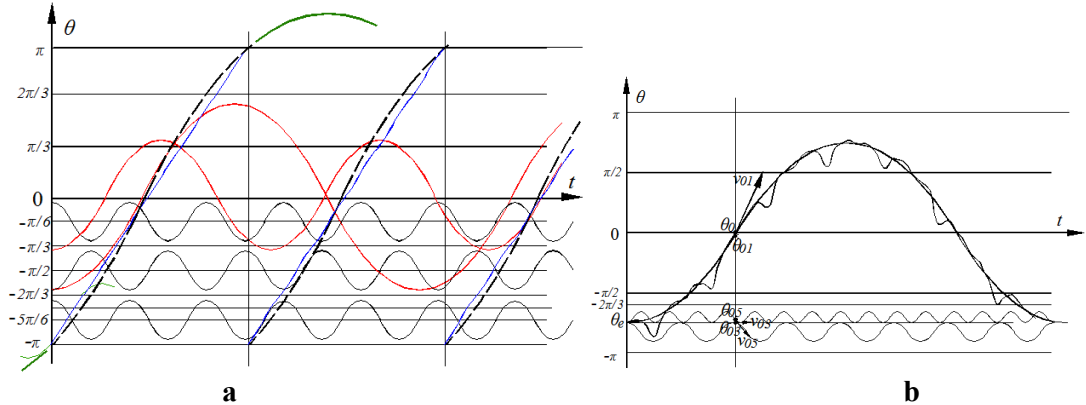

*Fig.4. **a**: Des solutions partielles pour $\theta_e < \pi$: rang 1 (ligne rouge), rang 3 (ligne noire).*
*Pour $\theta_e > \pi$: rang 1 (ligne pointillée), rang 3 (ligne bleue) b: solution générale*

## 5. Conclusions

Nous avons essayé ici, à l'aide de quelques exemples simples, de prouver que la méthode de résolution des équations différentielles en déterminant les coefficients du développement en séries sinusoïdales ou non sinusoïdales de la fonction inconnue est une méthode particulièrement solide, applicable à tous les types d'équations différentielles et intégro-différentielles, linéaires et non linéaires, d'équations aux dérivées partielles, des systèmes de telles équations, quels que soient leur ordre et quels que soient la complexité des coefficients.

## 6. Bibliographie

[1] Török A., Petrescu S., Feidt M., Séries de Fourier périodiques non sinusoïdales, https://hal.archives-ouvertes.fr/hal-02485085

[2] Nagle R. K., Saff E. B., Snider A. D., Fundamentals of Differential Equations and Boundary Value Problems, 8th ed, 2012, Pearson Education, Inc.

[3] O'Neil P. V., Advanced Engineering Mathematics, Seventh Edition, part 1, Cengage Learning, 2012, Publisher: Global Engineering: Christopher M. Shortt

[4] Kreyszig E., Advanced Engineering Mathematics, 10-th edition, chapters 1, 2, 5, 6, 11, 12, John Wiley & Sons, Inc., 2011

[5] Polyanin A. D., A. V. Manzhirov, Handbook of Mathematics for Engineers and Scientists, Chapman & Hall/CRC Press, 2007

[6] Polyanin A. D., A. V. Manzhirov, Handbook of Exact Solutions for Ordinary Differential Equations, CRC Press, 2003

[7] Zwillinger D., CRC Standard Mathematical Tables and Formulae, 31st Edition, 2003, by CRC Press Company

[8] Spiegel M. R., Mathematical Handbook of Formulas and Tables, Schaum's Outline Series, Schaum Publishing Co. 1st.ed. 1968, McGraw-Hill Book Company https://archive.org/details/MathematicalHandbookOfFormulasAndTables

[9] Dourmashkin, Classical Mechanics: MIT 8.01 Course Notes, Chapter 23 Simple Harmonic Motion, pp 8, http://web.mit.edu/8.01t/www/materials/modules/chapter23.pdf

[10]  Jordan D. W., Smith P., Nonlinear Ordinary Differential Equations. An introduction for Scientists and Engineers,4-th edition, Keele University, Oxford University Press, 2007,

[11]  Olver P. J., Equations Introduction to Partial Differential Equations, Springer Cham Heidelberg New York Dordrecht London, Springer International Publishing Switzerland, 2014, 6056 ISSN - 2197 5604 (electronic)

[12]  Dourmashkin, Classical Mechanics: MIT 8.01 Course Notes, Chapter 23 Simple Harmonic Motion, pp 8, http://web.mit.edu/8.01t/www/materials/modules/chapter23.pdf

[13]  Tolstov G. P., Silverman R. A., Fourier Series, Courier Corporation, 1976

[14]  Al-Gwaiz M.A., Sturm-Liouville Theory and its Applications, Springer Undergraduate Mathematics Series ISSN 1615-2085, Springer-Verlag London Limited 2008